\documentclass{article}

\usepackage{amsmath}
\usepackage{amsfonts}
\usepackage{amssymb}
\usepackage{amscd}
\usepackage{amsthm}

\def\sp{[\![}
\def\spp{]\!]}
\def\dbar{\;\;\bar{}\!\!\!d}

\newtheorem{lemma}{Lemma}[section]
\newtheorem{corollary}[lemma]{Corollary}
\newtheorem{theorem}[lemma]{Theorem}
\newtheorem{proposition}[lemma]{Proposition}
\newtheorem{definition}[lemma]{Definition}
\newtheorem{remark}[lemma]{Remark}
\newtheorem{hypothesis}[lemma]{Hypothesis}
\newcommand{\eh}{\hfill}
\newlength{\sperrT}

\newenvironment{proofT}{{\settowidth{\sperrT}{\rm Proof of Theorem 3.1.}
\par\addvspace{0.3cm}\noindent\parbox[t]{1.3\sperrT}{\rm P\eh r\eh o\eh o\eh 
f\eh \ \eh o\eh f\eh \ \eh T\eh h\eh e\eh o\eh r\eh e\eh m\eh \ \eh 3\eh .\eh 1\eh 
.\eh }}}{\nopagebreak\mbox{}\hfill $\blacksquare $\par\addvspace{0.25cm}}
\newlength{\sperrTT}

\newenvironment{proofTT}{{\settowidth{\sperrTT}{\rm Proof of Theorem 3.1.}
\par\addvspace{0.3cm}\noindent\parbox[t]{1.3\sperrTT}{\rm P\eh r\eh o\eh o\eh 
f\eh \ \eh o\eh f\eh \ \eh T\eh h\eh e\eh o\eh r\eh e\eh m\eh \ \eh 7\eh .\eh 3\eh 
.\eh }}}{\nopagebreak\mbox{}\hfill $\blacksquare $\par\addvspace{0.25cm}}

\numberwithin{equation}{section}

\begin{document}

\title{Magnetic Pseudodifferential Operators}

\date{\today}

\author{Viorel Iftimie, Marius M\u antoiu and Radu Purice\footnote{Institute
of Mathematics ``Simion Stoilow'' of
the Romanian Academy, P.O.  Box 1-764, Bucharest, RO-70700, Romania,
Email: viftimie@math.math.unibuc.ro, mantoiu@imar.ro, purice@imar.ro}}

\maketitle

\begin{abstract}
In previous papers, a generalization of the Weyl calculus was introduced in connection with the quantization of a particle moving in $\mathbb R^n$ under the influence of a variable magnetic field $B$. It incorporates phase factors defined by $B$ and reproduces the usual Weyl calculus for $B=0$. In the present article we develop the classical pseudodifferential theory of this formalism for the standard symbol classes $S^m_{\rho,\delta}$. Among others, we obtain properties and asymptotic developments for the magnetic symbol multiplication, existence of parametrices, boundedness and positivity results, properties of the magnetic Sobolev spaces. In the case when the vector potential $A$ has all the derivatives of order $\ge 1$ bounded, we show that the resolvent and the fractional powers of an elliptic magnetic pseudodifferential operator are also pseudodifferential. As an application, we get a limiting absorption principle and detailed spectral results for self-adjoint operators of the form $H=h(Q,\Pi^A)$, where $h$ is an elliptic symbol, $\Pi^A=D-A$ and $A$ is the vector potential corresponding to a short-range magnetic field. 
\end{abstract}



\section*{Introduction}

There are several different, but related, points of view on the usual Weyl calculus (for which we refer to \cite{Ho1}, \cite{Ho2}, \cite{Fo}, \cite{Sh}). One of them says that the correspondence symbol $\mapsto$ operator, $f\mapsto \mathfrak {Op}(f)$, is a functional calculus for the family of operators $Q_1,\dots,Q_n;D_1,\dots,D_n$ on $L^2(\mathbb R^n)$, where $Q_j$ is the multiplication with the variable $x_j$ and $D_j=-i\partial_j$. The familiar notation $\mathfrak {Op}(f)=f(Q,D)$ keeps track  of this fact. The relative sophistication of this formalism has its roots in the non-commutativity of the basic operators:

\begin{equation}\label{CCR}
i[Q_j,Q_k]=0=i[D_j,D_k],\ \ \ \ i[D_j,Q_k]=\delta_{j,k}.
\end{equation}
In particular, this is the reason why the symbol multiplication $(f,g)\mapsto f\circ g$, fitted to fulfil $f(Q,D)g(Q,D)=(f\circ g)(Q,D)$, is more complicated than pointwise multiplication. 

This remark suggests that generalizations of the Weyl calculus could be motivated by considering more general commutation relations as the starting point of a functional calculus. For the case of a constant magnetic field such a calculus has been developped in \cite{BMGH}, or in \cite{B} for the case of a lattice.

For a nonrelativistic quantum particle in $\mathbb R^n$ placed in a magnetic field $B$ deriving from a vector potential $A$, the basic self-adjoint operators (quantum observables) are the positions $Q_1,\dots,Q_n$ and the magnetic momenta $\Pi^A_1:=D_1-A_1,\dots,\Pi^A_n:=D_1-A_n$, satisfying the commutation relations

\begin{equation}\label{commag}
i[Q_j,Q_k]=0,\ \ \ \ i[\Pi^A_j,Q_k]=\delta_{jk},\ \ \ \ i[\Pi^A_j,\Pi^A_k]=B_{jk},
\end{equation}
where $B_{jk}$ is (the operator of multiplication by) the component $(jk)$ of the magnetic field. These relations are much more complicated than those of Heisenberg (\ref{CCR}). Actually, they are a representation by unbounded operators of a Lie algebra that has infinite dimension if $B_{jk}$ are not all polynomial functions. 

It is natural to look for a pseudodifferential calculus adapted to such a situation. At first sight, a procedure could be to replace in the explicit formula for $\mathfrak {Op}(f)$ the symbol $f(x,\xi)$ by $f(x,\xi-A(x))$, obtaining an operator $\mathfrak {Op}_A(f)$.  Although largely used in the literature (see \cite{GMS}, \cite{Ic1}, \cite{Ic2}, \cite{II}, \cite{IT1}, \cite{IT2}, \cite{ITs1}, \cite{ITs2}, \cite{NU1}, \cite{NU2}, \cite{Pa}, \cite{Um}), this point of view does not seem adequate, due to the fact that the operators $\mathfrak {Op}_A(f)$, although representing physical observables, are not gauge covariant. Two vector potentials $A$ and $A'$ connected by $A'=A+d\varphi$ for some smooth real function $\varphi$, being assigned to the same magnetic field $B=dA=dA'$, should produce unitarily equivalent operators $\mathfrak {Op}_A(f)$ and $\mathfrak {Op}_{A'}(f)$ for all reasonable $f$. In section 6 we are going to exibit large classes of symbols $f$ (including the third order monomial $f(x,\xi)=\xi_j\xi_k\xi_l$) for which the expected equality $\mathfrak {Op}_{A+d\varphi}(f)=e^{i\varphi}\mathfrak {Op}_A(f)e^{-i\varphi}$ fails. 

The right formalism was proposed independently and with different emphases in \cite{KO1}, \cite{KO2} and \cite{MP1}, \cite{MP2}. It was generalized and related to a $C^*-$algebraic formalism in \cite{MPR1}, and applied to the strict deformation quantization in the sense of Rieffel for systems in a magnetic field in \cite{MP3} and \cite {MP4}. The $C^*$-algebraic setting is also related to the canonical commutation relations (\ref{commag}), which can be reformulated by saying that the group $\mathbb R^n$ admits an action on itself that is twisted by a $2$-cocycle defined by the flux of the magnetic field. The $C^*$-algebras canonically assigned to this twisted dynamical system are an essentially isomorphic version of the magnetic pseudodifferential calculus.

We shall sketch very briefly only the pseudodifferential point of view in Section 1, while the other sections will be dedicated to our actual purposes: a development of the calculus for H\"ormander symbol classes $S^m_{\rho,\delta}$ and applications. For the moment, let us only state that the core of the theory consists of two formulae: 

\noindent
{\it 1. The quantization rule.} The input is a magnetic field in $\mathbb R^n$, given by a smooth $2$-form $B$. This $2$-form is closed, by one of the Maxwell equations. Thus we can write it as $B=dA$ for some highly non-unique $1$-form $A$, called vector potential. For a vector potential $A$ and points $x,y$ in $\mathbb R^n$, one can calculate the circulation $\int_{[x,y]}A$. The operators associated to some suitable ``symbol classes'' will be of the form:
$\quad\left[\mathfrak{Op}^A(f)u\right](x)=$
\begin{equation}\label{OpA}
=(2\pi)^{-n}\int\limits_{\mathbb R^n}\int\limits_{\mathbb R^n}\,dy\,d\xi\, e^{i<x-y,\xi>}\exp\left(-i\int\limits_{[x,y]}A\right)f\left(\frac{x+y}{2},\xi\right)u(y)
\end{equation}
for any $u\in\mathcal{S}(\mathbb{R}^n)$. The imaginary exponential of the circulation of the vector potential enters (\ref{OpA}) as a modification of the well-known Weyl formula for $\mathfrak{Op}(f)$. One should think of $\mathfrak{Op}^A(f)$ as the quantum observable (ideally a self-adjoint operator) corresponding to the classical observable given by the (real, smooth) function $f$ defined on the phase space $\mathbb R^{2n}$. In \cite{MP2} we indicated a sort of derivation of this formula, as done in \cite{Ho1}, \cite{Ho2} or \cite{Fo} for $B=0$. The point is to get first (\ref{OpA}) for exponential functions and then to use superposition for the general case. But quantizing exponential functions requires a calculation of the exponential of a linear combination of the operators $Q_j$ and $\Pi^A_j$, magnetic translations replace usual translations in the output, and this is the source of the extra factor $\exp\left(-i\int_{[x,y]}A\right)$.

\noindent
{\it 2. The composition law.} When no magnetic field is present, the classical mechanical stage for our particle is the phase space $\mathbb R^{2n}=\mathbb R^{n}\times\mathbb R^{n}$, which is a symplectic vector space with the canonical symplectic form $\sigma[(y,\eta),(z,\zeta)]=\sum_{j=1}^n(z_j\eta_j-y_j\zeta_j)$. When $B$ is turned on, the same $\mathbb R^{2n}$ should be seen as a symplectic manifold with the perturbed symplectic form

\begin{equation}\label{sigmaB}
\sigma_{B,(x,\xi)}[(y,\eta),(z,\zeta)]=\sum_{j=1}^n(z_j\eta_j-y_j\zeta_j)+\sum_{j,k}B_{jk}(x)y_jz_k.
\end{equation}
(The $2$-form $B$ on $\mathbb R^n$ can be pulled-backed to a $2$-form on $\mathbb R^{2n}$ - seen as the cotangent bundle of $\mathbb R^n$.) We refer to  \cite{DR}, \cite{MaR} and \cite{MP4} for details. If one belives that (\ref{OpA}) is justified and wants to have $\mathfrak{Op}^A(f)\mathfrak{Op}^A(g)=\mathfrak{Op}^A(f\circ^B g)$, then the right answer is $\quad\left(f\circ^B g\right)(X)=$
\begin{equation}\label{circB}
=\pi^{-2n}\int\limits_{\mathbb R^{2n}}\int\limits_{\mathbb R^{2n}}\,dY\,dZ\,\exp\left(-i\int\limits_{\mathcal T(X,Y,Z)}\sigma_B\right)f(X-Y)g(X-Z),
\end{equation}
with $\int_{\mathcal T(X,Y,Z)}\sigma_B$ the integral of the $2$-form $\sigma_B$ through the plane triangle $\mathcal T(X,Y,Z)$ in $\mathbb R^{2n}$ defined by the corners $X-Y-Z,\ X+Y-Z,\ X+Z-Y$. Actually (\ref{circB}), which we call {\it the magnetic Moyal product}, reduces for $B=0$ to the usual multiplication law in Weyl calculus, often called {\it the Moyal product}. For further discussions concerning this phase function, especially of a geometrical nature, see \cite{KO1} and \cite{KO2}. 
For the same relations containing also Planck's constant $\hbar$ and a study of the semiclassical limit $\hbar\rightarrow 0$ in a $C^*$-algebraic formalism we refer to \cite{MP3} and \cite{MP4}; here $\hbar=1$.

The formal deformation quantization (in the sense of star products) of cotangent bundles endowed with a magnetic symplectic form can be found in \cite{BNPW}.

The nice feature is that (\ref{circB}) is completely intrinsic. It makes use directly of the magnetic field $B$; no choice of a vector potential is needed. One may say that introducing $B$ leads to $B$-dependent function algebras, and only when we want linear representations of these abstract algebras vector potentials are needed. It is already obvious that these representations behave coherently with respect to equivalent choices. Vector potentials $A$ and $A'$ related by $A'=A+d\varphi$ give the same magnetic field $B$. Then the operators $\mathfrak{Op}^A(f)$ and $\mathfrak{Op}^{A'}(f)$, belonging to different representations, are unitarily equivalent: $\mathfrak{Op}^{A'}(f)=e^{i\varphi}\mathfrak{Op}^A(f)e^{-i\varphi}$.

The pseudodifferential calculus with a magnetic field has been used in several papers dealing with the Peierls substitution (\cite{DS}, \cite{HS1}, \cite{HS2}, \cite{N}, \cite{PST}, \cite{T}). Although gauge covariance is not essential for the technical arguments used in this context, it is possible that our formalism may bring some new insight and even technical advantages.

In the previously quoted articles \cite{MP2}, \cite{MP3}, \cite{MP4} and \cite{MPR1}, many technical difficulties were avoided by sometimes restricting to small classes of symbols $f,g$. This was enough for what was aimed and allowed very general magnetic fields. When general symbols were considered (by extending the validity of the equations by duality techniques), very often refined properties of the resulting objects are hidden. We are now at a point where these difficulties must be faced and the powerful graded framework of pseudodifferential theory should emerge. The presence of the two extra phase factors in (\ref{OpA}) and (\ref{circB}) makes all the picture rather complicated. For the time being, {\it we have only succeeded to treat smooth, bounded magnetic fields having bounded derivatives of all orders}. Since, however, no decay at infinity is requested, we belive that the theory we develop is general enough to support nontrivial applications.

Having indicated roughly the motivations of our topic, let us now describe the content of the article.

First we state some notations and conventions.

In a first section we sum up briefly, with few justification, several simple facts from \cite{MP2} and \cite{MPR1}. For the convenience of the reader we show that (\ref{OpA}) and (\ref{circB}) are tied together, by proving that $\mathfrak{Op}^A(f)\,\mathfrak{Op}^A(g)=\mathfrak{Op}^A(f\circ^Bg)$.

In Section 2 we show first that for $f\in S^m_{\rho,\delta}$, ($\rho\geq0$, $\delta<1$), $\mathfrak{Op}^A(f)$ leaves the Schwartz space invariant. Afterwards, we study the product $f\circ^B g$ for $f,g$ belonging to H\"ormander's classes of  symbols $S^m_{\rho,\delta}$, ($0\leq\delta<\rho\leq1$). This is basic for the rest of the article. We give an asymptotic series for $f\circ^B g$, the terms being calculated by recurrence. The development of the commutator $f\circ^B g-g\circ^B f$ starts with the Poisson bracket $\{f,g\}_B$ assigned to the magnetic symplectic form $\sigma_B$. The existence of a parametrix for elliptic magnetic pseudodifferential operators is also proved, as a consequence of the asymptotic development of $f\circ^B g$; this is one of the most important tools of the theory.

In Section 3 we prove that $\mathfrak{Op}^A(f)$ is bounded in $L^2(\mathbb{R}^n)$ if $f\in S^0_{\rho,\delta}$ and $0\leq\delta\leq\rho\leq1$, $\delta<1$. The case $\delta=\rho$ is a magnetic version of the Calderon-Vaillancourt theorem transcripted for the Weyl calculus (cf. \cite{Fo}).

Section 4 is dedicated to the study of magnetic Sobolev spaces. Previously, they were considered only in situations when a vector potential can be chosen with bounded derivatives of strictly positive order, cf. \cite{GMS} and \cite{Pa}.

In Section 5 we show that an elliptic magnetic pseudodifferential operator is self-adjoint on the corresponding Sobolev spaces. For convenient vector potentials, the Schwartz space is a core. As a consequence of a G\aa rding-type inequality, we also treat semiboundedness.

Throughout the paper we suppose the magnetic field $B$ to have bounded components together with all their derivatives. In Section 6 we shall moreover assume that $B=dA$ for some smooth vector potential $A$ having bounded derivatives of any strictly positive order. This facilitates certain arguments; in particular it leads to a connection between our magnetic calculus and the Weyl calculus for a certain $A$-dependent H\"ormander-type metric, and this allows the transcription of certain classical results (\cite{Ho1}, \cite{Ho2}, \cite{Bo3}) to our framework.

As said before, many authors use for a symbol $p(x,\xi)$ the magnetic quantization $\mathfrak{Op}_A(p):=\mathfrak{Op}(p_A)$ with $p_A(x,\xi):=p(x,\xi-A(x))$, that does not provide a gauge covariant calculus. In Section 6, as a continuation of the analysis in Subsection IV D of \cite{MP2}, we shall compare this procedure with our gauge covariant quantization and prove that $\mathfrak{Op}^A(p)-\mathfrak{Op}_A(p)$ is a pseudodifferential operator of strictly smaller order. 
In fact we prove a little bit more, showing that for any symbol $p$, one can find a symbol $q$ of the same order (the difference $p-q$ having a strictly inferior order) such that $\mathfrak{Op}^A(p)=\mathfrak{Op}_A(q)$. Thus, under the above hypothesis on the magnetic field $B$, one can  pass from the functional calculus $\mathfrak{Op}^A$ to the functional calculus $\mathfrak{Op}_A$ and in the opposite direction. Using the Weyl-H\"{o}rmander-Bony calculus leads to more precise results; for example we obtain a Fefferman-Phong type theorem and prove that the resolvent and the powers of a magnetic self-adjoint elliptic pseudo-differential operator are also pseudodifferential. In particular we are able to compare three candidates for the relativistic Schr\"{o}dinger Hamiltonian with magnetic field: $\sqrt{(D-A)^2+1}$, $\mathfrak{Op}_A(<\xi>)$ and $\mathfrak{Op}^A(<\xi>)$.

The last section is devoted to the spectral analysis (obtaining a limiting absorption principle) for a class of elliptic pseudodifferential operators obtained through a quantization (either by Weyl calculus, or by $\mathfrak{Op}^A$, or $\mathfrak{Op}_A$), for a symbol of the form $p=p_0+p_S+p_L$, where $p_0$ does not depend on $x$, $p_S$ is a symbol with "short range" behaviour and $p_L$ is a symbol with "long range" behaviour. We assume that all the derivatives of the magnetic field $B$ verify conditions of type "short range" at infinity, and our exemple 3 shows that these hypothesis are in some sense optimal. The spectral analysis of $\mathfrak{Op}_A(<\xi>)$ has been done in \cite{Um} but without considering the problem of a limiting absorption principle; moreover, as shown in exemple 2, our hypothesis are more general.

The main ingredient in proving our Theorem \ref{89} is an abstract result from the conjugate operator theory (see \cite{ABG}); verifying the necessary hypothesis involves the elaborate pseudodifferential calculus developed in the present article.

\section*{Notations and conventions}

We denote $\mathbb R^n$ by $\mathcal X$, with elements $x,y,z$. $\mathcal X^*$ will be the dual of $\mathcal X$, with elements $\xi,\eta,\zeta$. We also denote by $<\cdot,\cdot>$ the duality form: $<\xi,x>=\xi(x)=<x,\xi>$.

The phase space will be $\Xi=\mathbb R^{2n}=\mathcal X\oplus\mathcal X^*$, with elements $X=(x,\xi),\,Y=(y,\eta),\,Z=(z,\zeta)$. In fact these notations will be used in a rigid manner: if the contrary is not explicitly stated, when one encounters $X\in \Xi$, one should think that its components in $\mathcal X$ and $\mathcal X^*$, respectively, are called $x$ and $\xi$; the same for $Y$ and $Z$. $\Xi$ is a symplectic space with the canonical symplectic form $\sp Y, Z\spp=\sp (y,\eta), (z,\zeta)\spp=<\eta,z>-<\zeta,y>$.

On $\mathcal X$ we consider the usual Lebesgue measure. But on $\mathcal X^*$ and $\Xi$, respectively, it will be convenient to use $\dbar \xi=(2\pi)^{-n}d\xi$ and $\dbar X=\pi^{-n}dX$. 

If $\mathcal Y$ is one of the spaces $\mathcal X$, $\mathcal X^*$ or $\Xi$ (an $\mathbb R^m$ essentially), we set $C^\infty_0(\mathcal Y)=\{f\in C^\infty(\mathcal Y)\mid \text{supp}f\ \text{is compact}\}$, with $C^\infty(\mathcal Y)$ the space of infinitely derivable complex functions on $\mathcal Y$.

We use standard multi-index notations: $\alpha=(\alpha_1,\dots,\alpha_m)\in\mathbb N^m$, $\vert\alpha\vert=\alpha_1+\dots+\alpha_m$, $\alpha !=\alpha_1 !\dots\alpha_m !$, $\partial_y^{\alpha}=\partial_{y_1}^{\alpha_1}\dots\partial_{y_m}^{\alpha_m}$ or $\partial^{\alpha}=\partial_1^{\alpha_1}\dots\partial_m^{\alpha_m}$, 
$D^\alpha=i^{-|\alpha|}\partial^\alpha$, where $m=\dim\mathcal{Y}$.

$\mathcal S(\mathcal Y)$ will be the Schwartz space on $\mathcal Y$, with antidual $\mathcal S^*(\mathcal Y)$; we denote by $(u,v):=u(v)$, ($u\in\mathcal{S}^*(\mathcal{Y}),\;v\in\mathcal{S}(\mathcal{Y})$) the application of anti-duality. $L^p(\mathcal Y)$ are the standard Lebesgue spaces for $p\in[1,\infty]$. We set $BC^\infty(\mathcal Y)=\{f\in C^\infty(\mathcal Y)\mid \partial^\alpha f\ \text{is bounded for any\ }\alpha\in\mathbb N^m\}$. $C^\infty_{\text{pol}}(\mathcal Y)$ is the space of all $C^\infty$ functions on $\mathcal Y$ with the absolute value of each derivative dominated by an (arbitrary) polynomial. $C^\infty_{\text{pol,u}}(\mathcal Y)=\{f\in C^\infty(\mathcal Y)\mid\exists k\ge 0\ \text{such that}\ \vert(\partial^\alpha f)(y)\vert\le<y>^k\ \forall y\in\mathcal Y\;\alpha\in\mathbb N^m\}$ is the subspace of $C^\infty_{\text{pol}}(\mathcal Y)$ consisting of elements whose all derivatives are dominated by a polynomial of fixed (arbitrary) degree. 

We frequently consider integrals as converging in $\mathcal S^*$, in particular as oscillatory integrals.

We denote by $\mathcal B(\mathcal H_1,\mathcal H_2)$ the Banach space of all linear bounded operators $T:\mathcal H_1\rightarrow\mathcal H_2$, with $\mathcal H_1$, $\mathcal H_2$ Hilbert (or Banach) spaces. In fact we preserve this notation even if $\mathcal H_1$, $\mathcal H_2$ are topological vector spaces, to signify continuous, linear operators. For $\mathcal B(\mathcal H,\mathcal H)$ we abbreviate $\mathcal B(\mathcal H)$. $\mathcal K(\mathcal H_1,\mathcal H_2)$ will denote compact operators. $\mathcal B_p(\mathcal H)$ will be the Schatten-von Neumann class of order $p\in[1,\infty]$ on $\mathcal H$. For $p=1$ we have trace-class operators, for $p=2$ we get Hilbert-Schmidt operators, and $\mathcal B_\infty(\mathcal H)=\mathcal K(\mathcal H)$.

Given a Riemannian metric $g_{_X}$ on $\Xi$ and a positive function $M:\Xi\rightarrow\mathbb{R}_+^*$, we define the {\it symbol space}
$S(M,g)$ to be the space of $C^\infty$ functions $f:\Xi\rightarrow\mathbb{C}$ such that
$$
\underset{g_{_X}(T_j)\leq1}{\underset{X,T_j\in\Xi}{\sup}}
\left(M(X)^{-1}\left|\partial_{T_1}\ldots\partial_{T_k}f(X)\right|\right)
<\infty,\;\forall k\in\mathbb{N},
$$
where we denote by $\partial_Tf$ the derivative of $f$ with respect to the direction $T\in\Xi$. 
We denote by $\boldsymbol{\Psi}(M,g)$ the family of Weyl operators $\mathfrak{Op}(f)$ with $f\in S(M,g)$.

If $M(X)=<\xi>^m$ and the metric has the form
$$
g_{_X}=<\xi>^{2\delta}|dx|^2\,+\,<\xi>^{-2\rho}|d\xi|^2
$$
for $\rho$, $\delta$ and $m$ real numbers and $<\xi>:=(1+|\xi|^2)^{1/2}$, we denote $S(M,g)$ by $S^m_{\rho,\delta}(\Xi)$.  We stil use the notations
$$
S^m(\Xi):=S^m_{1,0}(\Xi),\qquad S^{-\infty}(\Xi)=\underset{m\in\mathbb{R}}{\bigcap}S^m_{\rho,\delta}(\Xi).
$$
Explicitly, a function $f\in C^\infty(\Xi)$ belongs to $S^m_{\rho,\delta}(\Xi)$ if for any multi-indices $\alpha$ and $\beta$ in $\mathbb{N}^n$ there exists a finite constant $C_{\alpha\beta}$ such that
$$
\left|\left(\partial^\alpha_x\partial^\beta_\xi f\right)(X)\right|\,\leq\,C_{\alpha\beta}\,<\xi>^{m-\rho|\beta|+\delta|\alpha|},
\quad\forall X=(x,\xi)\in\Xi.
$$

\section{Preliminaries}  

\subsection{The magnetic field}

The mathematical framework that we consider is supposed to model a quantum particle without internal structure moving in $\mathcal X=\mathbb{R}^n$, in the presence of a non-uniform magnetic field. The {\it magnetic field} is described by a closed 2-form $B$ on $\mathcal X\equiv\mathbb{R}^n$. In the standard coordinate system on $\mathbb{R}^n$, it is represented by a function taking real antisymmetric matrix values $B=(B_{jk})$ with $1\leq j\leq n$, $1\leq k\leq n$ and verifying the relation $\partial_jB_{kl}+\partial_kB_{lj}+\partial_lB_{jk}=0$. We shall always assume that $B_{jk}\in C^\infty_{\text{pol}}(\mathcal X)$, although this is not necessary for all constructions or assertions. Anyhow, later on, even stronger assumptions on $B$ will be imposed.

Any such field $B$ may be written as the exterior differential $dA$ of a 1-form $A$, {\it the vector potential}; by using coordinates, one has $B_{jk}=\partial_jA_k-\partial_kA_j$ for each $j,k=1,\cdots,N$. The components of the vector potential will always be taken of class $C^\infty_{\text{pol}}(\mathcal X)$, in order to define multipliers for $\mathcal{S}(\mathcal{X})$ and $\mathcal{S}^*(\mathcal{X})$. This is, indeed, always possible, as can be seen by considering {\it the transversal gauge}

\begin{equation}\label{transversal}
A_j(x)=-\sum_{k=1}^n\int_0^1ds\;B_{jk}(sx)sx_k.
\end{equation}

In the magnetic pseudodifferential calculus that we shall develop there are two phase factors that play an important role, one defined by $B$ and the other by $A$. Given a $k$-form $C$ on $X$ and a compact piecewise smooth $k$-surface $\gamma\subset X$, we denote by
$$
\Gamma^C(\gamma):=\int_\gamma C
$$
the usual invariant integral. We shall encounter circulations of the 1-form $A$ along linear segments $\gamma=[x,y]$ defined by points (ends) $x,y$ and fluxes of the 2-form $B$ through triangles $\gamma=<x,y,z>$ defined by points (corners) $x,y,z$. By Stokes' theorem, one has 

\begin{equation}\label{St-1}
\Gamma^B(<x,y,z>)=\Gamma^A([x,y])+\Gamma^A([y,z])+\Gamma^A([z,x]).
\end{equation}

\subsection{The magnetic functional calculus.}

In a former paper \cite{MP2} we have shown that for $f\in\mathcal{S}^*(\Xi)$ and $u\in\mathcal{S}(\mathcal X)$, the formula (properly interpreted)
\begin{equation}\label{OpAA}
\left[\mathfrak{Op}^A(f)u\right](x)=\int\!\!\!\int_{\mathcal X\times\mathcal{X}^*}dy\,\dbar\eta\, e^{i<x-y,\eta>}e^{-i\Gamma^A([x,y])}f\left(\frac{x+y}{2},\eta\right)u(y)
\end{equation}
defines an integral operator $\mathfrak{Op}^A(f)\in\mathcal{B}(\mathcal{S}(\mathcal{X}),
\mathcal{S}^*(\mathcal{X}))$, and in fact $\mathfrak{Op}^A$ gives an isomorphism between $\mathcal{S}^*(\Xi)$ and $\mathcal{B}(\mathcal{S}(\mathcal{X}),\mathcal{S}^*(\mathcal{X}))$ (as linear topological spaces) that restricts to an isomorphism between $\mathcal{S}(\Xi)$ and $\mathcal{B}(\mathcal{S}^*(\mathcal{X}),\mathcal{S}(\mathcal{X}))$.
Let us remark here that for any test functions $u$ and $v$ in $\mathcal{S}(\mathcal{X})$ and any distribution $f\in\mathcal{S}^*(\mathcal{X})$, we have the relation
$$
\left(\mathfrak{Op}^A(f)u\,,\,v\right)\,=\,\left(u\,,\,\mathfrak{Op}^A
(\overline{f})v\right),
$$
where $(\overline{f},u):=\overline{(f,\overline{u})}$. In particular, if $f$ is a real distribution (i.e. $\overline{f}=f$), then $\mathfrak{Op}^A(f)$ is a symetric operator in $\mathcal{B}(\mathcal{S}(\mathcal{X}),\mathcal{S}^*(\mathcal{X}))$.

An important property is {\it gauge covariance}.
Let $A$ and $A'$ be two vector potentials of class $C_{\text{pol}}^\infty$, defining the same magnetic field, $dA=B=dA'$. Then there exists a real function $\varphi\in C_{\text{pol}}^\infty(X)$ such that $A'=A+\nabla\varphi$ and $e^{i\varphi(Q)}\mathfrak{Op}^A(f)e^{-i\varphi(Q)}=\mathfrak{Op}^{A+\nabla\varphi}(f)\;$ for any $f\in\mathcal S'(\Xi)$ and all such functions $\varphi$; this second identity is valid in $\mathcal{B}\left[\mathcal{S}(\mathcal X),\mathcal{S}^*(\mathcal X)\right]$.

In \cite{MP2} we show that $\mathfrak{Op}^A$ induces a unitary map from $L^2(\Xi)$ to $\mathcal{B}_2(L^2(\mathcal{X}))$, the ideal of all Hilbert-Schmidt operators. The family of operators $\mathfrak{Op} ^A(f)$, $f$ being the Fourier tranform of an arbitrary function in $L^1(\Xi)$, is dense in the closed ideal $\mathcal{K}(L^2(\mathcal X))$ of all compact operators.

\subsection{The magnetic Moyal product}

Let $\tilde{f}$ be the Fourier transformation of $f\in\mathcal{S}^*(\Xi)$ with respect to the second variable. If $\Lambda^A(x,y):=\exp\{-i\Gamma^A([x,y])\}$, we can write
$$
\left[\mathfrak{Op}^A(f)u\right](x)=(2\pi)^{-n}\int_{\mathcal X}dy\,\Lambda^A(x,y)\,\tilde{f}\left(\frac{x+y}{2},y-x\right)\,u(y).
$$

For $f,g\in\mathcal{S}(\Xi)$, the associated magnetic Weyl operators $\mathfrak{Op}^A(f)$ and $\mathfrak{Op}^A(g)$ are smoothing operators and thus the product $\mathfrak{Op}^A(f)\mathfrak{Op}^A(g)$ is also smoothing, consequently of the form $\mathfrak{Op}^A(f\circ^B g)$ with $f\circ^B g$ in $\mathcal{S}(\Xi)$ depending linearly on $f$ and $g$. In order to obtain an explicit form for the magnetic Moyal product $f\circ^B g$, we remark that for $\phi\in\mathcal{S}(\mathcal{X})$:
$$
\left[\mathfrak{Op}^A(f)\,\mathfrak{Op}^A(g)\,\varphi\right](x)\;
=(2\pi)^{-n}\;\int_{\mathcal X}\,dz\,\Lambda^A(x,z)\times
$$
$$
\times\left\{\int_{\mathcal X}\,dy\;\Omega^B(x,y,z)\,\tilde{f}\left(\frac{x+y}{2},y-x\right)\,\tilde{g}\left(\frac{y+z}{2},z-y\right)\right\}\varphi(z),
$$
where $\Omega^B(x,y,z):=\exp\{-i\Gamma^B(<x,y,z>)\}$ and we have used (\ref{St-1}).

Thus we have obtained $(\widetilde{f\circ^B g})\left(\frac{x+z}{2},z-x\right)=$
$$
=(2\pi)^{-n}\;\int_{\mathcal X}\,dy\;\Omega^B(x,y,z)\,\tilde{f}\left(\frac{x+y}{2},y-x\right)\,\tilde{g}\left(\frac{y+z}{2},z-y\right).
$$
Let us remark that all the integrals above are absolutely convergent due to our assumptions on $f$ and $g$. Moreover we have seen that $f\circ^B g\in\mathcal{S}(\Xi)$, so that we can compute its partial Fourier transform (in the second variable) by the usual integral formula and obtain again an element in $\mathcal{S}(\Xi)$. Thus, after making the change of variables $u:=(x+z)/2$, $v:=z-x$, we can write (for any $\theta\in\mathcal{X}^*$)
$$
(f\circ^B g)(u,\theta)\;=(2\pi)^{-n}\;\int_{\mathcal X}dv\;e^{i<\theta,v>}(\widetilde{f\circ^B g})(u,v)
$$
$$
=(2\pi)^{-2n}\;\int_{\mathcal X}\int_{\mathcal X}\int_{\mathcal X^*}\int_{\mathcal X^*}dv\,dy\,d\eta\,d\zeta\;e^{i<\theta,v>}\;\Omega^B\left(u-(v/2),y,u+(v/2)\right)\,\times
$$
$$
\exp\{-i<\eta,y-u+(v/2)>\}\,f\left(\frac{u+y-(v/2)}{2},\eta\right)\times
$$
$$
\times\;\exp\{-i<\zeta,u+(v/2)-y>\}\,g\left(\frac{u+y+(v/2)}{2},\zeta\right).
$$
We shall use the Fubini Theorem and the change of variables
$
2y^\prime:=u-y+(v/2)$, $2z^\prime:=u-y-(v/2)$,  $\eta^\prime:=\theta-\eta$, $\zeta^\prime:=\theta-\zeta$.
In the sequel we shall use the notation
$$
\omega_B(x,y,z):=\exp\{-4iF_B(x,y,z)\}=\Omega^B(x-y+z,x-y-z,x+y-z),
$$
where
$$
F_B(x,y,z):=\frac{1}{4}\Gamma^B(<x-y+z,x-y-z,x+y-z>)=
$$
$$
=\sum_{j,k=1}^ny_j(z_k-y_k)\int_0^1\int_0^1ds\,dt\,\,
sB_{jk}(x-y-z+2sy+2st(z-y)).
$$
By easy computations we get {\it the Moyal product} $\quad (f\circ^B g)(X)\;=$
\begin{equation}\label{circB'}
=\;\int\limits_\Xi \int\limits_\Xi \dbar Y\dbar Z\; e^{-2i\sp Y,Z\spp}\,\omega^B(x,y,z)\,f(X-Y)\,g(X-Z)\;=
\end{equation}
$$
=\;\int\limits_\Xi \int\limits_\Xi  \dbar Y\dbar Z \; e^{-2i\sp X-Y,X-Z\spp}\,\Omega^B(x-z+y,-x+y+z,x-y+z)\,f(Y)\,g(Z).
$$

\subsection{The magnetic Moyal $^*$-algebras}

One can extend the validity of (\ref{circB'}) by duality, using the fact that for any functions $f$ and $g$ in $\mathcal{S}(\Xi)$ we have
$$
\int_\Xi dX\;(f\circ^B g)(X)=\int_\Xi dX\;(g\circ^B f)(X)=\int_\Xi dX\;f(X)g(X)=<f,g>\equiv (f,\overline g).
$$
Thus for $f,g,h\in\mathcal{S}(\Xi)$ we have
$
<f\circ^B g,h>=<f,g\circ^B h>= <g,h\circ^B f>.
$
Considering $<\cdot,\cdot>$ as duality between $\mathcal{S}^\prime(\Xi)$ and $\mathcal{S}(\Xi)$, we define for $F\in\mathcal{S}^\prime(\Xi)$ and $f\in\mathcal{S}(\Xi)$
$$
<F\circ^B f,h>:=<F,f\circ^B h>,\qquad <f\circ^B F,h>:=<F,h\circ^B f>,\qquad \forall h\in\mathcal{S}(\Xi),$$
getting two bilinear continuous mappings $\mathcal{S}^\prime(\Xi)\times\mathcal{S}(\Xi)\rightarrow\mathcal{S}^\prime(\Xi)$, resp. $\mathcal{S}(\Xi)\times\mathcal{S}^\prime(\Xi)\rightarrow\mathcal{S}^\prime(\Xi)$, with good associativity properties.

A substantial extension of the magnetic Moyal product is obtained in \cite{MP2} on the following class of distributions

\begin{equation}\nonumber
\mathcal{M}^B(\Xi):=\left\{F\in\mathcal{S}^\prime(\Xi)\;\mid\;F\circ^B f\in\mathcal{S}(\Xi),\ f\circ^B F\in\mathcal{S}(\Xi),\;\;\forall f\in\mathcal{S}(\Xi)\right\},
\end{equation}
called {\it the magnetic Moyal algebra}. For any $F,G\in\mathcal{M}^B(\Xi)$, we define

$$
<F\circ^B G,h>:=<F,G\circ^B h>,\qquad \forall h\in\mathcal{S}(\Xi).
$$

The set $\mathcal{M}^B(\Xi)$ together with the composition law $\circ^B$ defined above and the complex conjugation $F\mapsto \overline F$ is an unital $^*$-algebra, containing $\mathcal S(\Xi)$ as a self-adjoint two-sided ideal.
Maybe the most important fact is that $\mathfrak{Op}^A$ is an isomorphism of $\ ^*$-algebras betweeen $\mathcal{M}^B(\Xi)$ and $\mathcal B[\mathcal S(\mathcal X)]\cap\mathcal B[\mathcal S'(\mathcal X)].$ 

Simple examples show that $\mathcal{M}^B(\Xi)$ is much larger than $\mathcal S(\Xi)$. For instance, Fourier transforms of bounded, complex measures on $\Xi$, as well as $C^\infty_{\text{pol}}$-functions depending only on the variable in $\mathcal X$ are in the magnetic Moyal algebra. A less evident and very useful fact is that $C^\infty_{\text{\rm pol,u}}(\Xi)\subset\mathcal{M}^B(\Xi)$. Lemma \ref{43} will show that $S^m_{\rho,\delta}(\Xi)$ is also contained in $\mathcal{M}^B(\Xi)$.

\section{The magnetic composition of symbols}

Our first task is to extend the magnetic composition law to classes of symbols and obtain a precise asymptotic development for this composition, generalizing the well known formulae from usual pseudodifferential calculus. This will be a key technical ingredient in the following developments, in particular leading to the existence of parametrices for elliptic operators.

\subsection{H\"ormander symbols are in the magnetic Moyal algebra}
\noindent
To make use of the natural extension to the Moyal algebra (as discussed in [1.4]), we begin by shawing  that the classical symbol spaces $S^m_{\rho,\delta}(\Xi)$ are contained in the magnetic Moyal algebra. Only the case $\delta\leq0$ is covered by our previous result in \cite{MP2}, saying that $C^\infty_{\text{\sl pol, u}}(\Xi)\subset\mathcal M^B(\Xi)$.
\begin{lemma}\label{43}
If $B$ is a magnetic field having components of class $C^\infty_{\text{\rm
pol}}(\mathcal{X})$, then for any $m\in\mathbb{R}$, any $\rho\geq 0$ and any $\delta<1$ we
have $S^m_{\rho,\delta}(\Xi)\subset\mathcal{M}^B(\Xi)$.
\end{lemma}
\begin{proof}
We must prove that, for any couple $(f,\varphi)\in
S^m_{\rho,\delta}(\Xi)\times\mathcal{S}(\Xi)$, we have $f\circ^B\!\varphi\in\mathcal{S}(\Xi)$. 
Thus we have to study the following integral, that will exist as an oscillatory integral:
$$
(f\circ^B\!\varphi)( X)=\int_{\Xi}\int_{\Xi}
\dbar Y\,\dbar Z\,\,\exp\{-2i\sp Y, Z\spp\}
\,\omega_B(x,y,z)\,f( X- Y)
\varphi( X- Z).
$$

We choose $\chi\in C^\infty_0(\Xi)$ with $\chi(0)=1$, and for any $\epsilon>0$ we define $f_\epsilon( X):=\chi(\epsilon X)f( X)$.
Then we shall show that the limit $\underset{\epsilon\rightarrow
0}{\lim}(f_\epsilon\circ^B\!\varphi)$ exists pointwise and is
independent of the choice of $\chi$. 
By integration by parts we have
$$
(f_\epsilon\circ^B\!\varphi)( X)=
\int_{\Xi}\int_{\Xi}
\dbar Y\,\dbar Z\,e^{-2i\sp Y, Z\spp}
<y>^{-2N_\zeta}<\eta>^{-2N_z}\,f_{\epsilon}( X- Y)\times
$$
\begin{equation}\label{a1}
\times
\mathfrak{L}_z^{N_z}
\left[\omega^B(x,y,z)\mathfrak{L}_{\zeta}^{N_\zeta}
\varphi( X- Z)\right],
\end{equation}
with differential operators $\mathfrak{L}_z=1-(1/4)\Delta_z$ and $\mathfrak{L}_{\zeta}=1-(1/4)\Delta_\zeta$.
The integrals are well defined, due to the decay assumptions on $\varphi$. We choose first
$N_z\geq(1/2)(m+n+1)$ and then $N_\zeta\geq(1/2)(q(2N_z)+n+1)$, where
$q(N):=\underset{|\gamma|\leq N}{\max}p(\gamma)$ and $p(\gamma)$ is defined by
the following estimations (implied by the assumptions on $B$):
$$
\left|\partial^\gamma_{x,y,z}\omega_B(x,y,z)\}\right|\leq
C_\gamma(1+|x|+|y|+|z|)^{p(\gamma)}.
$$
We conclude that we can take the limit $\epsilon\rightarrow 0$ and obtain for 
$(f\circ^B\!\varphi)( X)$ an identity similar to (\ref{a1}).
Moreover, this equation is clearly independent on the choice of
$\chi$ and also of the exact choices of $N_z$ and $N_\zeta$ (by integration by
parts).
For any $k\in\mathbb{N}$ we may choose $N_z\geq(1/2)(m+k(|\delta|+|\rho|)+n+1)$
and $N_\zeta\geq(1/2)(q(2N_z+k)+n+1)$ in order to prove (by further integration
by parts with respect to $y$ and $\eta$) that in fact
$f\circ^B\!\varphi\in C^k(\Xi)$. In conclusion we proved that 
$f\circ^B\!\varphi\in C^\infty(\Xi)$.

Moreover, if we consider a multiindex
$\alpha=(\alpha_x,\alpha_\xi)\in\mathbb{N}^{2n}$ and integrate by parts with
respect to $y$ and $\eta$, we prove that 
$\partial_{ X}^\alpha(f\circ^B\!\varphi)$ is a finite linear combination
of terms of the form
$$
I( X)=\int_\Xi\int_\Xi\,\dbar Y\,\dbar Z\,e^{-2i\sp
 Y, Z\spp}<z>^{-2N_\eta}<\zeta>^{-2N_y}\times
$$
$$
\times
\left(\partial^{\beta^\prime}_\eta<\eta>^{-2N_z}\right)
\left(\partial^{\gamma^\prime}_y<y>^{-2N_\zeta}\right)
\left\{\partial^{\gamma^{\prime\prime}}_y\partial^{\delta^\prime}_z
\partial^{\alpha_x^\prime}_x\omega^B(x,y,z)\right\}\times
$$
$$
\times
\left\{\partial^{\alpha_x^{\prime\prime}}_x\partial^{\alpha_\xi^\prime}_\xi
\partial^{\delta^{\prime\prime}}_z\partial^{\lambda}_\zeta
\varphi( X- Z)
\right\}
\left\{\partial^{\alpha_x^{\prime\prime\prime}}_x
\partial^{\alpha_\xi^{\prime\prime}}_\xi
\partial^{\gamma^{\prime\prime\prime}}_y\partial^{\beta^{\prime\prime}}_\eta
\,f( X- Y)\right\},
$$
with
$\alpha_x^\prime+\alpha_x^{\prime\prime}+\alpha_x^{\prime\prime\prime}=\alpha_x$, 
$\alpha_\xi^\prime+\alpha_\xi^{\prime\prime}=\alpha_\xi$, $|\alpha|=k$,
$|\beta^\prime|+|\beta^{\prime\prime}|\leq2N_\eta$, $|\gamma^\prime|+
|\gamma^{\prime\prime}|+|\gamma^{\prime\prime\prime}|\leq2N_y$, 
$|\delta^\prime|+|\delta^{\prime\prime}|\leq2N_z$ and $|\lambda|\leq2N_\zeta$. 
Then we can choose $N_y,N_\eta,N_z,N_\zeta$ as functions of $k$ and $l$ so that
$|I( X)|\leq C_l< X>^{-l}$. 
In fact, taking into account that
$\varphi\in\mathcal{S}(\Xi)$, we have for any $(s,t)\in\mathbb{N}\times\mathbb{N}$
the inequalities
$$
\left|\left(\partial^{\gamma^{\prime\prime}}_y\partial^{\delta^\prime}_z
\partial^{\alpha_x^\prime}_x\omega^B\right)(x,y,z) \right|\leq C_1(<x>+<y>+<z>)^{q(2N_y+2N_z+k)},
$$
$$
\left|\left(\partial^{\alpha_x^{\prime\prime}}_x\partial^{\alpha_\xi^\prime}_\xi
\partial^{\delta^{\prime\prime}}_z\partial^{\lambda}_\zeta
\varphi\right)( X- Z)
 \right|\leq C_2<x-z>^{-t}<\xi-\zeta>^{-s},
$$
$$
\left| \left(\partial^{\alpha_x^{\prime\prime\prime}}_x
\partial^{\alpha_\xi^{\prime\prime}}_\xi
\partial^{\gamma^{\prime\prime\prime}}_y\partial^{\beta^{\prime\prime}}_\eta
\,f\right)( X- Y)\right|\leq C_3
<\xi-\eta>^{m-\rho|\alpha_\xi^{\prime\prime}+\beta^{\prime\prime}|+
\delta|\alpha_x^{\prime\prime\prime}+\gamma^{\prime\prime\prime}|}.
$$
We may suppose that $\delta\in(0,1)$. Let us remark that
$$
<\xi-\eta>^{m-\rho|\alpha_\xi^{\prime\prime}+\beta^{\prime\prime}|+
\delta|\alpha_x^{\prime\prime\prime}+\gamma^{\prime\prime\prime}|}\;\leq
C<\xi>^{m+\delta(k+2N_y)}
<\eta>^{|m|+\delta(k+2N_y)},
$$
$$
<x-z>^{-t}\;\leq\;C<x>^{-t}<z>^{t},\quad <\xi-\zeta>^{-s}\;
\leq\;C<\xi>^{-s}<\zeta>^{s}.
$$
We choose first $N_y\geq(1/2)(1-\delta)^{-1}(n+1+l+m+k\delta)$ in order to
verify the integrability condition with respect to $\zeta\in \mathcal{X}$. Then we choose $s=l+m+\delta(k+2N_y)$, and obtain a factor $<\xi>^{-l}$. We can also choose $N_z\geq(1/2)(n+1+|m|+(k+2N_y)\delta)$ in order to verify the
integrability condition with respect to $\eta\in \mathcal{X}$ and $t=l+q(2N_y+2N_z+k)$ to obtain a factor $<x>^{-l}$. We end up by choosing $N_\eta\geq(1/2)(n+1+l+2q(2N_y+2N_z+k))$ and $N_\zeta\geq(1/2)(n+1+q(2N_y+2N_z+k))$ in order to get the
integrability with respect to $(y,z)\in \mathcal{X}\times \mathcal{X}$.
\end{proof}

\subsection{Estimations on the magnetic flux}

If the magnetic field $B$ has components of class
$BC^\infty(\mathcal{X})$, by arguments similar to those above, for any $f\in S^{m_1}(\Xi)$ and $g\in
S^{m_2}(\Xi)$,
the magnetic Moyal product $f\circ^B\! g$ belongs to $S^{m_1+m_2}(\Xi)$. Some sharp estimations on the flux of the magnetic field will make possible in [2.4] a precise result concerning the asymptotic development of $f\circ^B\! g$.

\begin{lemma}\label{45}
If the components of the magnetic field $B$ are of class $C^\infty(\mathcal{X})$,
then
$
\;\partial_{x_j}F_B=\sum_{k=1}^n\left(D_{jk}y_k+E_{jk}z_k\right)
$,
$
\;\partial_{y_j}F_B=\sum_{k=1}^n\left(D_{jk}^\prime y_k+E_{jk}^\prime z_k\right)
$,
$
\;\partial_{z_j}F_B=\sum_{k=1}^n\left(D_{jk}^{\prime\prime}y_k+
E_{jk}^{\prime\prime}z_k\right)
$,
where the coefficients $D_{jk}$, \ldots , $E_{jk}^{\prime\prime}$ are of class $BC^\infty(\mathcal{X}^3)$. In particular
$$
D_{jk}^\prime(x,y,z):=\frac{1}{4}\int_{-1}^1tdt\,\left[B_{jk}(x+ty-z)+
B_{jk}(x-ty+tz)\right],
$$
$$
E_{jk}^\prime(x,y,z):=\frac{1}{4}\int_{-1}^1dt\,\left[tB_{jk}(x-ty+tz)-
B_{jk}(x-y+tz)\right],
$$
$$
D_{jk}^{\prime\prime}(x,y,z):=-\frac{1}{4}\int_{-1}^1dt\,
\left[tB_{jk}(x-ty+tz)+B_{jk}(x+ty-z)\right],
$$
$$
E_{jk}^{\prime\prime}(x,y,z):=\frac{1}{4}\int_{-1}^1tdt\,
\left[B_{jk}(x-ty+tz)-B_{jk}(x-y+tz)\right].
$$
\end{lemma}
\begin{proof} 
Using the condition $dB=0$, i.e. $\partial_lB_{jk}+\partial_jB_{kl}+\partial_kB_{lj}=0$, we obtain
$$
\partial_{x_l}F_B(x,y,z)=
$$
$$
=-\frac{1}{2}\sum_{k=1}^n(z_k-y_k)\int_0^1ds\,
\int_0^sdt\,\frac{d}{ds}B_{kl}(x-y-z+2sy+2t(z-y))-
$$
$$
-\frac{1}{2}\sum_{j=1}^ny_j\int_0^1ds\,
\int_0^sdt\,\frac{d}{dt}B_{lj}(x-y-z+2sy+2t(z-y)).
$$
The identity
$
\int_0^1ds\,\int_0^sdt\,f(s,t)=\int_0^1dt\,
\int_t^1ds\,f(s,t)
$
and some straightforward computations lead to the first equality in the statement of the Lemma. The others are then proved in a similar way.
\end{proof}

The next Corollary follows from Lemma \ref{45} by direct calculation. A crucial fact in our arguments is the boundedness in $x$ of all the derivatives of $\omega_B$.

\begin{corollary}\label{46}
If we assume that all the components of the magnetic field $B$ are of class
$BC^\infty(\mathcal{X})$, then we have
$$
\left|\left(\partial^\alpha_x\partial^\beta_y\partial^\gamma_z\omega_B\right)(x,y,z)\right|\leq C_{\alpha,\beta,\gamma}(<y>+<z>)^{|\alpha|+|\beta|+|\gamma|},\;\;
\forall(\alpha,\beta,\gamma)\in[\mathbb{N}^n]^3,
$$
where $C_{\alpha,\beta,\gamma}$ are positive constants.
\end{corollary}
\begin{remark}\label{46-r}
The same proof as above shows that if the components of the magnetic field are of class $BC^\infty(\mathcal{X})$, then
$$
\left|\partial^\alpha_z \left[e^{-i\Gamma^B(<x,y,z>)}\right] \right|\,\leq\,C_\alpha(<x-z>+<y-z>)^{|\alpha|},\quad\forall\alpha
\in\mathbb{N}^n,
$$
where $C_\alpha$ are positive constants.
\end{remark}

\subsection{Stationary phase result}

An important ingredient for the estimation of the integral appearing in (\ref{circB'}) is a 'stationary phase' result, for which we introduce some notations.
For any $\varphi\in C^\infty(\mathcal{X}^3)$ we define the following first order differential operator (with respect to the variables $U=(u,\mu)\in\Xi$ and $V=(v,\nu)\in\Xi$), having coefficients that only depend on $(x,y,z)$:
$$
M_B(\varphi):=[\omega_B(x,y,z)]^{-1}\sum_{j=1}^n\left[
\partial_{y_j}\left(\omega_B\varphi\right)\partial_{\nu_j}-
\partial_{z_j}\left(\omega_B\varphi\right)\partial_{\mu_j}
\right].
$$
We shall denote
$$
\sp\partial_{ U},\partial_{ V}\spp:=<\partial_\mu,\partial_v>
-<\partial_u,\partial_\nu>\equiv\sum_{j=1}^n\left(\partial_{\mu_j}
\partial_{v_j}-\partial_{u_j}\partial_{\nu_j}\right)
$$
and for $t\in\mathbb{R}\setminus\{0\}$ and $\varphi\in C^\infty_{\text{\rm pol}}(\mathcal{X}^3)$ we define:
$$
L_{\varphi}(t):=\frac{1}{2i}\left\{2\varphi\sp
\partial_{ U},\partial_{ V}\spp+t^{-1}M_0(\varphi)\right\}.
$$
\begin{lemma}\label{48}
Let $\varphi\in C^\infty_{\text{\rm pol}}(\mathcal{X}^3)$. 
\begin{enumerate}
\item We have the following equality:
$$
\int_{\Xi}\int_{\Xi}\dbar Y\,\dbar Z\,\,\exp\{-2i\sp Y, Z\spp\}\varphi(x,y,z)=\varphi(x,0,0)
$$
(the integral being interpreted as an oscillatory integral).
\item If $h\in\mathcal{S}(\Xi\times\Xi)$, for any $t\in\mathbb{R}^*$ we have
$$
\int_{\Xi}\int_{\Xi}\,\dbar Y\,\dbar Z\,\,
e^{-2i\sp Y, Z\spp}\varphi(x,y,z)h( X-t Y,
 X-t Z)\,=\,\varphi(x,0,0)h( X, X)\,+
$$
$$
+t^2\int_0^1sds\int_{\Xi}\int_{\Xi}\,\dbar Y\,\dbar Z\,\,
e^{-2i\sp Y, Z\spp}\left[L_\varphi(st)h\right]( X-
st Y, X-st Z).
$$
\end{enumerate} 
\end{lemma}
\begin{proof}
(1) Let us fix $\chi\in C^\infty_0(\mathcal{X})$ such that $\chi(0)=1$. For any $\epsilon>0$ we define
$$
I_\epsilon:=\int_{\Xi}\int_{\Xi}
\dbar Y\,\dbar Z\,\,
\chi(\epsilon y)\chi(\epsilon \eta)\chi(\epsilon z)\chi(\epsilon \zeta)
\exp\{-2i\sp Y,
 Z\spp\}\varphi(x,y,z)=
$$
$$
=\pi^{-2n}\int_{\mathcal{X}}\int_{\mathcal{X}}dy\,dz\,\,
\chi(\epsilon y)\chi(\epsilon z)\varphi(x,y,z)\epsilon^{-2n}\hat{\chi}(-(2/\epsilon)y)
\hat{\chi}((2/\epsilon)z)=
$$
$$
=(2\pi)^{-2n}\int_{\mathcal{X}}\int_{\mathcal{X}}dy\,dz\,\,
\chi(-(\epsilon^2/2) y)\chi((\epsilon^2/2) z)\hat{\chi}(y)\hat{\chi}(z)
\varphi(x,-(\epsilon/2)y,(\epsilon/2)z).
$$
By the Lebesgue dominated convergence theorem, this converges for $\epsilon\rightarrow 0$ to $\varphi(x,0,0)$.

(2) We perform a first order Taylor expansion of $h( X-t Y, X-t Z)$ with respect to $t$. The first term on the right-hand side is given by 1. For the second term we integrate by parts, using the identities
$$
y_je^{-2i\sp Y, Z\spp}=\frac{1}{2i}\partial_{\zeta_j}
e^{-2i\sp Y, Z\spp}\;,\qquad
\eta_je^{-2i\sp Y, Z\spp}=-\frac{1}{2i}\partial_{z_j}
e^{-2i\sp Y, Z\spp},
$$
$$
z_je^{-2i\sp Y, Z\spp}=-\frac{1}{2i}\partial_{\eta_j}
e^{-2i\sp Y, Z\spp}\;,\qquad
\zeta_je^{-2i\sp Y, Z\spp}=\frac{1}{2i}\partial_{y_j}
e^{-2i\sp Y, Z\spp}.
$$
\end{proof}

\subsection{Magnetic composition of symbols}

To any differential operator of order $m$ with respect to the variables $ U$ and $ V$, $P:=\sum\;c_{\alpha\beta}(x,y,z)\partial^\alpha_{ U}
\partial^\beta_{ V},$
we associate another differential operator $M_B(P)$ defined by
$
M_B(P):=\sum\;M_B(c_{\alpha\beta})\partial^\alpha_{ U}
\partial^\beta_{ V}
$.
This operator will evidently have the same form, but will be of order $m+1$.

For any sequence of positive numbers $\{t_j\in\mathbb{R}^*\,\mid\,j\in\mathbb{N}^*\}$, let $t_{(k)}:=t_1\cdot\ldots\cdot t_k$ (for $k\in\mathbb{N}^*$) and
let us define by recurrence the sequence of differential operators:
$$
L_0:=\boldsymbol{1},
$$
$$
L_1(t_1):=\omega_B^{-1}L_{\omega^B}(t_1),
$$
$$
L_{j+1}(t_1,\ldots,t_{j+1}):=L_1(t_{(j+1)})L_j(t_1,\ldots,t_j)
+\frac{M_0\left(L_j(t_1,\ldots, t_j)\right)}{2it_{(j+1)}}.
$$

\begin{theorem}\label{50}
Let us assume that all the components of the magnetic field $B$ are of class $BC^\infty(\mathcal{X})$. If $f\in S^{m_1}_{\rho,\delta}(\Xi)$ and $g\in S^{m_2}_{\rho,\delta}(\Xi)$, with $m_1\in\mathbb{R}$, $m_2\in\mathbb{R}$, $0\leq\delta < \rho\leq 1$, then $f\circ^B\!g\in S^{m_1+m_2}_{\rho,\delta}(\Xi)$ and we have the following asymptotic development:
$$
f\circ^B\!g\;\sim \; \sum_{j=0}^\infty h_j,\qquad h_j\in S^{m_1+m_2-j(\rho-\delta)}_
{\rho,\delta}(\Xi), \qquad h_0(X)=f(X)g(X),
$$
$$
h_j( X)=\int_0^1\int_0^1\cdots
\int_0^1dt_1dt_2\ldots dt_j\;t_1^{2j-1}t_2^{2j-3}\cdot\ldots\cdot t_j\times
$$
$$
\times[L_j(t_1,\ldots,t_j)(f( U)g( V))]_{_{}}\big|
\underset{y=z=0}{\underset{U= V= X}{}}.
$$
\end{theorem}

\begin{proof}
We shall verify by induction that for any $k\geq 1$ we have
\begin{equation}\label{adMp}
f\circ^B\!g=\sum_{j=0}^{k-1}h_j+R_k,
\end{equation}
where $\quad R_k( X)=$
$$
=\int_0^1\int_0^1\cdots
\int_0^1dt_1dt_2\ldots dt_k\;t_1^{2k-1}t_2^{2k-3}\cdot\ldots\cdot t_k
\int_{\Xi}\int_{\Xi}\,\dbar Y\,
\dbar Z\,\,e^{-2i\sp Y, Z\spp}\times
$$
$$
\times\omega_B(x,y,z)
\left[L_k(t_1,\ldots,t_k)\left(f( U)\otimes g( V)
\right)\right]\big|
\underset{V= X-t_{(k)} Z}{\underset{U= X-t_{(k)} Y}{}}.
$$

For $k=1$ we apply Lemma \ref{48}, taking $\varphi=\omega_B$, $h=f\otimes g$ and $t=1$, $s=t_1$. We just have to remark that in this case
$$
L_{\omega_{_B}}(t_1)=\frac{1}{2i}\left\{2\omega_B\sp\partial_{ U},
\partial_{ V}\spp+\omega_Bt_1^{-1}M_B(1)\right\}=\omega_BL_1(t_1).
$$
Regarding the integrals as oscillatory integrals, we may assume that the functions $f$ and $g$ belong to $\mathcal{S}(\Xi)$.

Let us suppose now that formula (\ref{adMp}) is verified for some $k\geq 1$. In order to prove it for $k+1$ we shall rewrite the integral defining the rest $R_k$. We notice that
$$
[L_k(t_1,\ldots,t_k)(f\otimes g)](U,V)=\sum_{\alpha,\beta}a_{_{\alpha\beta}}(x,y,z,t_1,\ldots,t_k)
\left(\partial^\alpha f\right)( U)\left(\partial^\beta g\right)(V).
$$
For each term we apply Lemma \ref{48}, taking $\varphi=\omega_Ba_{_{\alpha\beta}}$, $h=(\partial^\alpha f)\otimes
(\partial^\beta g)$ and $t=t_1\cdots t_k$, $s=t_{k+1}$. 
We remark first that
$$
\sum_{\alpha,\beta}a_{_{\alpha\beta}}(x,0,0,t_1,\ldots,t_k)
\left(\partial^\alpha f\right)( X)\left(\partial^\beta g\right)( X)=
$$
$$
=\left[L_k(t_1,\ldots,t_k)(f( U)g( V))
\right]\big|
\underset{y=z=0}{\underset{U= V= X}{}}.
$$
Then
$$
L_{(\omega_{_B}a_{_{\alpha\beta}})}(t_{(k+1)})=\frac{1}{2i}
\left\{2a_{_{\alpha\beta}}\omega_B\sp\partial_{ U},\partial_{ V}\spp + (t_{(k+1)})^{-1}M_0(a_{_{\alpha\beta}}\omega_B)\right\}=
$$
$$
=\frac{1}{2i}
\omega_B\left\{2a_{_{\alpha\beta}}\sp\partial_{ U},\partial_{ V}\spp + (t_{(k+1)})^{-1}M_B(a_{_{\alpha\beta}})\right\},
$$
and moreover
$$
M_B(a_{_{\alpha\beta}})=\omega_B^{-1}\sum_{j=1}^n\left[\partial_{y_j}
(\omega_Ba_{_{\alpha\beta}})\partial_{\nu_j}-\partial_{z_j}
(\omega_Ba_{_{\alpha\beta}})\partial_{\mu_j}\right]=
a_{_{\alpha\beta}}M_B(1)+M_0(a_{_{\alpha\beta}}).
$$
Putting all these together we get
$$
L_{(\omega_Ba_{_{\alpha\beta}})}(t_{(k+1)})=
\omega_B\left\{a_{_{\alpha\beta}}
L_1(t_{(k+1)})+\frac{1}{2i}(t_{(k+1)})^{-1}M_0(a_{_{\alpha\beta}})
\right\}.
$$
We remark further that $L_1(t_{(k+1)})$ is a differential operator with respect to the variables $ U$ and $ V$ and thus commutes with multiplication with the function $a_{_{\alpha\beta}}$; thus
$
\sum M_0(a_{_{\alpha\beta}})
\partial^\alpha_{ U}\partial^\beta_{ V}=
M_0(L_k(t_1,\ldots,t_k))
$.
Finally we get
$
R_k( X)
=h_k( X)+R_{k+1}( X)
$
and the proof of (\ref{adMp}) is finished.

Let us show now that, for each $j\in\mathbb N$, we have $h_j\in S^{m_1+m_2-j(\rho-\delta)}_{\rho,\delta}(\Xi)$. This is evident for $j=0$. For $j\geq 1$ notice that
$$
L_j(t_1,\ldots,t_j)=\sum a_{_{\alpha^\prime,\alpha^{\prime\prime},\beta^\prime,\beta^{\prime\prime}}}
(x,y,z;t_1,\ldots,t_j)\partial_u^{\alpha^\prime}\partial_\mu^{\alpha^{\prime\prime}}
\partial_v^{\beta^\prime}\partial_\nu^{\beta^{\prime\prime}},
$$
where $|\alpha^{\prime\prime}|+|\beta^{\prime\prime}| = j$, $|\alpha^\prime|+
|\beta^\prime|\leq j$ and $a_{_{\alpha^\prime,\alpha^{\prime\prime},\beta^\prime,\beta^{\prime\prime}}}$ are linear combinations of products of derivatives of $F_B$ with respect to $y$ and $z$ and monomials in $t_1^{-1}$, \ldots , $t_j^{-1}$, with exponents that do not exceed those of $t_1$, \ldots, $t_j$ appearing in the integrals. Lemma \ref{45} shows that for $y=z=0$ the derivatives of $F_B$ are either vanishing or at least bounded functions of $x$. We conclude that $h_j$ is a linear combination of terms of the type
$$
b(x)\left[\left(\partial_x^{\alpha^\prime}\partial_\xi^{\alpha^{\prime\prime}}f
\right)
( X)\right]\left[\left(\partial_x^{\beta^\prime}
\partial_\xi^{\beta^{\prime\prime}}g\right)( X)\right]
\int_0^1t_1^{p_1}\,dt_1\cdots\int_0^1t_j^{p_j}\,dt_j,
$$
where $b\in BC^\infty(\mathcal{X})$, $p_l\geq 0$ for any $l\in\{1,\ldots,j\}$ and $|\alpha^{\prime\prime}|+|\beta^{\prime\prime}| = j$, $|\alpha^\prime|+
|\beta^\prime|\leq j$. It follows immediatly, from the hypothesis on $f$ and $g$, that $h_j\in S^{m_1+m_2-j(\rho-\delta)}_{\rho,\delta}(\Xi)$.

We shall end the proof of the Theorem by checking that for any $k\geq 1$ one has
$
R_k\in S^{m_1+m_2+2n-k(\rho-\delta)}_{\rho,\delta}(\Xi).
$
We shall use once again the structure of the operator $L_j(t_1,\ldots,t_j)$ that was described above. It follows that $R_k$ can be written as a linear combination of terms of the form
$$
I( X)=C\int_0^1\cdots\int_0^1dt_1\ldots dt_k\;
t_1^{p_1}\cdots t_k^{p_k}\int_{\Xi}\int_{\Xi}
\dbar Y\dbar Ze^{-2i\sp Y, Z\spp}
\omega_B(x,y,z)\times
$$
$$
\times
a_{_{\alpha^\prime,\alpha^{\prime\prime},\beta^\prime,\beta^{\prime\prime}}}
(x,y,z;t_1,\ldots,t_k)\times
$$
$$
\times\left[\left(\partial_x^{\alpha^\prime}\partial_\xi^{\alpha^{\prime\prime}}f
\right)( X-t_1\cdots t_k Y)\right]
\left[\left(\partial_x^{\beta^\prime}
\partial_\xi^{\beta^{\prime\prime}}g\right)
( X-t_1\cdots t_k Z)\right].
$$
Now we no longer restrict to $y=z=0$ and thus factors of the type $y^\beta z^\gamma$ may appear in the functions $a_{_{\alpha^\prime,\alpha^{\prime\prime},\beta^\prime,\beta^{\prime\prime}}}$ (that contain derivatives of $F_B$). These factors may be handled by integration by parts, using the exponential $e^{-2i\sp Y, Z\spp}$, and will generate operators of the form $\partial_\eta^\gamma\partial_\zeta^\beta$ applied to the functions $\partial_x^{\alpha^\prime}\partial_\xi^{\alpha^{\prime\prime}}f$ and $\partial_x^{\beta^\prime}\partial_\xi^{\beta^{\prime\prime}}g$, but this will not alter their decay. We proceed as in the proof of Lemma \ref{43}, and after a number of integrations by parts we write $I( X)$ as a linear combination of terms of the type
$$
J( X)=\int_0^1\cdots\int_0^1dt_1\cdots dt_k\;
t_1^{q_1}\cdots t_k^{q_k}\int_{\Xi}\int_{\Xi}
\dbar Y\dbar Ze^{-2i\sp Y, Z\spp}\times
$$
$$
\times
<z>^{-2N_\eta}<\zeta>^{-2N_y}\left(\partial_\eta^{\gamma^\prime}<\eta>^{-2N_z}
\right)\left(\partial_y^{\delta^\prime}<y>^{-2N_\zeta}\right)\times
$$
$$
\times\partial_y^{\delta^{\prime\prime}}\partial_z^{\epsilon^{\prime}}\left(
\omega^Ba_{_{\alpha^\prime,\alpha^{\prime\prime},\beta^\prime,\beta^{\prime\prime}}}\right)(x,y,z)
\times
$$
$$
\times
\left(\partial_x^{\alpha^\prime+\delta^{\prime\prime\prime}}
\partial_\xi^{\alpha^{\prime\prime}+\gamma^{\prime\prime}}f\right)( X-t_{(k)} Y)
\left(\partial_x^{\beta^\prime+\epsilon^{\prime\prime}}\partial_\xi^
{\beta^{\prime\prime}+\lambda}g\right)( X-t_{(k)} Z),
$$
where $|\gamma^\prime|+|\gamma^{\prime\prime}|\leq2N_\eta$, $|\delta^\prime|+|\delta^{\prime\prime}|+|\delta^{\prime\prime\prime}|\leq2N_y$, $|\epsilon^{\prime}|+|\epsilon^{\prime\prime}|\leq2N_z$, $|\lambda|\leq2N_\zeta$ and $q_j\geq 0$ for any $j\in\{1,\ldots,k\}$.

We fix $N_\eta=N_\zeta=n$ in order to have integrability in the variables $y$ and $z$. Then we decompose the $\eta$-integral with respect to the two domains $\{|\eta|\leq\kappa<\xi>\}$ and $\{|\eta|\geq\kappa<\xi>\}$ for some small fixed $\kappa>0$, and similarly for the $\zeta$-integration, and thus write $J( X)$ as a sum of 4 terms $J_a( X)$ (with $a=1,2,3,4$). In order to estimate each of these 4 terms separately we remark that we may choose different values for the pair $(N_y,N_z)$ in each of these terms, due to the fact that these choices are made by integration by parts in the variables $y$ and $z$. Moreover, for any $r\in\mathbb{R}$ and for $\forall N_y\geq0$, $\forall N_z\geq0$
$$
\int_{\{|\eta|\leq\kappa<\xi>\}}d\eta
<\eta>^{-2N_z}<\xi-t\eta>^{r+2\delta N_y}\; \leq\; C<\xi>^{r+2\delta N_y+n}.
$$
Taking $N_y\geq0$ and $2N_z>|r|+2\delta N_y+n$, we have
$$
\int_{\{|\eta|\geq\kappa<\xi>\}}d\eta
<\eta>^{-2N_z}<\xi-t\eta>^{r+2\delta N_y}\; \leq
$$
$$
\leq\; C\int_{\{|\eta|\geq\kappa<\xi>\}}d\eta\,<\eta>^{|r|+2\delta N_y-2N_z}
\;\leq\;C<\xi>^{|r|+2\delta N_y-2N_z+n}.
$$
Let us set $r_1:=m_1-\rho|\alpha^{\prime\prime}|+\delta|\alpha^\prime|$ and $r_2:=m_2-\rho|\beta^{\prime\prime}|+\delta|\beta^\prime|$. We get the following upper bound for $|J_1( X)|$
$$
\underset{0\leq t\leq1}{\sup}
\int_{\{|\eta|\leq\kappa<\xi>\}}d\eta
\frac{<\eta>^{2N_z}}{
<\xi-t\eta>^{r_1+2\delta N_y}}
\int_{\{|\zeta|\leq\kappa<\xi>\}}d\zeta\frac{<\zeta>^{2N_y}}{
<\xi-t\zeta>^{r_2+2\delta N_z}}\;\leq
$$
$$
\leq\;C<\xi>^{m_1+m_2-k(\rho-\delta)+2n},
$$
by choosing for this domain $N_y=N_z=0$. Then, we get for the next type of domain the upper bound:
$$
\underset{0\leq t\leq1}{\sup}
\int_{\{|\eta|\leq\kappa<\xi>\}}d\eta\frac{<\eta>^{2N_z}}{
<\xi-t\eta>^{r_1+2\delta N_y}}
\int_{\{|\zeta|\geq\kappa<\xi>\}}d\zeta\frac{<\zeta>^{2N_y}}{
<\xi-t\zeta>^{r_2+2\delta N_z}}\;\leq
$$
$$
\leq\;C<\xi>^{r_1+2\delta N_y+n}<\xi>^{|r_2|-2N_y+n}\;=\;C<\xi>^
{r_1+|r_2|+2n-2(1-\delta)N_y},
$$
and we have to choose on this domain $N_z=0$ and $N_y$ large enough. On the similar domain with $\eta$ and $\zeta$ interchanged we have to choose $N_y=0$ and $N_z$ large enough. For the fourth domain we obtain the upper bound
$$
\underset{0\leq t\leq1}{\sup}
\int_{\{|\eta|\geq\kappa<\xi>\}}d\eta\frac{<\eta>^{2N_z}}{
<\xi-t\eta>^{r_1+2\delta N_y}}
\int_{\{|\zeta|\geq\kappa<\xi>\}}d\zeta\frac{<\zeta>^{2N_y}}{
<\xi-t\zeta>^{r_2+2\delta N_z}}\;\leq
$$
$$
\leq\;C<\xi>^{|r_1|+|r_2|+2n-2(1-\delta)(N_y+N_z)}
$$
and thus we have to choose both $N_y$ and $N_z$ large enough. We conclude that
$$
|R_k( X)|\;\leq\;C<\xi>^{m_1+m_2+2n-k(\rho-\delta)}.
$$
For the derivatives of $R_k$ we may proceed in a similar way, since all the terms obtained by differentiating the factor $\omega_B$ with respect to $x$ are bounded by monomials in $y$ and $z$ and thus can be dealt with by integration by parts in $y$ and $z$. We get
$$
\left|\left(\partial_x^\alpha\partial_\xi^\beta R_k\right)( X)\right|
\;\leq\;C<\xi>^{m_1+m_2+2n-k(\rho-\delta)-\rho|\beta|+\delta|\alpha|},
\;\;\forall\alpha\in\mathbb{N}^n,\;\forall\beta\in\mathbb{N}^n.
$$

In order to end our proof we fix some $p\in\mathbb{N}$ and write the identity
$$
f\circ^B\!g\;-\;\sum_{j=0}^ph_j\;=\;\sum_{j=p+1}^kh_j\;+\;R_k,
$$
where $k\geq p+1$ is chosen large enough to have $2n-k(\rho-\delta)\leq-(p+1)(\rho-\delta)$. This gives
$$
f\circ^B\!g\;-\;\sum_{j=0}^ph_j\;\in S^{m_1+m_2-(p+1)(\rho-\delta)}(\Xi),
$$
and thus $f\circ^B\!g\sim\sum_jh_j$.
\end{proof}

\subsection{The first terms of the asymptotic development}

Let us explicitly compute the first terms of the asymptotic development of the
magnetic Moyal product.

\noindent
1. We know from the statement of Theorem \ref{50} that
$
h_0=fg.
$

\noindent
2. In order to compute $h_1$, we remark that
$$
L_1(t_1)(f\otimes g)=
-i\sum_{j=1}^n(\partial_{\mu_j}f)\otimes(\partial_{v_j}g)-
(\partial_{u_j}f)\otimes(\partial_{\nu_j}g)-
$$
$$
-2t_1^{-1}\sum_{j=1}^n(\partial_{y_j}F_B)[f\otimes(\partial_{v_j}g)]-
(\partial_{z_j}F_B)[(\partial_{\mu_j}f)\otimes g].
$$
For $y=z=0$ the first order derivatives of $F_B$ vanish and thus
$$
\left.[L_1(t_1)(f\otimes g)]\frac{}{}\right|_{_{y=z=0}}\hspace{-0.8cm}( X, X)
=-i\sum_{j=1}^n\left[(\partial_{\xi_j}f)( X)(\partial_{x_j}g)
( X)-(\partial_{x_j}f)( X)(\partial_{\xi_j}g)
( X)\right].
$$
In conclusion we have
$
h_1=-\frac{i}{2}\{f,g\}.
$

\noindent
3. For $h_2$ we need to compute the explicit form of the operator
$
L_2(t_1,t_2)$. Using Lemma \ref{45} we obtain:
$$
\left.\left[M_0(L_1(t_1))(f\otimes g)\right]\frac{}{}\right|_{_{y=z=0}}\hspace{-0.8cm}( X, X)=
-\frac{2}{t_1}\sum_{j,k=1}^nB_{jk}(x)[(\partial_{\xi_j}f)\otimes(\partial_{\xi_k}g)]
( X, X).
$$
$$
L_1(t_1\cdot t_2)L_1(t_1)(f\otimes g)=
$$
$$
=
-\sum_{j,k=1}^n\left(\partial_{\mu_j}\partial_{\mu_k}\partial_{v_j}\partial_{v_k}+
\partial_{u_j}\partial_{u_k}\partial_{\nu_j}\partial_{\nu_k}-2
\partial_{\mu_j}\partial_{u_k}\partial_{v_j}\partial_{\nu_k}\right)(f\otimes g)\;+
$$
$$
+\;\{\text{\sl terms vanishing for y=z=0}\}.
$$
Finally, we put everything together to obtain
$$
h_2( X)=\frac{1}{8}\sum_{j,k=1}^n\left[
2\partial_{\xi_j}\partial_{x_k}f
\otimes\partial_{x_j}\partial_{\xi_k}g-
\partial_{\xi_j}\partial_{\xi_k}f\otimes\partial_{x_j}\partial_{x_k}g-\right.
$$
$$
\left.-\partial_{x_j}\partial_{x_k}f\otimes
\partial_{\xi_j}\partial_{\xi_k}g\right]( X, X)-\frac{1}{2i}\sum_{j,k=1}^nB_{jk}(x)[(\partial_{\xi_j}f)\otimes(\partial_{\xi_k}g)]
( X, X).
$$

In particular we have
\begin{equation}\label{mMcom}
f\circ^B\!g-g\circ^B\!f\cong\frac{1}{i}\{f,g\}-\frac{1}{i}\sum_{j,k=1}^nB_{jk}
(\partial_{\xi_j}f)(\partial_{\xi_k}g)=\frac{1}{i}\{f,g\}_B
\end{equation}
$(\text{\sl mod.}\;\; S^{m_1+m_2-3(\rho-\delta)}_{\rho,\delta}(\Xi))$, with $\{\cdot,\cdot\}_B$ the Poisson bracket associated to the symplectic form $\sigma_B$ defined in (\ref{sigmaB}).

\subsection{The parametrix}
One of the main tools in pseudodifferential theory is the parametrix.
\begin{definition}
A symbol $a\in S^m_{\rho,\delta}(\Xi)$ is called {\rm elliptic} when there exists two positive constants $R$ and $C$ such that
$$
C<\xi>^m\leq |a(x,\xi)|,\quad\forall x\in\mathcal{X},\quad\forall\xi\in\mathcal{X}^*\;\text{\sl  with }|\xi|\geq R.
$$
The operator $\mathfrak{Op}^A(a)$ will also be called {\rm elliptic}.
\end{definition}

\begin{theorem}\label{parametrica}
Let $a\in S^m_{\rho,\delta}(\Xi)$, $0\le\delta<\rho\le1$ be an elliptic symbol.
Then there exists $b\in S^{-m}_{\rho,\delta}(\Xi)$ such that $a\circ^B b-1,\,b\circ^B a-1\in S^{-\infty}(\Xi)$. Thus for any vector potential $A$ with components of class $C^\infty_\text{\sl pol}(\mathcal{X})$ associated to $B$, $\mathfrak{Op}^A(b)$ is a parametrix for $\mathfrak{Op}^A(a)$.
\end{theorem}

\begin{proof}
We construct first an approximation for $b$. 

Let $\chi\in C^\infty(\mathcal X^*)$, $\chi(\xi)=0$ if $\vert\xi\vert\le R$, $\chi(\xi)=1$ if $\vert\xi\vert\ge 2R$. We define $b_0(x,\xi):=\chi(\xi)a(x,\xi)^{-1}\in C^{\infty}(\Xi)$. 
It is easy to see by recurrence that for $\alpha,\beta\in\mathbb N^n$ and $(x,\xi)\in\Xi$
\begin{equation}\label{ineg}
\vert\left(\partial^\beta_x\partial^\alpha_\xi b_0\right)(x,\xi)\vert\le C_{\alpha,\beta}<\xi>^{-m-\rho\vert\alpha\vert+\delta\vert\beta\vert},
\end{equation}
and thus $b_0\in S^{-m}_{\rho,\delta}(\Xi)$. 

Thus, by Theorem \ref{50}, we have an asymptotic development $\ b_0\circ^B a\;\sim\;\sum_{j=0}^\infty h_j(b_0,a)$. Since $h_0(b_0,a)=b_0a=1+(\chi-1)$, one has $\ r_0:=b_0\circ^B a-1\in S^{-(\rho-\delta)}_{\rho,\delta}(\Xi)$. We denote by $r\in S^{0}_{\rho,\delta}(\Xi)$ an asymptotic sum of the series $1+\sum_{j=1}^{\infty}(-1)^jr_0\circ^B\dots\circ^B r_0$ (the $j$'th term contains $j$ factors and it belongs to $S^{-j(\rho-\delta)}_{\rho,\delta}(\Xi)$). Then $b:=r\circ^B b_0\in S^{-m}_{\rho,\delta}(\Xi)$, and from the equality $1+r_0=b_0\circ^B a$ it follows immediatly that $b\circ^B a-1\in S^{-\infty}(\Xi)$. In the same way one constructs $b'\in S^{m'}_{\rho,\delta}(\Xi)$ such that $a\circ^B b'-1\in S^{-\infty}(\Xi)$. Then 
\begin{equation*}
b-b'=b\circ^B(1-a\circ^B b')-(1-b\circ^B a)\circ^B b'\in S^{-\infty}(\Xi),
\end{equation*} 
thus $\,a\circ^B b-1\in S^{-\infty}(\Xi)$.
\end{proof}

\section{$\boldsymbol{L^2}$-continuity}

\begin{theorem}\label{109}
Suppose that the magnetic field $B$ has components of class $BC^\infty$. Let $f\in S^0_{\rho,\rho}(\Xi)$ for some $\rho\in[0,1)$. Then $\mathfrak{Op}^A(f)\in\mathcal{B}(L^2(\mathcal{X}))$ and we have the inequality
\begin{equation}\label{inequality}
\left\|\mathfrak{Op}^A(f)\right\|_{_{\mathcal{B}(L^2(\mathcal{X}))}}\,\leq\,c(n)\,\underset{|\alpha|\leq p(n)}{\sup}\;\underset{|\beta|\leq p(n)}{\sup}\;\underset{X\in\Xi}{\sup}<\xi>^{\rho(|\beta|-|\alpha|)}\left|\partial^\alpha_x\partial^\beta_\xi f(X)\right|,
\end{equation}
where $c(n),p(n)$ are constants depending only on $n$, that can be determined explicitly.
\end{theorem}

The proof we give makes use of an idea of L. Boutet de Monvel (\cite{LBM-1}); more precisely we have to use the Cotlar-Knapp-Stein lemma in a way adapted to our specific situation. We shall need the following remark.

\begin{remark}\label{111}
Let linear operators $\{T_j\}_{j\in\mathbb{N}}$ and $T$ in $\mathcal{B}(\mathcal{S}(\mathcal{X}),\mathcal{S}^*(\mathcal{X}))$ be given, with integral kernels $\{K_j\}_{j\in\mathbb{N}}$ and $K$, respectively, in $\mathcal{S}^*(\mathcal{X}\times\mathcal{X})$. Assume that
\begin{enumerate}
\item[i)] $T_j\in\mathcal{B}(L^2(\mathcal{X})),\;\forall j\in\mathbb{N}$ and there exists $C>0$ such that $\|T_j\|\leq C,\;\forall j\in\mathbb{N}$;
\item[ii)] $\underset{j\rightarrow\infty}{\lim}K_j=K$ in $\mathcal{S}^*(\mathcal{X}\times\mathcal{X})$.
\end{enumerate}
Then $T\in\mathcal{B}(L^2(\mathcal{X}))$ and $\|T\|\leq C$ (for the same constant $C$ as above).
\end{remark}

\begin{proofT}
By Remark \ref{111} we may suppose that $f\in\mathcal{S}(\Xi)$. Indeed, let us choose $\chi\in C^\infty_0(\Xi)$ such that $\chi(X)=1$ for $|x|+|\xi|\leq1$, and for any $\epsilon\in(0,1)$ let us define $\chi_\epsilon(X):=\chi(\epsilon X)$ and $f_\epsilon:=\chi_\epsilon f\in\mathcal{S}(\Xi)$. For any multiindices $\alpha$ and $\beta$ let us define
$$
C_f(\alpha,\beta)\;:=\;\underset{\alpha^\prime\leq\alpha}{\sup}\;\underset{\beta^\prime\leq\beta}{\sup}\;\underset{X\in\Xi}{\sup}<\xi>^{\rho(|\beta^\prime|-|\alpha^\prime|)}\left|\partial^{\alpha^\prime}_x\partial^{\beta^\prime}_\xi f(X)\right|.
$$
Then it is easy to check that there exist $C_{\alpha\beta}\in\mathbb{R}_+$ such that
$$
<\xi>^{\rho(|\beta|-|\alpha|)}\left|\partial^\alpha_x\partial^\beta_\xi f_\epsilon(X)\right|\;\leq
$$
$$
\leq\;C_{\alpha\beta}C_f(\alpha,\beta)\underset{\alpha^\prime\leq\alpha}{\max}\;\underset{\beta^\prime\leq\beta}{\max}\left[\,\epsilon^{|\beta-\beta^\prime|}<\xi>^{\rho(|\beta|-|\alpha|)}<\xi>^{\rho(|\alpha^\prime|-|\beta^\prime|)}\right].
$$
As $\epsilon|\xi|\leq c^\prime$ on the support of $f_\epsilon$, we obtain the estimation $\epsilon\leq c<\xi>^{-1}$ (valid on the support of $f_\epsilon$) and thus the right-hand side above is bounded by
$$
C^\prime_f(\alpha,\beta)<\xi>^{-|\beta-\beta^\prime|}<\xi>^{\rho(|\beta|-|\beta^\prime|)}<\xi>^{\rho(|\alpha^\prime|-|\alpha|)}\;\leq\;C''_f(\alpha,\beta).
$$
Consequently, if we prove (\ref{inequality}) for $f_\epsilon$ (for any $\epsilon\in(0,1]$), applying Remark \ref{111} to the family $\{\mathfrak{Op}^A(f_\epsilon)\}_{\epsilon\in(0,1]}$ gives the conclusion of Theorem \ref{109} for $f$. Thus we suppose $f\in\mathcal{S}(\Xi)$. 
Choose $u\in\mathcal{S}(\mathcal{X})\subset L^2(\mathcal{X})$ and define $\lambda(\xi):=<\xi>^{\rho}$.

\noindent
{\sl 1-st Step.} Given $N\in\mathbb{N}$, that we shall fix later on, we consider the differential operator
$$
L_\xi\;:=\;\left[1\,+\,\lambda(\xi)^{2N}|x-y|^{2N}\right]^{-1}\left[1\,+\,\lambda(\xi)^{2N}(-\Delta_\xi)^{2N}\right],
$$
that satisfies $L_\xi\exp\{i<x-y,\xi>\}=\exp\{i<x-y,\xi>\}$. After some integrations by parts we get
$$
[\mathfrak{Op}^A(f)u](x)\,=\,\int\!\!\!\int_{\mathcal{X}\times\mathcal{X}^*}dy
\dbar\xi\,e^{i<x-y,\xi>-i\Gamma^A([x,y])}g(x,y;\xi)u(y),
$$
where
$
g(x,y;\xi)\,:=\,[(^tL_\xi)f]((x+y)/2,\xi)
$.
One easily checks the inequalities
$$
\left|\left(\partial^\alpha_\xi\partial^\beta_x\partial^\gamma_yg\right)(x,y;\xi)\right|\,\leq\,
C_{\alpha\beta\gamma}\lambda(\xi)^{|\beta|+|\gamma|}\left(1\,+\,\lambda(\xi)|x-y|
\right)^{-2N}\times
$$
$$
\times\underset{|\alpha^\prime|\leq|\alpha|+2N}{\sup}\;\underset{|\beta^\prime|
\leq|\beta|+|\gamma|}{\sup}\;\underset{(y,\xi)\in\Xi}{\sup}\;<\xi>^{\rho(|\alpha^\prime|
-|\beta^\prime|)}\left|\left(\partial^{\alpha^\prime}_\xi\partial^{\beta^\prime}_yf
\right)(y,\xi)\right|.
$$

\noindent
{\sl 2-nd Step.} Let us write
$
[\mathfrak{Op}^A(f)u](x)$ as the integral over $\xi\in\mathcal{X}^*$ of
$$
[P_\xi u](x)\;:=\;\int_{\mathcal{X}}dy\,\exp\{i<x-y,\xi>-i\Gamma^A([x,y])\}g(x,y;\xi)u(y).
$$
We denote by $K_\xi$ the integral kernel of the operator $P_\xi$; we have the bound 
$$
|K_\xi(x,y)|\leq C(1+|x-y|)^{-2N}.
$$
Choosing $N>n/2$, we find a constant $C_0$ such that $\|P_\xi\|_{_{\mathcal{B}(L^2(\mathcal{X}))}}\leq C_0$ for any $\xi\in\mathcal{X}^*$.

\noindent
{\sl 3-rd Step.} For $v\in\mathcal{S}(\mathcal{X})$ and $\eta\in\mathcal X^*$ we have
$$
[P_\eta^*v](y)\;=\;\int_{\mathcal{X}}dz\,\exp\{i<y-z,\eta>-i\Gamma^A([y,z])\}
\overline{g(z,y;\eta)}v(z).
$$
Thus the integral kernel of the operator $P_\xi P^*_\eta$ is given by
$$
K_{\xi,\eta}(x,y)\,:=\,\exp\{i[<x,\xi>-<y,\eta>-\Gamma^A([x,y])]\}\times
$$
$$
\times
\int_{\mathcal{X}}dz\,\exp\{i<z,\eta-\xi>-i\Gamma^B(<x,z,y>)\}\,g(x,z;\xi)
\overline{g(y,z;\eta)}.
$$

\noindent
{\sl 4-th Step.} We consider the following differential operator
$$
M_z\;:=\;\left[1\,+\,\left(\frac{|\xi-\eta|}{\lambda(\xi)+\lambda(\eta)}\right)^2
\right]^{-1}\left]1\,-\,(\lambda(\xi)+\lambda(\eta))^{-2}\Delta_z\right],
$$
satisfying $M_z\exp\{i<z,\eta-\xi>\}=\exp\{i<z,\eta-\xi>\}$. Thus for any $k\in\mathbb{N}$ we get
$$
\exp\{-i[<x,\xi>-<y,\eta>-\Gamma^A([x,y])]\}K_{\xi,\eta}(x,y)\,:=
$$
$$
=\,\int_{\mathcal{X}}dz\,\exp\{i<z,\eta-\xi>\}\,M_z^k\left[
\exp\{-i\Gamma^B(<x,z,y>)\}\,
g(x,z;\xi)\overline{g(y,z;\eta)}\right].
$$
Using Remark \ref{46-r} and the bounds for the derivatives of $g$ (obtained above) we get for any $k\leq N$
$$
\left|M_z^k\left[\exp\{-i\Gamma^B(<x,z,y>)\}\,g(x,z;\xi)\overline{g(y,z;\eta)}
\right]\right|\,\leq
$$
$$
\leq\,C_k(f)\left(1+\frac{|\xi-\eta|}{\lambda(\xi)+\lambda(\eta)}\right)^{-2k}
[1\,+\,\lambda(\xi)|x-z|]^{-2(N-k)}[1\,+\,\lambda(\eta)|y-z|]^{-2(N-k)},
$$
where $C_k(f)$ is a linear combination of a finite number (depending only on $k$) of products of two seminorms of $f$ in $S^0_{\rho,\rho}$. Choosing now $N$ large enough in order to have $N-k>n/2$, we easily obtain the inequality
$$
\int_{\mathcal{X}}dx\left|K_{\xi,\eta}(x,y)\right|\,+\int_{\mathcal{X}}dy
\left|K_{\xi,\eta}(x,y)\right|\leq 
$$
$$
\leq C^\prime_k(f)[\lambda(\xi)
\lambda(\eta)]^{-n}\left(1+\frac{|\xi-\eta|}{\lambda(\xi)+\lambda(\eta)}
\right)^{-2k}\;=:\;h^2(\xi,\eta).
$$

\noindent
{\sl 5-th Step.} We have thus obtained the following estimations:
$$
\|P_\xi\|_{_{\mathcal{B}(L^2(\mathcal{X}))}}\;\leq\;C_0,\quad\forall\xi\in\mathcal{X}^*,\;\,\text{\sl if } N>\frac{n}{2},
$$
$$
\|P_\xi P^*_\eta\|_{_{\mathcal{B}(L^2(\mathcal{X}))}}\;\leq\;h^2(\xi,\eta),\quad\forall(\xi,\eta)\in(\mathcal{X}^*)^2,\;\,\text{\sl if } N>\frac{n}{2}+k.
$$
The conclusion of Theorem \ref{109} will follow now from the Cotlar-Knapp-Stein lemma, once we have proved that there exists $k\in\mathbb{N}$ such that
$$
\underset{\xi\in\mathcal{X}^*}{\sup}\int_{\mathcal{X}^*}d\eta\,\lambda(\xi)^{-n/2}\lambda(\eta)^{-n/2}\left(1\,+\,\frac{|\xi-\eta|}{\lambda(\xi)+\lambda(\eta)}\right)^{-k}\;<\;\infty.
$$
We shall decompose this integral on three subdomains of $\mathcal{X}^*$: $\{|\eta|\leq|\xi|/2\}=:D_1$, $\{|\eta|\geq2|\xi|\}=:D_2$ and $\{|\xi|/2<|\eta|<2|\xi|\}=:D_3$. We get
$$
\mathcal{I}_1\leq C\int_{D_1}d\eta\,\lambda(\xi)^{-n/2}\left(1+\frac{|\xi-\eta|}
{\lambda(\xi)}\right)^{-k}\leq\,C^\prime\lambda(\xi)^{-n/2}
(1+|\xi|\lambda(\xi)^{-1})^{-k}|\xi|^n\,\leq
$$
$$
\leq\,C^{\prime\prime}<\xi>^{n(1-\rho/2)-k(1-\rho)}\,\leq\,C^{\prime\prime\prime},
\quad\text{\sl if}\quad k\geq\frac{n(2-\rho)}{2(1-\rho)},
$$
$$
\mathcal{I}_2\leq\,C\int_{D_2}d\eta\,[\lambda(\xi)\lambda(\eta)]^{-n/2}
\left(1+\frac{|\xi-\eta|}
{\lambda(\xi)+\lambda(\eta)}\right)^{-k}\,\leq
$$
$$
\leq\,C^\prime(1+|\xi|)^{-n\rho/2}\int_{D_2}d\eta\,(1+|\eta|)^{-n\rho/2-k(1-\rho)}\,
\leq\,C^{\prime\prime},\quad\text{\sl if}\quad k>\frac{n(2-\rho)}{2(1-\rho)},
$$
$$
\mathcal{I}_3\leq\,C\lambda(\xi)^{-n}\int_{D_3}d\eta\,\left(1\,+\,\frac{|\xi-\eta|}
{\lambda(\xi)}\right)^{-k}\,\leq\,C^\prime\int_{\mathcal{X}^*}d\eta\,
(1+|\eta|)^{-k}<\infty,
$$
if we take $k>n$. As we have
$
n(2-\rho)/2(1-\rho)\;\geq\;n
$
for any $\rho\geq0$, we conclude that we can fix $k\in\mathbb{N}$ and $N\in\mathbb{N}$ such that
$
k>n(2-\rho)/2(1-\rho)$ and $N>k+n/2$. This finishes the proof.
\end{proofT}

\begin{remark}\label{3.8}
Theorem \ref{109} remains true also for symbols of class $S^0_{\rho,\delta}(\Xi)$ with $0\leq\delta<\rho\leq1$, due to the obvious inclusion $S^0_{\rho,\delta}(\Xi)\subset S^0_{\delta,\delta}(\Xi)$
\end{remark}

\section{Sobolev spaces}\label{S4}

In this section we shall suppose that the components of the magnetic field $B$ are of class $BC^\infty(\mathcal{X})$; we shall work in a Schr\"odinger representation $\mathfrak{Op}^A$ associated to a vector potential $A$ (such that $B=dA$) with components of class $C^\infty_{\text{\sl pol}}(\Xi)$. 
\begin{definition}\label{61}
$$
\boldsymbol{\Psi}^{A,m}_{\rho,\delta}(\mathcal{X}):=\left\{
\mathfrak{Op}^A(a)\;\mid\;a\in S^m_{\rho,\delta}(\Xi)\right\},\; \boldsymbol{\Psi}^{A,m}(\mathcal{X}):=
\boldsymbol{\Psi}^{A,m}_{1,0}(\mathcal{X}).
$$
\end{definition}
\begin{definition}\label{62}
For any $s\in\mathbb{R}_+$ we set
$$
p_s(\xi):=<\xi>^s\in S^s(\Xi),\qquad
\mathfrak{P}_s:=\mathfrak{Op}^A(p_s)\in\boldsymbol{\Psi}^{A,s}
(\Xi),
$$
$$
\boldsymbol{H}^s_A(\mathcal{X}):=\left\{u\in L^2(\mathcal{X})\;\mid\;
\mathfrak{P}_su\in L^2(\mathcal{X})\right\}.
$$
\end{definition}
\noindent
As a consequence of the Proposition below, {\it the magneticSobolev space} $\boldsymbol{H}^s_A(\mathcal{X})$ (for $s\in\mathbb{R}_+$) may be defined using any elliptic operator of order $s$.

\begin{proposition}\label{63}
For any $s\in\mathbb{R}_+$ we have:
\begin{enumerate}
\item $\boldsymbol{H}^s_A(\mathcal{X})$ is a Hilbert space for the scalar product
$$
(u,v)_{s,A}:=(\mathfrak{P}_su,\mathfrak{P}_sv)_{L^2}+(u,v)
_{L^2},\quad\forall u,v\in\boldsymbol{H}^s_A(\mathcal{X}).
$$
\item If $0\leq\delta<\rho\leq1$, then $T\in\boldsymbol{\Psi}^{A,s}_{\rho,\delta}(\mathcal{X})$ is bounded from $\boldsymbol{H}^s_A(\mathcal{X})$ to $L^2(\mathcal{X})$.
\item For any elliptic operator $T\in\boldsymbol{\Psi}^{A,s}_{\rho,\delta}(\mathcal{X})$, the map
$$
\boldsymbol{H}^s_A(\mathcal{X})\ni u\mapsto\|u\|^T_{s,A}:=
\left\{\|Tu\|^2_{L^2}+\|u\|^2_{L^2}\right\}^{1/2}
$$
defines an equivalent norm on $\boldsymbol{H}^s_A(\mathcal{X})$.
\end{enumerate}
\end{proposition}

\begin{proof}
The first conclusion is evident, because $\mathfrak{P}_s
\in\mathcal{B}[\mathcal{S}^*(\mathcal X)]$ (Lemma \ref{43}). 
Let us fix now an operator $T\in\boldsymbol{\Psi}^{A,s}_{\rho,\delta}(\mathcal{X})$. We notice that, being elliptic, $\mathfrak{P}_s\in\boldsymbol{\Psi}^{A,s}(\mathcal{X})$ has a parametrix $\mathfrak{Q}_s\in\boldsymbol{\Psi}^{A,-s}(\mathcal{X})$ (Theorem \ref{parametrica}), i.e. there exists $\mathfrak{R}_s\in\boldsymbol{\Psi}^{A,-\infty}(\mathcal{X})$ such that $\mathfrak{Q}_s\mathfrak{P}_s=\boldsymbol{1}+\mathfrak{R}_s$. Thus, for any $u\in \boldsymbol{H}^s_A(\mathcal{X})$ we have $Tu=(T\mathfrak{Q}_s)\mathfrak{P}_su-(T\mathfrak{R}_s)u$. Theorem \ref{50} implies that $T\mathfrak{Q}_s$ and $T\mathfrak{R}_s$ belong to $\boldsymbol{\Psi}^{A,0}_{\rho,\delta}(\mathcal{X})$. We use Remark \ref{3.8} to deduce that $Tu\in L^2(\mathcal{X})$ and $\|Tu\|^2_{L^2}\leq C_0\left(\|\mathfrak{P}_su\|^2_{L^2}+\|u\|^2_{L^2}
\right)$. Thus we get the second conclusion of the Proposition and the inequality $\|u\|^T_{s,A}\leq C\|u\|_{s,A}$ for any $u\in \boldsymbol{H}^s_A(\mathcal{X})$. 
In order to get the reversed inequality and thus the third conclusion of the Lemma, we have to suppose $T$ elliptic. Then there is a parametrix $S\in \boldsymbol{\Psi}^{A,-s}_{\rho,\delta}(\mathcal{X})$ for $T$ such that $ST=\boldsymbol{1}+R$, with $R\in \boldsymbol{\Psi}^{A,-\infty}(\mathcal{X})$. For $u\in L^2(\mathcal{X})$ such that $Tu\in L^2(\mathcal{X})$, we obtain $\mathfrak{P}_su=(\mathfrak{P}_sS)Tu-(\mathfrak{P}_sR)u\in L^2(\mathcal{X})$ and also $\|u\|_{s,A}\leq C\|u\|^\prime_{s,A}$.
\end{proof}

\begin{lemma}\label{68}
If $0\leq s\leq t$, we have a continuous embedding $\boldsymbol{H}^{t}_{A}(\mathcal{X})\hookrightarrow\boldsymbol{H}^{s}_{A}
(\mathcal{X})$.
\end{lemma}

\begin{proof}
Assume that  $u\in \boldsymbol{H}^{t}_{A}(\mathcal{X})$. By definition of the Sobolev space, it follows that $u\in L^2(\mathcal{X})$ and $\mathfrak{P}_tu\in L^2(\mathcal{X})$. Making once again use of the parametrix $\mathfrak{Q}_t\in\boldsymbol{\Psi}^{A,-t}(\mathcal{X})$ of $\mathfrak{P}_t$, we deduce that there exists some $\mathfrak{R}_t\in\boldsymbol{\Psi}^{A,-\infty}(\mathcal{X})$ such that $u=\mathfrak{Q}_t\mathfrak{P}_tu+\mathfrak{R}_tu$. Thus $\mathfrak{P}_su=
\mathfrak{P}_s\mathfrak{Q}_t\mathfrak{P}_tu+\mathfrak{P}_s\mathfrak{R}_tu$. Using Theorem \ref{50} we get that $\mathfrak{P}_s\mathfrak{Q}_t\in\boldsymbol{\Psi}^{A,-(t-s)}(\mathcal{X})$ and $\mathfrak{P}_s\mathfrak{R}_t\in\boldsymbol{\Psi}^{A,-\infty}(\mathcal{X})$, so that by Remark \ref{3.8} we deduce that $\|\mathfrak{P}_su\|_{L^2}\leq C(\|\mathfrak{P}_tu\|_{L^2}+\|u\|_{L^2})$.
\end{proof}

\begin{lemma}\label{64}
Suppose given $s\in\mathbb{R}_+$, $m\leq s$ and $T\in \boldsymbol{\Psi}^{A,m}_{\rho,\delta}(\mathcal{X})$. Then $T$ is a bounded operator from $\boldsymbol{H}^{s}_{A}(\mathcal{X})$ to $\boldsymbol{H}^{s-m}_{A}(\mathcal{X})$.
\end{lemma}
\begin{proof}
Consider $u\in \boldsymbol{H}^{s}_{A}(\mathcal{X})$. Since $m\leq s$ we also have $T\in \boldsymbol{\Psi}^{A,s}_{\rho,\delta}(\mathcal{X})$, thus $Tu\in L^2(\mathcal{X})$ and $\|Tu\|_{L^2(\mathcal{X})}\leq C\|u\|_{s,A}$. Moreover, due to our Theorem \ref{50}, $\mathfrak{P}_{s-m}T\in \boldsymbol{\Psi}^{A,s}_{\rho,\delta}(\mathcal{X})$, so that we have $\mathfrak{P}_{s-m}Tu\in L^2(\mathcal{X})$ and $\|\mathfrak{P}_{s-m}Tu\|_{L^2(\mathcal{X})}\leq C\|u\|_{s,A}$. We conclude that $Tu\in\boldsymbol{H}^{s-m}_{A}(\mathcal{X})$ and $\|Tu\|_{s-m,A}\leq C\|u\|_{s,A}$.
\end{proof}

\begin{lemma}\label{67}
If the vector potential $A$ has components of class $C^\infty_{\text{\sl pol}}(\Xi)$, then $R\in\boldsymbol{\Psi}^{A,-\infty}(\mathcal{X})$ defines a linear continuous operator from $\mathcal{S}^*(\mathcal{X})$ into $C^\infty(\mathcal{X})$.
\end{lemma}
\begin{proof}
Due to our hypothesis on $R$, it exists a symbol $r\in S^{-\infty}(\Xi)$ such that for any $\varphi\in\mathcal{S}(\mathcal{X})$ we have 
$
(R\varphi)(x)=\int_\mathcal{X}dy\,K_r(x,y)\varphi(y)
$,
with
$$
K_r(x,y):=e^{-i\Gamma^A([x,y])}\int_\mathcal{X}\dbar\xi\,e^{i<x-y,\xi>}
r\left(\frac{x+y}{2},\xi\right).
$$
In a similar way, the adjoint $R^*$ of $R$ is an integral operator with kernel
$$
\tilde{K}_r(x,y):=e^{-i\Gamma^A([x,y])}\int_\mathcal{X}\dbar\xi\,e^{i<x-y,\xi>}
\overline{r}\left(\frac{x+y}{2},\xi\right).
$$
By our hypothesis on $A$, $\tilde{K}_r\in C^\infty(\mathcal{X}\times\mathcal{X})$, and $\forall\alpha,\beta\in\mathbb{N}^n$ there exist natural numbers $N_\alpha$ and $N_\beta$ such that
$$
\left|\left(\partial^\alpha_x\partial^\beta_y\tilde{K}_r\right)(x,y)\right|
\leq C_{\alpha,\beta,p}<x>^{N_\alpha}<y>^{N_\beta}<x-y>^{-p},\quad\forall p\in\mathbb{N}.
$$
In particular, $\forall\beta\in\mathbb{N}^n$, $\partial^\beta_y\tilde{K}_r(\cdot,y)\in\mathcal{S}(\mathcal{X})$ uniformly, for $y$ in any compact subset. Let $\Omega$ be a relatively compact open subset of $\mathcal{X}$ and $\varphi\in C^\infty_0(\mathcal{X})$, $\text{\sl supp}\varphi\subset\Omega$. For any $u\in\mathcal{S}^*(\mathcal{X})$
$$
(Ru,\varphi)=(u,R^*\varphi)=\int_\mathcal{X}dy\,(u(\cdot),\tilde{K}_r(\cdot,y))
\overline{\varphi(y)}
=\int_\Omega dy\,f(y)\overline{\varphi(y)},
$$
where $f(y):=(u(\cdot),\tilde{K}_r(\cdot,y))$ is obviously in $C^\infty(\Omega)$. Thus $Ru\in C^\infty(\mathcal{X})$, and the continuity follows from the closed graph theorem.
\end{proof}

\begin{lemma}\label{69}
For any $s\in\mathbb{R}_+$ we have the continuous, dense embeddings $\mathcal{S}(\mathcal{X})\hookrightarrow\boldsymbol{H}^{s}_{A}(\mathcal{X})
\hookrightarrow\mathcal{S}^*(\mathcal{X})$.
\end{lemma}

\begin{proof}
The existence and continuity of the two embeddings is evident. Let us prove the density of $\mathcal{S}(\mathcal{X})$ in $\boldsymbol{H}^{s}_{A}(\mathcal{X})$. Take $u\in \boldsymbol{H}^{s}_{A}(\mathcal{X})$. Let us choose a sequence $\{v_j\}_{j\in\mathbb{N}}\subset\mathcal{S}(\mathcal{X})$ that converges to $\mathfrak{P}_su$ in $L^2(\mathcal{X})$. We consider once again the parametrix $\mathfrak{Q}_s$ of $\mathfrak{P}_s$, fix a cut-off function $\chi\in C^\infty_0(\mathcal{X})$ with $\chi(x)=1$ for $|x|\leq1$ and set $\chi_j(x):=\chi(x/j)$. We define $u_j:=\mathfrak{Q}_sv_j-\chi_j\mathfrak{R}_su$. Evidently $u_j\in\mathcal{S}(\mathcal{X})$ and $\mathfrak{Q}_sv_j$ converges to $\mathfrak{Q}_s\mathfrak{P}_su$ in $\boldsymbol{H}^{s}_{A}(\mathcal{X})$. But $\mathfrak{Q}_s\mathfrak{P}_su=u+\mathfrak{R}_su$, so that the density conclusion will follow if we prove that $\chi_j\mathfrak{R}_su$ converges to $\mathfrak{R}_su$ in $\boldsymbol{H}^{s}_{A}(\mathcal{X})$. Let us put $\nu_t(\xi):=\sum_{j=1}^n\xi^{2t}_j$ for $t\in\mathbb{N},\;2t\geq s$ and $T:=\mathfrak{Op}^A(\nu_t)=\sum_{j=1}^n(D_j-A_j)^{2t}.$
Then $T$ is a differential operator of order $2t$ and an elliptic operator in $\boldsymbol{\Psi}^{A,2t}(\mathcal X)$. A simple computation shows that $T\chi_j\mathfrak{R}_su-\chi_jT\mathfrak{R}_su$ is a finite sum of terms of the form $(D_k-A_k)^{m_k}\mathfrak{R}_su$, each one multiplied by a bounded function of $x$ (containing derivatives of $\chi$) and a strictly negative power of $j$. Thus $T\chi_j\mathfrak{R}_su$ converges to $T\mathfrak{R}_su$ in $L^2(\mathcal{X})$. Using Proposition \ref{63}(3), we deduce the convergence of $\chi_j\mathfrak{R}_su$ to $\mathfrak{R}_su$ in $\boldsymbol{H}^{2t}_{A}(\mathcal{X})$, and thus also in $\boldsymbol{H}^{s}_{A}(\mathcal{X})$.
\end{proof}

\begin{definition}\label{65}
For $s\in\mathbb{R}_+$, we denote by $\boldsymbol{H}^{-s}_{A}(\mathcal{X})$ the anti-dual of $\boldsymbol{H}^{s}_{A}(\mathcal{X})$ endowed with the natural norm (that induces a scalar product):
$$
\|u\|_{-s,A}:=\underset{\varphi\in\boldsymbol{H}^{s}_{A}\setminus\{0\}}{\sup}
\frac{|(u,\varphi)|}{\|\varphi\|_{s,A}}.
$$
\end{definition}
\begin{proposition}\label{72}
If $s_1\leq s_2$ are two real numbers, then we have a continuous embedding $\boldsymbol{H}^{s_2}_{A}(\mathcal{X})\hookrightarrow\boldsymbol{H}^{s_1}_{A}(\mathcal{X})$.
\end{proposition}
\begin{proof}
Just use Lemma \ref{68} and a duality argument.
\end{proof}

\begin{proposition}\label{66}
Let us fix $s\in\mathbb{R}_+\setminus\{0\}$.
\begin{enumerate}
\item If $u\in\mathcal{S}^*(\mathcal{X})$ is of the form $u=\mathfrak{P}_sv+w$, with $v$ and $w$ from $L^2(\mathcal{X})$, then $u\in\boldsymbol{H}^{-s}_{A}(\mathcal{X})$ and 
$
\|u\|_{-s,A}\leq\left(\|v\|^2_{L^2}+\|w\|^2_{L^2}\right)^{1/2}
$.
\item Reciprocally, if $u\in\boldsymbol{H}^{-s}_{A}(\mathcal{X})$, then there exists $v$ and $w$ in $L^2(\mathcal{X})$ such that $u=\mathfrak{P}_sv+w$ and
$
\left(\|v\|^2_{L^2}+\|w\|^2_{L^2}\right)^{1/2}\leq\|u\|_{-s,A}
$.
\end{enumerate}
In conclusion we have
$$
\boldsymbol{H}^{-s}_{A}(\mathcal{X})=\left\{u\in\mathcal{S}^*(\mathcal{X})\;
\mid\;\exists\,v,w\in L^2(\mathcal{X})\;\text{such that}
\;u=\mathfrak{P}_sv+w
\right\},
$$
and for $u\in\boldsymbol{H}^{-s}_{A}(\mathcal{X})$ and $v$, $w$ as above,
$
\|u\|_{-s,A}=\left(\|v\|^2_{L^2}+\|w\|^2_{L^2}\right)^{1/2}
$.
\end{proposition} 

\begin{proof}
1. For $u=\mathfrak{P}_sv+w\in\mathcal{S}^*(\mathcal{X})$, 
let us define the linear map
$$
L_u:\mathcal{S}(\mathcal{X})\rightarrow\mathbb{C},\qquad L_u(\varphi):=
(u,\varphi)=(\mathfrak{P}_sv+w,\varphi)=(v,\mathfrak{P}_s\varphi)+
(w,\varphi).
$$
We have
$
\left| L_u(\varphi)\right|\leq\|\mathfrak{P}_s\varphi\|\|v\|+
\|\varphi\|\|w\|\leq\left(\|v\|^2+\|w\|^2\right)^{1/2}\|\varphi\|_{s,A}
$.
It follows that $L_u$ is continuous for the $\boldsymbol{H}^{s}_{A}(\mathcal{X})$-norm on $\mathcal{S}(\mathcal{X})$, which is dense in $\boldsymbol{H}^{s}_{A}(\mathcal{X})$ (see Lemma \ref{69}). Thus $L_u\in\boldsymbol{H}^{-s}_{A}(\mathcal{X})$ and we have
$$
\|L_u\|_{-s,A}=\underset{\varphi\in\mathcal{S}\setminus\{0\}}{\sup}\frac
{|L_u(\varphi)|}{\|\varphi\|_{s,A}}\leq\left(
\|v\|^2+\|w\|^2\right)^{1/2}.
$$
We have got a map
$$
\mathfrak{i}:\mathcal{M}:=\left\{u\in\mathcal{S}^*(\mathcal{X})
\mid\exists\,v,w\in L^2(\mathcal{X})\;\text{such that}
\;u=\mathfrak{P}_sv+w\right\}\rightarrow\boldsymbol{H}^{-s}_{A}(\mathcal{X}),
$$
defined by the formula $\mathfrak{i}(u):=L_u$. It is well defined, since it does not depend on the representation $u=\mathfrak{P}_sv+w$. Moreover, we proved that it satisfies 
$
\|\mathfrak{i}(u)\|_{-s,A}\leq\left(
\|v\|^2+\|w\|^2\right)^{1/2}.
$
The map $\mathfrak{i}$ is injective, because
$
\mathfrak{i}(u)=0\,\text{for}\,u\in\mathcal{M}\;\Rightarrow\;L_u=0\;
\Rightarrow\;(u,\varphi)=0,\,\forall\varphi\in\mathcal{S}(\mathcal{X})
\;\Rightarrow\;u=0.
$

2. Let us prove that $\mathfrak{i}$ is also surjective. Take $L\in\boldsymbol{H}^{-s}_{A}(\mathcal{X})$; since $\mathcal{S}(\mathcal{X})$ is continuously embedded in $\boldsymbol{H}^{s}_{A}(\mathcal{X})$ (see Lemma \ref{69}), we have $L\in\mathcal{S}^*(\mathcal{X})$. We must find an element $u\in\mathcal{M}$ such that $\mathfrak{i}(u)\equiv L_u=L$.
We begin by defining the map $\Phi:\boldsymbol{H}^{s}_{A}(\mathcal{X})
\rightarrow L^2(\mathcal{X})\times L^2(\mathcal{X})$, $\,\Phi(\psi):=
\{\mathfrak{P}_s\psi,\psi\}$. Let $\mathcal{J}:=\text{\rm Range}\,\Phi$. It will be a closed subspace of $L^2(\mathcal{X})\times L^2(\mathcal{X})$,  $\mathfrak{P}_s$ being closed. If we consider on $L^2(\mathcal{X})\times L^2(\mathcal{X})$ the norm
$
\|\{f,g\}\|:=\left(\|f\|_{L^2}+\|g\|_{L^2}\right)^{1/2}
$,
it is evident that $\Phi:\boldsymbol{H}^{s}_{A}(\mathcal{X})
\rightarrow \mathcal{J}$ is an isometric isomorphism. Let us denote by $P$ the orthogonal projection $L^2(\mathcal{X})\times L^2(\mathcal{X})\rightarrow\mathcal{J}$, set
$
\tilde{L}:\mathcal{J}\rightarrow\mathbb{C}$, $\tilde{L}(\{p,q\}):=L(q)$
and $L^\prime:=\tilde{L}\circ P$. Notice that $L^\prime:=L\circ\Phi^{-1}\circ P:L^2(\mathcal{X})\times L^2(\mathcal{X})\rightarrow\mathbb{C}$ is antilinear and continuous and $\left.L^\prime\right|_{\mathcal{J}}=L\circ\Phi^{-1}=\tilde{L}$. Using the Riesz Theorem, we can find a unique pair $\{v,w\}\in L^2(\mathcal{X})\times L^2(\mathcal{X})$ such that $L^\prime(p,q)=(v,p)+(w,q)$ for any pair $\{p,q\}\in L^2(\mathcal{X})\times L^2(\mathcal{X})$. With this choice, we get for any $\varphi\in\mathcal{S}(\mathcal{X})$ the equality
$
L(\varphi)=L\left(\Phi^{-1}(\{
\mathfrak{P}_s\varphi,\varphi\})\right)=L^\prime(\{
\mathfrak{P}_s\varphi,\varphi\})=(v,\mathfrak{P}_s\varphi)+(w,\varphi)=(u,\varphi)
$
for $u=\mathfrak{P}_sv+w$. We also have 
$$
\left(\|v\|_{L^2}+\|w\|_{L^2}\right)^{1/2}=\|L^\prime\|_{(L^2\times L^2)^*}=
\underset{\{p,q\}\in(L^2\times L^2)\setminus\{0\}}{\sup}\frac{|L^\prime(\{p,q\})|}{\|\{p,q\}\|_{L^2\times L^2}} \leq
$$
$$
\leq \underset{\{p,q\}\in\mathcal{J}\setminus\{0\}}{\sup}
\frac{|\tilde{L}(\{p,q\})|}{\|\{p,q\}\|_{L^2\times L^2}}=
\underset{q\in\boldsymbol{H}^{s}_{A}\setminus\{0\}}{\sup}\frac{|L(q)|}{\|\Phi(q)\|_{L^2\times L^2}} =
$$
$$
=\underset{q\in\boldsymbol{H}^{s}_{A}\setminus\{0\}}{\sup}\frac{|L(q)|}{
\|q\|_{s,A}}=\|L\|_{-s,A}.
$$
\end{proof}

\begin{lemma}\label{70}
For any $s\in\mathbb{R}_+$ we have the continuous embeddings $\mathcal{S}(\mathcal{X})\hookrightarrow\boldsymbol{H}^{-s}_{A}(\mathcal{X})
\hookrightarrow\mathcal{S}^*(\mathcal{X})$, the space $\mathcal{S}(\mathcal{X})$ being dense in $\boldsymbol{H}^{-s}_{A}(\mathcal{X})$.
\end{lemma}

\begin{proof}
The continuous embeddings follow from Lemma \ref{69} and the definition of $\boldsymbol{H}^{-s}_{A}(\mathcal{X})$ as anti-dual of $\boldsymbol{H}^{s}_{A}(\mathcal{X})$. Let us fix now $u\in\boldsymbol{H}^{-s}_{A}(\mathcal{X})$. There exists a pair $\{v,w\}\in L^2(\mathcal{X})\times L^2(\mathcal{X})$ such that $u=\mathfrak{P}_sv+w$. Moreover, we may approach $v$ and $w$ in $L^2$-norm by sequences $\{v_j\}_{j\in\mathbb{N}}$ and $\{w_j\}_{j\in\mathbb{N}}$ from $\mathcal{S}(\mathcal{X})$. If we put $u_j:=\mathfrak{P}_sv_j+w_j$ and use Proposition \ref{66}, we deduce that
$$
\|u_j-u\|_{-s,A}\leq\left(\|v_j-v\|^2_{L^2}+\|w_j-w\|^2_{L^2}\right)^{1/2}
\underset{j\rightarrow\infty}{\longrightarrow}0.
$$
\end{proof}

\begin{proposition}\label{71}
For $s$ and $m$ real numbers, any $T\in\boldsymbol{\Psi}^{A,m}_{\rho,\delta}(\mathcal X)$ is bounded as operator from $\boldsymbol{H}^{s}_{A}(\mathcal{X})$ to $\boldsymbol{H}^{s-m}_{A}(\mathcal{X})$.
\end{proposition}

\begin{proof}
The case $s\geq 0$ and $m\leq s$ is the content of Lemma \ref{64}. By duality we obtain also the case $s\leq 0$ and $m\geq s$. If $s\geq 0$ and $m>s$, let us choose $u\in\boldsymbol{H}^{s}_{A}(\mathcal{X})$ and write once again $u=\mathfrak{Q}_s\mathfrak{P}_su+\mathfrak{R}_su$, with $\mathfrak{Q}_s\in\boldsymbol{\Psi}^{A,-s}$ and $\mathfrak{R}_s\in\boldsymbol{\Psi}^{A,-\infty}$. Thus $Tu=(T\mathfrak{Q}_s)\mathfrak{P}_su+T\mathfrak{R}_su$. We have $\mathfrak{P}_su\in L^2(\mathcal{X})$, $T\mathfrak{Q}_s\in\boldsymbol{\Psi}^{A,m-s}_{\rho,\delta}$ and $m-s\geq 0$, so that we conclude that $(T\mathfrak{Q}_s)\mathfrak{P}_su\in\boldsymbol{H}^{s-m}_{A}(\mathcal{X})$. We also remark that $T\mathfrak{R}_s\in\boldsymbol{\Psi}^{A,-\infty}$, so that $T\mathfrak{R}_su\in L^2(\mathcal{X})\subset\boldsymbol{H}^{s-m}_{A}(\mathcal{X})$. The case $s\leq 0$ and $m<s$ follows from Proposition \ref{66}.
\end{proof}

\begin{definition}\label{73-D}
We define $\boldsymbol{H}^{-\infty}_A(\mathcal{X}):=\cup_{s\in\mathbb{R}}\boldsymbol{H}^{s}_{A}(\mathcal{X})$ endowed with the inductive limit topology and
$\boldsymbol{H}^{\infty}_A(\mathcal{X}):=\cap_{s\in\mathbb{R}}\boldsymbol{H}^{s}_{A}(\mathcal{X})$ endowed with the projective limit topology.
\end{definition}

\begin{proposition}\label{73}
Let $T\in\boldsymbol{\Psi}^{A,m}_{\rho,\delta}(\mathcal X)$; then
\begin{enumerate}
\item $T$ induces linear continuous maps $\boldsymbol{H}^{-\infty}_A(\mathcal{X})\rightarrow\boldsymbol{H}^{-\infty}_A(\mathcal{X})$ and $\boldsymbol{H}^{\infty}_A(\mathcal{X})\rightarrow\boldsymbol{H}^{\infty}_A(\mathcal{X})$.
\item If $m=-\infty$, $T$ induces a linear continuous map $\boldsymbol{H}^{-\infty}_A(\mathcal{X})\rightarrow\boldsymbol{H}^{\infty}_A(\mathcal{X})$.
\item If $T$ is an elliptic operator and we have $u\in\boldsymbol{H}^{-\infty}_A(\mathcal{X})$ and $Tu\in\boldsymbol{H}^{s}_A(\mathcal{X})$, then $u\in\boldsymbol{H}^{s+m}_A(\mathcal{X})$.
\end{enumerate}
\end{proposition}

\begin{proof}
1. If $u\in\boldsymbol{H}^{-\infty}_A(\mathcal{X})$, there exists some real number $s$ such that $u\in\boldsymbol{H}^{s}_A(\mathcal{X})$. Thus $Tu\in\boldsymbol{H}^{s-m}_A(\mathcal{X})\subset\boldsymbol{H}^{-\infty}_A(\mathcal{X})$. If $u\in\boldsymbol{H}^{\infty}_A(\mathcal{X})$, then for any real number $s$ we have $u\in\boldsymbol{H}^{s}_A(\mathcal{X})$ and thus $Tu\in\boldsymbol{H}^{s-m}_A(\mathcal{X})$ for any $s\in\mathbb{R}$, i.e. $Tu\in\boldsymbol{H}^{\infty}_A(\mathcal{X})$. The continuity may be proved either directly or by using the closed graph theorem.

\noindent
2. If $T\in\boldsymbol{\Psi}^{A,-\infty}(\mathcal X)$ and $u\in\boldsymbol{H}^{-\infty}_A(\mathcal{X})$, then $u\in\boldsymbol{H}^{s}_A(\mathcal{X})$ for some real $s$ and $T\in\boldsymbol{\Psi}^{A,m}(\mathcal X)$ for any real $m$. We deduce that $Tu\in\boldsymbol{H}^{s-m}_A(\mathcal{X})$ for any $m\in\mathbb{R}$, i.e. $Tu\in\boldsymbol{H}^{\infty}_A(\mathcal{X})$.

\noindent
3. If $T$ is elliptic, there exists $S\in\boldsymbol{\Psi}^{A,-m}_{\rho,\delta}(\mathcal X)$ and $R\in\boldsymbol{\Psi}^{A,-\infty}(\mathcal X)$ such that $ST-\boldsymbol{1}=R$. If $u\in\boldsymbol{H}^{-\infty}_A(\mathcal{X})$ and $Tu\in\boldsymbol{H}^{s}_A(\mathcal{X})$, it follows that $u=S(Tu)-Ru\in\boldsymbol{H}^{s+m}_A(\mathcal{X})$.
\end{proof}
\begin{remark}\label{74}
The property (3) of Proposition \ref{73} may be completed as follows: If $T$ is an elliptic operator and $u\in\mathcal{S}^*(\mathcal{X})$,
then we have $\text{\rm sing supp}\,Tu=\text{\rm sing supp}\,u$; in
particular, if $Tu\in C^\infty(\mathcal{X})$ then $u\in
C^\infty(\mathcal{X})$.
\end{remark}
\noindent
The above statement follows from Lemma \ref{67} and from the fact that any operator $T\in\boldsymbol{\Psi}^{A,m}_{\rho,\delta}(\mathcal X)$ is {\it pseudo-local} (i.e. $\text{\rm sing supp}\,Tu\subset\text{\rm sing supp}\,u$). In fact, the integral kernel of $T$ is the product of the $C^\infty$ function $\exp\{-i\Gamma^A([x,y])\}$ and the distribution defined by the oscillatory integral
$
\int_\mathcal{X}\dbar\xi\,e^{i<x-y,\xi>}t\left(\frac{x+y}{2},\xi\right),
$
that is a $C^\infty$ function outside the diagonal of $\mathcal{X}\times\mathcal{X}$.

\begin{lemma}\label{75}
For any $m\in\mathbb{N}$ we have the equality
\begin{equation}\label{SobSp1}
\boldsymbol{H}^m_A(\mathcal{X})=\left\{u\in L^2(\mathcal{X})\,\mid\,(D-A)^\alpha u\in L^2(\mathcal{X}),\;\forall\alpha\in\mathbb{N}^n\;\text{with}\,|\alpha|\leq m\right\},
\end{equation}
where $(D-A)^\alpha=(D_1-A_1)^{\alpha_1}\cdots(D_n-A_n)^{\alpha_n}$. Moreover, we have the following equivalent norm on $\boldsymbol{H}^m_A(\mathcal{X})$:
$
\|u\|_{m,A}\;\sim\;\left(\sum_{|\alpha|\leq m}\left\|(D-A)^\alpha u\right\|^2_{L^2}\right)^{1/2}$.
\end{lemma}

\begin{proof}
Let us denote by $\mathcal{M}$ the linear space defined in (\ref{SobSp1}) endowed with the norm $\|\cdot\|_\mathcal{M}$ defined by the formula above. Remark that $D_j-A_j=\mathfrak{Op}^A(\xi_j)\in\boldsymbol{\Psi}^{A,1}(\mathcal{X})$, so that $(D-A)^\alpha\in\boldsymbol{\Psi}^{A,m}(\mathcal{X})$ for $|\alpha|\leq m$. In conclusion, for $u\in\boldsymbol{H}^m_A(\mathcal{X})$ we have
$
(D-A)^\alpha u\in L^2(\mathcal{X})$, and $\|(D-A)^\alpha u\|_{L^2}\leq C\|u\|_{m,A}
$.
Reciprocally, let $u\in\mathcal{M}$. We consider the operator 
$
E=\mathfrak{Op}^A(\nu_m)\in\boldsymbol{\Psi}^{A,2m}(\mathcal{X})
$ (with the notation introduced in the proof of Lemma \ref{69}),
that is elliptic. Thus we can find $F\in\boldsymbol{\Psi}^{A,-2m}(\mathcal{X})$ and $R\in\boldsymbol{\Psi}^{A,-\infty}(\mathcal{X})$ such that $FE-\boldsymbol{1}=R$. Due to our choice of $u$ and the definition of $\mathcal{M}$, we have $(D_j-A_j)^mu\in L^2(\mathcal{X})$ for any $j\in\{1,\ldots,n\}$ and thus $Eu\in\boldsymbol{H}^{-m}_A(\mathcal{X})$. We get
$
u\,=\,F(Eu)\,-\,Ru\,\in\,\boldsymbol{H}^{m}_A(\mathcal{X})
$,
and we finish the proof by the closed graph theorem, due to the fact that $\mathcal{M}$ is a Hilbert space.
\end{proof}
\begin{remark}
Let us point out here that for any real function $\phi\in C^\infty_{\text{\sl pol}}(\mathcal{X})$, the multiplication with $e^{i\phi}$ defines a unitary operator intertwining the Sobolev spaces $\boldsymbol{H}^{s}_{A+\nabla\phi}(\mathcal{X})$ and $\boldsymbol{H}^{s}_A(\mathcal{X})$.
\end{remark}

\section{Self-adjointness and semiboundedness}

\begin{theorem}\label{76} Suppose given a magnetic field $B$ with components of class $BC^\infty(\mathcal{X})$.
Let $p\in S^m_{\rho,\delta}(\Xi)$ be real with $m\geq0$, $0\leq\rho<\delta\leq 1$; if $m>0$ we also assume that $p$ is elliptic. Let us set $P:=\mathfrak{Op}^A(p)$ in a Schr\"{o}dinger representation defined by a vector potential $A$ associated to $B$ (i.e. $B=dA$), having components of class $C^\infty_{\text{\rm pol}}(\mathcal{X})$. Then $P$ defines a self-adjoint operator $\tilde{P}$ on the domain $\mathcal{D}(\tilde{P}):=\boldsymbol{H}^m_A(\mathcal{X})$ and $\mathcal{S}(\mathcal{X})$ is a core for $\tilde{P}$.
\end{theorem}
\begin{proof}
The operator $P$ is symmetric on $\mathcal{S}(\mathcal{X})$, which is dense in $\boldsymbol{H}^m_A(\mathcal{X})$. The case $m=0$ is clear, because $P$ is a bounded operator. 
For $m>0$ we can define $\tilde{P}$ on $\mathcal{D}(\tilde{P}):=\boldsymbol{H}^m_A(\mathcal{X})$ using Proposition \ref{71} and obtain a symmetric operator. If $v\in\mathcal{D}(\tilde{P}^*)$, there exists $f\in L^2(\mathcal{X})$ such that $(\tilde{P}u,v)_{L^2}=(u,f)_{L^2}$ for any $u\in\mathcal{S}(\mathcal{X})$. Thus $Pv=f$ (as distributions) and $v\in\boldsymbol{H}^m_A(\mathcal{X})=\mathcal{D}(P)$, proving self-adjointness.
The last statement follows from the fact that the topology of $\boldsymbol{H}^m_A(\mathcal{X})$ may be defined by the graph norm of $P$ (see Lemma \ref{63}(3)) and from the density of $\mathcal{S}(\mathcal{X})$ in $\boldsymbol{H}^m_A(\mathcal{X})$ (Lemma \ref{69}).
\end{proof}

\noindent
We intend to prove a magnetic version of {\it the G\aa rding inequality}. One needs first

\begin{lemma}\label{57}
Let $0\leq\delta<\rho\leq 1$ and suppose that $f$ is a real valued symbol of class $S^0_{\rho,\delta}(\Xi)$ such that there exists a real valued symbol $f_0\in S^0_{\rho,\delta}(\Xi)$ satisfying $\underset{X\in\Xi}{\inf}f_0(X)>0$ and $f-f_0\in S^{-(\rho-\delta)}_{\rho,\delta}(\Xi)$. 
Then there exists a real valued symbol $g\in S^0_{\rho,\delta}(\Xi)$ such that
$
f\, -\, g\circ^B\!g\,\in\, S^{-\infty}(\Xi)
$.
\end{lemma}
\begin{proof}
Let us define $g_0:=\sqrt{f_0}$ and observe that, due to our hypothesis, it is a real valued symbol of class $S^0_{\rho,\delta}(\Xi)$. Then our hypothesis and Theorem \ref{50} imply that $f-g_0\circ^B\!g_0\,\in\, S^{-(\rho-\delta)}_{\rho,\delta}(\Xi)$. We shall prove by induction that for any $j\in\mathbb{N}$ there exists a real valued symbol $g_j\in S^{-j(\rho-\delta)}_{\rho,\delta}(\Xi)$ such that
$$
r_{j+1}:=f\, -\,\left(\sum_{l=0}^{j}g_l\right)\circ^B\,\left(\sum_{l=0}^{j}g_l\right)
\in S^{-(j+1)(\rho-\delta)}_{\rho,\delta}(\Xi).
$$
The case $j=0$ has just been proved above. Suppose that we have chosen $g_1,\ldots g_{j-1}$ satisfying the stated relations, so that
$
r_{j}:=f-(\sum b_l)\circ^B(\sum g_l)
\in S^{-j(\rho-\delta)}_{\rho,\delta}(\Xi)
$.
Then $g_j$ should be a real valued symbol of class $S^{-j(\rho-\delta)}_{\rho,\delta}(\Xi)$ satisfying the relation
$$
r_{j+1}
=\left\{f\, -\,\left(\sum_{l=0}^{j-1}g_l\right)\circ^B\,\left(\sum_{l=0}^{j-1}g_l\right)
\right\}-
g_j\circ^B\,\left(\sum_{l=0}^{j}g_l\right)-
\left(\sum_{l=0}^{j-1}g_l\right)\circ^B\,g_j=
$$
$$
r_j-g_j\circ^B\!g_0-g_0\circ^B\!g_j+s
\in S^{-(j+1)(\rho-\delta)}_{\rho,\delta}(\Xi),
$$
where $s\in S^{-(j+1)(\rho-\delta)}_{\rho,\delta}(\Xi)$. Using once again Theorem \ref{50}, it is enough to choose $g_j:=r_j/(2g_0)$. Finally, we know that we can find a real valued symbol $g\in S^{0}_{\rho,\delta}(\Xi)$ such that $g\,\sim\,\sum_jg_j$.
\end{proof}

\begin{theorem}\label{77} Let $B$ be a magnetic field with components of class $BC^\infty(\mathcal{X})$.
Let $m\in\mathbb{R}$, $0\leq\delta<\rho\leq 1$, $p\in S^m_{\rho,\delta}(\Xi)$. Suppose that there exist two constants $R$ and $C$ such that $\text{\rm Re}\,p(x,\xi)\geq C|\xi|^m$ for $|\xi|\geq R$. Let us set $P:=\mathfrak{Op}^A(p)$ in any Schr\"{o}dinger representation defined by a vector potential $A$ associated to $B$, whose components are of class $C^\infty_{\text{\sl pol}}(\mathcal{X})$. Then $\forall s\in\mathbb{R}$ there exist two finite positive constants $C_0$ and $C_1$ such that
$$
\text{\rm Re}(Pu,u)_{L^2}\ \geq\ C_0\|u\|^2_{m/2,A}\,-\,C_1\|u\|^2_{s,A},
\quad\forall u\in\boldsymbol{H}^{\infty}_A(\mathcal{X}).
$$
\end{theorem}

\begin{proof}
We assume first that $m=0$. We can choose a positive constant $d$ and a cut-off function $\chi\in C^\infty_0(\mathcal{X}^*)$ such that $\chi(\xi)\geq 0$ and $\chi(\xi)=1$ on a given neighbourhood (large enough) of the origin $0\in\mathcal{X}^*$, so that for $\tilde{p}:=p+d\chi\in S^0_{\rho,\delta}(\Xi)$ we have $\text{\rm Re}(\tilde{p}(x,\xi))\geq c>0$. Hence it is evident that we can take from the begining $\text{\rm Re}\,p(x,\xi)\geq c>0$.
Using Lemma \ref{57}, we deduce the existence of $g\in S^0_{\rho,\delta}(\Xi)_\mathbb R$ and real $r_0\in S^{-\infty}(\Xi)$ such that
$
\text{\rm Re}\,p\,-\,\frac{c}{2}-g\circ^Bg\,=\,r_0
$.
We get
$$
\text{\rm Re}\,(Pu,u)_{L^2}=\frac{c}{2}\|u\|^2_{L^2}+
\left\|\mathfrak{Op}^A(g)u\right\|^2_{L^2}+\left(\mathfrak{Op}^A(r_0)u,u\right)_{L^2},
\;\;\forall u\in\boldsymbol{H}^{\infty}_A(\mathcal{X}),
$$
and thus the inequality for $m=0$.

For the case $m\ne0$, notice that the operator $Q:=\mathfrak{P}_{-m/2}P\mathfrak{P}_{-m/2}$ satisfies the conditions of the case $m=0$. Thus $\forall s^\prime\in\mathbb{R}$
$$
\text{\rm Re}(Qv,v)_{L^2}\,\geq\,C^\prime_0\|v\|^2_{L^2}\,-\,C^\prime_1
\|v\|^2_{s^\prime,A},\quad\forall v\in\boldsymbol{H}^{\infty}_A(\mathcal{X}).
$$
For $u\in\boldsymbol{H}^{\infty}_A(\mathcal{X})$ we denote $v=\mathfrak{Q}_{-m/2}u\in\boldsymbol{H}^{\infty}_A(\mathcal{X})$ where $\mathfrak{Q}_{-m/2}\in\boldsymbol{\Psi}^{A,m/2}(\mathcal{X})$ and $\mathfrak{P}_{-m/2}\mathfrak{Q}_{-m/2}=\boldsymbol{1}+R^\prime$, $\mathfrak{Q}_{-m/2}\mathfrak{P}_{-m/2}=\boldsymbol{1}+R^{\prime\prime}$, with $R^\prime$ and $R^{\prime\prime}$ belonging to $\boldsymbol{\Psi}^{A,-\infty}(\mathcal{X})$. We conclude that
$$
\text{\rm Re}\left(P(\boldsymbol{1}+R^\prime)u,(\boldsymbol{1}+R^\prime)u\right)_{L^2}
\,\geq\,C^\prime_0\left\|\mathfrak{Q}_{-m/2}u\right\|^2_{L^2}\,-\,C^\prime_1
\left\|\mathfrak{Q}_{-m/2}u\right\|^2_{s^\prime,A}.
$$
To obtain the stated inequality, we remark that for functions $f$ and $g$ in $\boldsymbol{H}^{\infty}_A(\mathcal{X})$ and $t\in\mathbb{R}$ we have $|(f,g)_{L^2}|\leq C\|f\|_{t,A}\,\|g\|_{-t,A}$, and that $\forall s\in\mathbb{R}$
$$
\|u\|_{m/2,A}\,\leq\,\left\|\mathfrak{P}_{-m/2}\mathfrak{Q}_{-m/2}u\right\|
_{m/2,A}\,+\,\|R^\prime u\|_{m/2,A}\,\leq
$$
$$
\leq\,C\left(\left\|
\mathfrak{Q}_{-m/2}u\right\|_{L^2}\,+\,\|u\|_{s,A}\right).
$$
\end{proof}

\begin{corollary}\label{79}
Under the hypothesis of Theorem \ref{76}, if $p\geq0$ for $|\xi|\geq R$, the self-adjoint operator $P$ is lower semibounded.
\end{corollary}
\begin{proof}
The single non-trivial case is $m>0$. Taking $s=0$ in Theorem \ref{77} one gets
$
(Pu,u)_{L^2}\,\geq\,C_0\|u\|^2_{m/2,A}\,-\,C_1\|u\|^2_{L^2}\,
\geq\,-C_1\|u\|^2_{L^2}$, $\forall u\in\boldsymbol{H}^{\infty}_A(\mathcal{X})
$.
\end{proof}

\section{Vector potentials with bounded derivatives of strictly positive order}\label{S6}

In this Section we are going to assume (even when it is not stated explicitly) that the magnetic field $B$ can be deduced from a vector potential $A$ satisfying
$$
\vert\left(\partial^\alpha A_j\right)(x)\vert\le C_\alpha,\ \ \forall j=1,\ldots,n,\ \ \forall \alpha\in\mathbb N^n,\ \vert\alpha\vert\ge1,
$$
that implies evidently that all the components of $B$ are of class $BC^\infty$.

\subsection{General facts}

We illustrate our assumption on $B$ by two examples.

\paragraph{Exemple 1}\label{93}

Assume that the components $\{B_{jk}\}_{_{1\leq j<k\leq n}}$ of the magnetic field are $C^\infty$ real valued functions, periodic with respect to a lattice $\Gamma\subset\mathcal{X}$. It has been proved (see \cite{HH} and \cite{If}) that there exists a constant magnetic field $B^\circ=\{B^\circ_{jk}\}_{_{1\leq j<k\leq n}}$ and a potential vector $\tilde{A}$ of class $C^\infty$ and $\Gamma$-periodic, such that $B-B^\circ=d\tilde{A}$. If $A^\circ$ is a linear vector potential defining the magnetic field $B^\circ$, then $A:=A^\circ+\tilde{A}$ has all the derivatives (of strictly positive order) bounded and $B=dA$.

\paragraph{Exemple 2}\label{94}

Let us assume $B_{jk}\in C^\infty\cap L^\infty$, and $\left|(\partial^\alpha B_{jk})(x)\right|\leq\,c_\alpha<x>^{-1}$ for all multiindices $\alpha$ with $|\alpha|\geq1$. We define the associated transversal gauge vector potential
\begin{equation}\label{transv-gauge}
A_j(x)\;:=\;-\sum_{k=1}^n\,x_k\,\int_0^1ds\,s\;B_{jk}(sx).
\end{equation}
Then for any $\alpha$ with $|\alpha|\geq1$ we get
$$
\left(\partial^\alpha A_j\right)(x)\;=\;-\sum_{k=1}^n\,x_k\,\int_0^1ds\,s^{1+|\alpha|}\;\left(\partial^\alpha B_{jk}\right)(sx)\;+\;\{\text{\sl bounded termes}\}.
$$
Outside the ball of radius 1 we have (making the change of variable $s|x|:=t$)
$$
\left|\left(\partial^\alpha A_j\right)(x) \right|\;\leq\;C|x|\int_0^1ds\,s^{1+|\alpha|}\,(1+s|x|)^{-1}\;+\;C_1\;
\leq\;C_2.
$$

\begin{definition}\label{95-D}
For $0\leq\delta<\rho\leq1$ we consider the following metric on $\Xi$:
$$
g^A_{_X}\;\equiv\;g^A_{_{(x,\xi)}}\;:=\;<\xi-A(x)>^{2\delta}|dx|^2\,+
\,<\xi-A(x)>^{-2\rho}|d\xi|^2,
$$
and its symplectic inverse (with respect to the canonical symplectic form $\sp.,.\spp$)
$$
g^{A,\sigma}_{_X}\;\equiv\;g^{A,\sigma}_{_{(x,\xi)}}\;:=\;
<\xi-A(x)>^{2\rho}|dx|^2\,+\,<\xi-A(x)>^{-2\delta}|d\xi|^2.
$$
Let $\mu^A(X):=<\xi-A(x)>\;\;\geq\,1$, $\nu^A(X):=\xi-A(x)\in\mathcal{X}^*$.
\end{definition}

\begin{lemma}\label{95}
The metric defined in Definition \ref{95-D} has the following properties:
\begin{enumerate}
\item[a)] It is a H\"{o}rmander metric, i.e.:
\begin{itemize}
\item {\sl (slow variation)} there exists $C>0$ such that $\;g^A_{_X}(X-Y)\leq C^{-1}$ implies $ (g^A_{_X}/g^A_{_Y})^{\pm 1}\leq C$
\item {\sl (temperedness condition)} there exist $C>0$ and $N\in\mathbb{N}$ such that $(g^A_{_X}/g^A_{_Y})^{\pm 1}\leq C\left(1+g^{A,\sigma}_{_X}(X-Y)\right)^N$ for any $X,Y$ in $\Xi$,
\item {\sl (uncertainty condition)} $g^A_{_X}\,\leq\,g^{A,\sigma}_{_X}$.
\end{itemize}
\item[b)] It is a conformal metric, i.e. $g^{A,\sigma}_{_X}=\lambda^A(X)^2g^A_{_X}$,
 where we have defined $\lambda^A(X)^2:=\inf\left\{g^{A,\sigma}_{_X}(T)\mid g^A_{_X}(T)=1\right\}$.
\item[c)] It is geodesically tempered, i.e. it verifies the temperedness condition with respect to the geodesic distance associated to $g^{A,\sigma}_{_X}$.
\end{enumerate}
\end{lemma}

\begin{proof}
A direct computation gives $g^{A,\sigma}_{_X}=\lambda^A(X)^2g^A_{_X}$, with $\lambda^A(X)=\mu^A(X)^{\rho-\delta}$. We have got the last condition in (a) and condition (b). The first two conditions in (a) follow from \cite{Ho1} once we notice that $\mu^A$ is a 'basic weight function' (see \cite{NU1}), so that the pair $\{\mu^A(X)^\rho,\mu^A(X)^{-\delta}\}$ is a 'pair of weight functions' in the sense of Beals (cf. \cite{Be}).

In order to verify condition (c), let us denote by $d^{A,\sigma}(Y,Z)$ the geodesic distance from $Y=(y,\eta)\in\Xi$ to $Z=(z,\zeta)\in\Xi$ associated to the metric $g^{A,\sigma}$: 
\begin{equation}\label{geod-dist}
L\;\equiv\;d^{A,\sigma}(Y,Z)\;:=\;\underset{X(t)\text{ of class }\;C^1}{\underset{X(0)=Y}{\underset{X(1)=Z}{\inf}}}\,\left\{\int_0^1dt\;
g^{A,\sigma}_{_{X(t)}}(\dot{X}(t))^{1/2}\right\}.
\end{equation}
Thus there exists $\epsilon\in(0,1]$ and a path $X(t)$ such that
$$
\int_0^1dt\left[\mu^A(X(t))^{2\rho}|\dot{x}(t)|^2+\mu^A(X(t))^{-2\delta}|\dot{\xi}(t)|^2 \right]^{1/2}\;=\;L+\epsilon.
$$
Changing the parameter $t\in[0,1]$ with the arc-length $s\in[0,L+\epsilon]$ along the path $X(t)$, we may suppose that $g^\sigma_{_{X(s)}}(\dot{X}(s))=1$ for all $t\in[0,L+\epsilon]$, and we also have $X(0)=Y$ and $X(L+\epsilon)=Z$. Now let us set $\nu^A(s):=\xi(s)-A(x(s))$ and $m^A(s):=<\nu^A(s)>$. We have $\dot{m}^A(s)=m^A(s)^{-1}<\nu^A(s),\dot{\nu}^A(s)>$ and
$$
(1-\delta)^{-1}\left|\mu^A(Z)^{1-\delta}-\mu^A(Y)^{1-\delta}\right|=
\left|\int_0^{L+\epsilon}ds\;m^A(s)^{-\delta}\dot{m}^A(s)\right|\leq
$$
$$
\leq\;C \int_0^{L+\epsilon}ds\;\mu^A(X(s))^{-\delta}\left[|\dot{\xi}(s)|^2+|\dot{x}(s)|^2
\right]^{1/2}\;\leq
$$
$$
\leq\;C \int_0^{L+\epsilon}ds\;\left[\mu^A(X(s))^{-2\delta}|\dot{\xi}(s)|^2+
\mu^A(X(s))^{2\rho}|\dot{x}(s)|^2\right]^{1/2}=\;C(L+\epsilon).
$$
Taking into account that $\mu^A(X)\geq1$, we obtain 
$
\left(\frac{\mu^A(Z)}{\mu^A(Y)}\right)^{\pm 1}\;\leq\;C_1\,(1+L)^N
$
for some suitable constants $C_1>0$ and $N\in\mathbb{N}$. Considering successively the situations $\mu^A(Y)\leq\mu^A(Z)$ and $\mu^A(Z)\leq\mu^A(Y)$ and the above estimation, we finally obtain 
$
\left(\frac{g^A_{_Z}}{g^A_{_Y}}\right)^{\pm 1}\;\leq\;C_0\,(1+d^{A,\sigma}(Y,Z))^{N_0}
$
for some suitable constants $C_0>0$ and $N_0\in\mathbb{N}$.
\end{proof}
\begin{lemma}\label{96}
For any $m\in\mathbb{R}$ the function $M^A_m(X):=<\xi-A(x)>^m$ is a $g^A$-weight for the metric $g^A$ defined in Definition \ref{95-D}, i.e. it is
\begin{itemize}
\item {\sl ($g^A$-continuous)} There exists $C>0$ such that $g^A_{_X}(X-Y)\leq\,C^{-1}$ implies $(M^A_m(X)/M^A_m(Y))^{\pm 1}\leq\,C$,
\item {\sl ($g^A$-tempered)} There exists $C>0$ and $N\in\mathbb{N}$ such that\\ $\left[M^A_m(X)/M^A_m(Y)\right]^{\pm 1}\leq\,C\left[1+g^{A,\sigma}_{_X}(X-Y)\right]^N$.
\end{itemize}
\end{lemma}

\begin{proof}
Suppose $g^A_{_X}(X-Y)\leq\,C^{-1}$; then by the slow variation of $g^A$ we get 
$
\left(\frac{g^A_{_X}(T)}{g^A_{_Y}(T)}\right)^{\pm 1}\leq C$, $\forall T=(t,\tau)\in\Xi
$.
But $g^A_{_X}(0,\tau)/g^A_{_Y}(0,\tau)=(\mu^A(X)/\mu^A(Y))^{-2\rho}$, so that we get the claimed inequality for the function $M^A_m$.
For the second condition, by using the temperedness of the metric $g^A$, we just remark that:
$$
\left(\frac{M^A_m(X)}{M^A_m(Y)}\right)^{\pm 1}\;=\;\left(\frac{\mu^A(X)^{-2\rho}}{\mu^A(Y)^{-2\rho}}\right)^{\mp m/2\rho}\;=\;\left(\frac{g^A_{_X}(0,\tau)}{g^A_{_Y}(0,\tau)}\right)^{\mp m/2}\;\leq
$$
$$
\leq\;C(1+g^{A,\sigma}_{_X}(X-Y))^{N|m|/2}.
$$
\end{proof}

\begin{definition}\label{97}
We consider the following spaces of symbols associated to the metric $g^A$ of Definition \ref{95-D} and a $g^A$-weight $M$:
\begin{itemize}
\item $S^A_{\rho,\delta}(M)\equiv S(M,g^A)$ the symbols $q\in\,C^\infty(\Xi)$ such that
$\forall(\alpha,\beta)\in\mathbb{N}^n\times\mathbb{N}^n$, $\left|(\partial^\alpha_x\partial^\beta_\xi q)(x,\xi)\right|\leq\,C_{\alpha\beta}\,M(X)\mu^A(X)^{-\rho|\beta|+\delta|\alpha|}$.
\item $S^{A,+}_{\rho,\delta}$ the symbols $q\in\,C^\infty(\Xi)$ such that $\forall\alpha,\beta\in\mathbb{N}^n$, with
$|\alpha|+|\beta|\geq1$, we have $\left|(\partial^\alpha_x\partial^\beta_\xi q)(x,\xi)\right|\leq\,C_{\alpha\beta}\,\mu^A(X)^{\rho-\delta-\rho|\beta|+\delta|\alpha|}$,
\item $S^{A,m}_{\rho,\delta}:=S^A_{\rho,\delta}(\mu^m)$ for $m\in\mathbb{R}$ (we call this $m$ the order of the Weyl operator associated to a symbol of this class).
\end{itemize}
If $\rho=1$ and $\delta=0$ the indices $\rho$ and $\delta$ will be omitted from the above notations.
By a slight abuse, for any $p\in C^\infty(\Xi)$ we set $(p\circ\nu^A)(X):=p(x,\nu^A(X))$.
\end{definition}

\begin{remark}\label{98}
For $p\in C^\infty(\Xi)$ it is clear that $p\in S^m_{\rho,\delta}(\Xi)$ if and only if $p\circ\nu^A\in S^{A,m}_{\rho,\delta}$. This allows us to define asymptotic sums of symbols from $S^{A,m}_{\rho,\delta}$. In fact, for a sequence $\{q_j\}_{_{j\in\mathbb{N}}}$ with $q_j\in S^{A,m_j}_{\rho,\delta}$ and $\{m_j\}_{_{j\in\mathbb{N}}}$ decreasing, with $\underset{j\rightarrow\infty}{\lim}m_j=-\infty$, there exists $q\in S^{A,m_0}_{\rho,\delta}$, uniquely defined {\sl modulo} $S^{A,-\infty}:=\underset{m\in\mathbb{R}}{\cap}S^{A,m}_{\rho,\delta}$, such that
$
q-\sum_{j=0}^{k-1}q_j\in S^{A,m_k}_{\rho,\delta}$, $\forall k \geq\,1
$.
We shall write $q\sim\sum_{j=0}^\infty q_j$.
\end{remark} 
\begin{remark}
The symbol $p\in S^m_{\rho,\delta}(\Xi)$ is elliptic if and only if $p\circ\nu^A$ is elliptic for the metric $g^A$ (i.e. $1+|(p\circ\nu^A)(x,\xi)|\geq c[\mu^A(x,\xi)]^m$).
\end{remark}

\subsection{Comparison of two quantizations}

We shall use the notation $\mathfrak{Op}(p)\equiv\mathfrak{Op}^0(p)$ for the usual Weyl quantization. As mentioned in the Introduction, $\mathfrak{Op}_A(p):=\mathfrak{Op}(p\circ\nu^A)$ is sometimes used as the quantization of the symbol $p$. We show now explicitly that it lacks gauge covariance, 
completing the discussion in \cite{MP2}. 

For $f\in S^m(\Xi)$, $A\in C^\infty_{\text{pol}}(\mathcal X,\mathcal X^*)$ and $\varphi\in C^\infty_{\text{pol}}(\mathcal X)$ real valued, set 
$$
F(f,A,\varphi):=e^{i\varphi}\mathfrak{Op}(f\circ\nu^A)e^{-i\varphi}-\mathfrak{Op}(f\circ\nu^{A+\nabla\varphi}).
$$ 
It is an operator with distribution kernel
$$
[K(f,A,\varphi)](x,y)=\exp\left\{i\left<x-y,A\left(\frac{x+y}{2}\right)\right>\right\}\Phi(x,y)\tilde f\left(\frac{x+y}{2},x-y\right),
$$
where $\tilde f$ is the Fourier transform of $f$ in the second variable and 
$$
\Phi(x,y):=\exp\left\{i[\varphi(x)-\varphi(y)]\right\}-\exp\left\{i\left<x-y,(\nabla\varphi)\left(\frac{x+y}{2}\right)\right>\right\}.
$$
Thus gauge covariance is equivalent with the vanishing of the tempered distribution $\exp\{i<y,A(x)>\}\Phi(x+y/2,x-y/2)\tilde f(x,y)$. An easy argument proves that $\phi$ vanishes identically if and only if $\varphi$ is a polynomial of degree $\le 2$. This can easily be used to prove the lack of gauge covariance for a very large class of symbols $f$. Let us consider the monomial $f(x,\xi)=\xi^\alpha$, $\alpha\in\mathbb N$; one has $F(f,A,\varphi)=0$ if and only if $[i(\partial_y+iA(x))]^\beta\left[\Phi(x+y/2,x-y/2)\right]\vert_{y=0}=0$ for any $\beta\le\alpha$. Simple calculations show that this holds if $\vert\beta\vert\le 2$ but is no longer true for $\vert\beta\vert\ge 3$, (one checks easily that $f(x,\xi)=\xi_j\xi_k\xi_l$ is indeed a counterexemple because $(\partial_{y_j}+iA_j(x))(\partial_{y_k}+iA_k(x))(\partial_{y_l}+iA_l(x))
\left[\Phi(x+y/2,x-y/2)\right]\vert_{y=0}=(\partial_j\partial_k\partial_l\varphi)(x)\ne0$ for at least one triple $(j,k,l)$ if $\varphi$ is not a second order polynomial). Let us also notice that the Fourier transform of $f(x,\xi)=<\xi>$ is a distribution $\tilde f$ with singular support $\mathcal{X}\times\{0\}$ and analytic on its complement. In fact $\tilde{f}$ is rotation invariant and some straightford computation proves that it verifies an ordinary differential equation (in the radial variable) with coefficients analytic outside $\{0\}$. 
Thus it is nonzero on a dense set in $\mathcal{X}\times\mathcal{X}$ and in order to have gauge covariance, the function $\Phi$ should be identically zero, but this is not the case if $\varphi$ is not a second order polynomial. We conclude that $\mathfrak{Op}_A(<\xi>)$ dos not provide a gauge covariant quantization.

In spite of all these, it is useful to express $\mathfrak{Op}^A(p)$ as $\mathfrak{Op}_A(q)$ for some symbol $q$, but keeping in mind that this operator is the magnetic quantization of $p$ and not of $q$. We are going to explore this in the sequel.
We define $\Gamma^A([x,y]):=<x-y,\Gamma^A(x,y)>$ so that $\Gamma^A(x,y)=\int_0^1ds\,A((1-s)x+sy)$.
\begin{proposition}\label{100}
For any $p\in S^m_{\rho,\delta}(\Xi)$ there exists a unique $q\in S^m_{\rho,\delta}(\Xi)$ such that
$
\mathfrak{Op}^A(p)=\mathfrak{Op}(q\circ\nu^A)
$.
Besides, we have
$
q\circ\nu^A\sim\sum_{j=0}^\infty q^A_j
$,
where
$$
q^A_j(X):=\sum_{|\alpha|=j}\left.\frac{1}{\alpha!}\left\{(-D_y)^\alpha
\partial^\alpha_\xi\left[p(x,\xi-\Gamma^A(x+y/2,x-y/2)\right]\right\}\right|_{_{y=0}}.
$$
In particular $q^A_0=p\circ\nu^A$, $q^A_1=0$, $q^A_j\in S^{A, m-(j+1)\rho}_{\rho,\delta}$ ($\forall j\geq1$), $q-p\in S^{m-3\rho}_{\rho,\delta}(\Xi)$.
\end{proposition}
\begin{proof}
Let us write down the distribution kernel of the operator $\mathfrak{Op}^A(p)$ (using oscillatory integrals with values in $\mathcal{S}^*(\mathcal{X}\times\mathcal{X})$):
$$
K^A(x,y)\;=\;\int_{\mathcal{X}^*}\dbar\eta\;e^{i<x-y,\eta>}\,p\left(\frac{x+y}{2},\eta-\Gamma^A(x,y)\right).
$$
Thus the usual Weyl symbol of the operator $\mathfrak{Op}^A(p)$ is
$$
\tilde{q}^A(X)=\int_{\mathcal{X}}dy\;e^{-i<y,\xi>}\,K^A\left(x+(y/2),x-
(y/2)\right)
=
$$
$$
=\int\!\!\!\int_{\mathcal{X}\times\mathcal{X}^*}dy\dbar\eta\;e^{i<y,\eta>}\,p\left(
x,\xi+\eta-\Gamma^A(x+(y/2),x-(y/2))\right).
$$
We use the Taylor expansion
$
p(x,\zeta+\eta)=\sum\limits_{|\alpha|<N}(\alpha!)^{-1}\eta^\alpha\,
(\partial^\alpha_\eta p)(\zeta)+r_N(x,\zeta,\eta)
$,
$$
r_N(x,\zeta,\eta)\;:=\;\sum\limits_{|\alpha|=N}\frac{\eta^\alpha}{(N-1)!}
\int_0^1dt\,(1-t)^{N-1}\;(\partial^\alpha_\eta p)(x,\zeta+t\eta).
$$
Inserting this development in the definition of $\tilde{q}^A$, we get
$$
\tilde{q}^A(X)\;=\;\sum\limits_{j=0}^{N-1}q^A_j(X)\;+\;\tilde{r}^A_N(X),
$$
with $q^A_j$ given exactly by the formula in the statement of Proposition \ref{100}. For $\tilde{r}^A_N$, we get the explicit formula
$$
\tilde{r}^A_N(X)\;=\;\frac{1}{(N-1)!}\sum\limits_{|\alpha|=N}
\int_0^1dt\,(1-t)^{N-1}\int\!\!\!\int_{\mathcal{X}\times\mathcal{X}^*}dy
\dbar\eta\;e^{i<y,\eta>}\times
$$
$$
\times\left\{(-D_y)^\alpha\partial^\alpha_\xi\left[p(x,\xi+t\eta-
\Gamma^A(x+y/2,x-y/2)\right]\right\}.
$$
It is clear now that $q^A_0=p\circ\nu^A$ and $q_j\in S^{A,m-(j+1)\rho}_{\rho,\delta},\ \forall j\geq1$. Moreover
$$
q^A_1(X)\;=\;\frac{i}{2}\sum\limits_{j,k=1}^n(\partial_{\xi_j}\partial_{\xi_k}p)(\xi-A(x))\left(\frac{\partial\Gamma^A_k}{\partial y_j}-\frac{\partial\Gamma^A_k}{\partial x_j}\right)(x,x),
$$
$$
\left(\frac{\partial\Gamma^A_k}{\partial y_j}-\frac{\partial\Gamma^A_k}{\partial x_j}\right)(x,x)\;=\;\int_0^1ds\,(2s-1)\,(\partial_jA_k)(x)\;=\;0,
$$
and we get $q^A_1(X)=0$.
We shall estimate now the derivatives of $\tilde{r}^A_N$, by some integration by parts, using the identities
$$
<y>^{-2N}(1-\Delta_\eta)^N\,e^{i<y,\eta>}=e^{i<y,\eta>},
$$
$$
<\eta>^{-2N}(1-\Delta_y)^N\,e^{i<y,\eta>}=e^{i<y,\eta>}.
$$
For $X\in\Xi$ fixed and any $\epsilon>0$, we decompose the integral over $\eta\in\mathcal{X}^*$ into two parts corresponding to the two regions
$$
\mathcal{D}_\epsilon:=\left\{\eta\in\mathcal{X}^*\mid <\eta>\,\leq\,\epsilon\,
<\xi-\Gamma^A(x+y/2,x-y/2)>\right\},\;\;\widetilde{\mathcal{D}_\epsilon}:=
\mathcal{X}^*\setminus\mathcal{D}_\epsilon.
$$
But
$
\;<A(x)-\Gamma^A(x+y/2,x-y/2)>\;\leq\;C\int_{-1/2}^{1/2}ds\,<A(x)-A(x-sy)>$
is bounded by $C<y>$
and thus
$
<\xi-\Gamma^A(x+y/2,x-y/2)>^m$ is bounded by $C\,<\xi-A(x)>^m<y>^m$, $\forall m\in\mathbb{R}
$.
Then, for $\alpha,\beta$ in $\mathbb{N}^n$, we have $|(\partial_x^\alpha\partial_\xi^\beta\tilde{r}^A_N)(X)|\leq C_{\alpha\beta}(R^\prime_N(X)+R^{\prime\prime}_N(X))$, where $R^\prime_N(X)$ is by definition
$$
\int_0^1dt\int_{\mathcal{X}}dy\int_{\mathcal{D}_\epsilon}\dbar\eta\,<y>^{-2N_1}
<\eta>^{-2N_2}<\xi+t\eta-\Gamma(x+y/2,x-y/2)>^k,
$$
$R^{\prime\prime}_N(X)$ is the same integral but on the domain $\widetilde{\mathcal{D}_\epsilon}$ instead of $\mathcal{D}_\epsilon$  and  $k=m-(\rho-\delta)N-\rho|\beta|+\delta|\alpha|+2\delta N_2\;=\;k_0+2\delta N_2$. $\mathcal{D}_\epsilon$ being a bounded domain in $\mathcal{X}^*$, we shall take $N_2=0$ and obtain
$$
R^\prime_N(X)\leq\;C\,
<\xi-A(x)>^{k_0+n}\int_{\mathcal{X}}dy\,<y>^{-2N_1+|k_0|+n}
$$
that is finite if we choose $2N_1\,\geq\,|k_0|+2n+1$. Then for 
$R^{\prime\prime}_N(X)$ we get the bound
$$
C\,\int_{\mathcal{X}}dy\,<y>^{-2N_1}\int_{\widetilde{\mathcal{D}_\epsilon}}
\dbar\eta\,<\eta>^{|k_0|-2(1-\delta)N_2}
$$
$$
\leq\;C\,<\xi-A(x)>^{|k_0|+n-2(1-\delta)N_2}\int_{\mathcal{X}}dy\,
<y>^{-2N_1+||k_0|+n-2(1-\delta)N_2|}.
$$
We have to choose
$
2N_1>\max\{|k_0|+2n+1,||k_0|+n-2(1-\delta)N_2|+n+1\}
$ and also
$
|k_0|+n-2(1-\delta)N_2\leq k_0+n
$
in order to obtain  the conclusion $\tilde{r}^A_N\in S^{A,m-(\rho-\delta)N+n}_{\rho,\delta}$. We apply the argument at the end of the proof of Theorem \ref{50}, to get that $\tilde{q}^A\in S^{A,m}_{\rho,\delta}$.
We end the proof by taking $q(X):=\tilde{q}^A(x,\xi+A(x))$.
\end{proof}

\begin{remark}\label{101}
If $p\in S^2(\Xi)$ is a polynomial of degree less then or equal to 2 in the variable $\xi\in\mathcal{X}^*$ (with coefficients depending on the variable $x\in\mathcal{X}$), we have $q=p$ in the above Proposition and thus we get $\mathfrak{Op}^A(p)=\mathfrak{Op}(p\circ\nu^A)$.
\end{remark}

\begin{proposition}\label{102} {\sl (Converse of Proposition \ref{100})}
For any $q\in S^m_{\rho,\delta}(\Xi)$, there exists a unique $p\in S^m_{\rho,\delta}(\Xi)$ such that $\mathfrak{Op}^A(p)=\mathfrak{Op}(q\circ\nu^A)$.
\end{proposition}
\begin{proof}
For any tempered distribution $p\in\mathcal{S}^*(\Xi)$ we can define the operator $\mathfrak{Op}^A(p)\in\mathcal{B}(\mathcal{S}(\mathcal{X}),\mathcal{S}^*(\mathcal{X}))$, as an integral operator with distribution kernel
$$
K^A(x,y)\;=\;\int_{\mathcal{X}^*}\dbar\eta\;e^{i<x-y,\eta>}\,e^{i<x-y,\Gamma^A(x,y)>}\,p\left(\frac{x+y}{2},\eta\right),
$$
so that
$$
\int_{\mathcal{X}^*}\dbar\eta\;e^{i<y,\eta>}\,p\left(x,\eta\right)\,=\,e^{-i<y,\Gamma^A(x+y/2,x-y/2)>}K^A(x+y/2,x-y/2),
$$
and finally
$$
p(x,\xi)\;=\;\int_{\mathcal{X}}dy\;e^{-i<y,\xi>}e^{-i<y,\Gamma(x+y/2,x-y/2)>}K(x+y/2,x-y/2).
$$
One can also write the distribution kernel of the Weyl operator $\mathfrak{Op}(q\circ\nu)$
$$
K(x,y)\;=\;\int_{\mathcal{X}^*}\dbar\eta\;e^{i<x-y,\eta>}\,q\left(\frac{x+y}{2},
\eta-A\left(\frac{x+y}{2}\right)\right).
$$
We may have the equality $\mathfrak{Op}^A(p)=\mathfrak{Op}(q\circ\nu^A)$ if and only if
$$
p(X)=\int\!\!\!\int_{\mathcal{X}\times\mathcal{X}^*}dy\dbar\eta\;e^{i<y,\eta-\xi>}
e^{-i<y,\Gamma^A(x+y/2,x-y/2)>}q(x,\eta-A(x))\;=
$$
$$
=\;\int\!\!\!\int_{\mathcal{X}\times\mathcal{X}^*}dy\dbar\eta\;e^{i<y,\eta>}q\,
(x,\xi+\eta+\Gamma^A(x+y/2,x-y/2)-A(x)).
$$
Proceeding then as in the proof of Proposition \ref{100} we show that $p\in S^m_{\rho,\delta}(\Xi)$.
\end{proof}

\begin{remark}
The Propositions \ref{100} and \ref{102} imply that, under the hypothesis of Section \ref{S6}, the properties of the magnetic pseudodifferential operators may be obtained through the usual Weyl functional calculus associated to the metric $g^A$ ( \cite{Bo3}, \cite{Ho1}, \cite{Ho2}). An exemple is the following Fefferman-Phong theorem:
\end{remark}
\begin{corollary}\label{103}
Let us choose $p\in S^{2(\rho-\delta)}_{\rho,\delta}(\Xi)$ with $p\geq0$. Then there exists a constant $C>0$ such that
$
\left(\mathfrak{Op}^A(p)u,u\right)_{L^2}\;\geq\;-C\|u\|^2_{L^2}$, $\forall u\in\mathcal{S}(\mathcal X).
$
\end{corollary}
\begin{proof}
Choosing $p$ as in the statement of the Corollary and using Proposition \ref{100}, we conclude that there exist $q\in S^{2(\rho-\delta)}_{\rho,\delta}(\Xi)$ and $r\in S^{-\rho-2\delta}_{\rho,\delta}(\Xi)$ such that $q=p+r$ and $\mathfrak{Op}^A(p)=\mathfrak{Op}(q\circ\nu)$. The condition $p\geq0$ implies that $\mathfrak{Op}^A(p)$ is symetric and thus $q$ and $r$ will be real. Thus we can write
$
\mathfrak{Op}^A(p)=\mathfrak{Op}(p\circ\nu^A)+\mathfrak{Op}(r\circ\nu^A)
$
and $p\circ\nu^A\in S^{A,2(\rho-\delta)}_{\rho,\delta}$. As a consequence of the Fefferman-Phong inequality (\cite{Ho2}, T.18.6.8), there exists a constant $C_0>0$ such that
$
\left(\mathfrak{Op}(p\circ\nu^A)u,u\right)_{L^2}\geq -C_0\,\|u\|^2_{L^2}$, $\forall u\in\mathcal{S}(\mathcal{X})
$.
Using Proposition \ref{102} we deduce the existence of a symbol $r_0\in S^{-\rho-2\delta}_{\rho,\delta}(\Xi)$ such that $\mathfrak{Op}(r\circ\nu^A)=\mathfrak{Op}^A(r_0)$, that is a bounded operator in $L^2(\mathcal{X})$ due to the fact that $\rho+2\delta\geq0$ and to Remark \ref{3.8}. We conclude that there exists a constant $C_1>0$ such that
$
\left(\mathfrak{Op}(r\circ\nu^A)u,u\right)_{L^2}\;\geq\;-C_1\,\|u\|^2_{L^2}$, $\forall u\in\mathcal{S}(\mathcal{X})
$.
\end{proof}

\subsection{Resolvents and fractional powers of elliptic magnetic pseudodifferential operators}

Due to the fact that the H\"{o}rmander metric $g^A$ is conformal and geodesically temperate we can use a Theorem of Bony (\cite{Bo1}) characterizing pseudodifferential operators by commutators and prove that the resolvent and the powers of an elliptic magnetic self-adjoint pseudodifferential operator are also of this type.

\begin{theorem}\label{104} {\sl (Bony)}
Let $q\in\mathcal{S}^*(\Xi)$ and $Q:=\mathfrak{Op}(q)$. Then $q\in S^{A,m}_{\rho,\delta}$ if and only if $Q\in\mathcal{B}(\boldsymbol{H}^m_A,L^2(\mathcal{X}))$ and for any finite family $\{b_j\}_{_{1\leq j\leq k}}\subset S^{A,+}_{\rho,\delta}$ we have $\text{\sl ad}(\mathfrak{Op}(b_1))\cdots\text{\sl ad}(\mathfrak{Op}(b_k))Q\in\mathcal{B}(\boldsymbol{H}^m_A,L^2(\mathcal{X}))$.
\end{theorem}

\begin{corollary}\label{105}
Under the hypothesis of Theorem \ref{76} let $P:=\mathfrak{Op}^A(p)$. We also denote by $P$ the induced self-adjoint operator in $L^2(\mathcal{X})$ (with domain $\boldsymbol{H}^m_A$). Then for any $z\in\mathbb{C}\setminus\sigma(P)$ we have $(P-z)^{-1}=\mathfrak{Op}^A(\tilde{p}_z)$ with $\tilde{p}_z\in S^{-m}_{\rho,\delta}(\Xi)$.
\end{corollary}

\begin{proof}
Obviously $(P-z)^{-1}\in\mathcal{B}(\boldsymbol{H}^{-m}_A,L^2(\mathcal{X}))$. Using Proposition \ref{100}, there exists $q\in S^{A,m}_{\rho,\delta}$ such that $\mathfrak{Op}^A(p)=\mathfrak{Op}(q)$. For a finite family $\{b_j\}_{_{1\leq j\leq k}}\subset S^{A,+}_{\rho,\delta}$, the arguments in \cite{Bo1} imply that $\text{\sl ad}(\mathfrak{Op}(b_1))\cdots\text{\sl ad}(\mathfrak{Op}(b_k))\mathfrak{Op}(q)$ is a Weyl operator having a symbol of class $S^{A,m}_{\rho,\delta}$. A simple computation shows that the operator $\text{\sl ad}(\mathfrak{Op}(b_1))\cdots\text{\sl ad}(\mathfrak{Op}(b_k))(P-z)^{-1}$ is a finite sum of terms of the form 
$
\pm(P-z)^{-1}K_{a_1}(P-z)^{-1}\cdots(P-z)^{-1}K_{a_l}(P-z)^{-1},
$
where $l\le k$ and each factor $K_a$ is of the form 
$$
K_a\;=\;\underset{j\in J_a}{\Pi}\text{\sl ad}(\mathfrak{Op}(b_j))P\;=\;\underset{j\in J_a}{\Pi}\text{\sl ad}(\mathfrak{Op}(b_j))\mathfrak{Op}(q),
$$
with $J_a$ finite subset of $\{1,\ldots,k\}$. We conclude that $K_a\in\mathcal{B}(L^2(\mathcal{X}),\boldsymbol{H}^{-m}_A)$, and thus 
$$
\text{\sl ad}(\mathfrak{Op}(b_1))\cdots\text{\sl ad}(\mathfrak{Op}(b_k))(P-z)^{-1}\in\mathcal{B}(\boldsymbol{H}^{-m}_A,L^2(\mathcal{X})).
$$
Using Theorem \ref{104}, we conclude the existence of a symbol $\tilde{q}_z\in S^{A,-m}_{\rho,\delta}$ such that $(P-z)^{-1}=\mathfrak{Op}(\tilde{q}_z)$. By Proposition \ref{102}, we deduce the existence of a symbol $\tilde{p}_z\in S^{-m}_{\rho,\delta}(\Xi)$ such that
$
\mathfrak{Op}^A(\tilde{p}_z)=\mathfrak{Op}(\tilde{q}_z)=(P-z)^{-1}
$.
\end{proof}

\begin{remark}\label{106}
From Theorem \ref{104} it follows directly that an operator $\mathfrak{Op}^A(p)$ (with $p\in\mathcal{S}^*(\Xi)$) is a ``smoothing'' one, i.e. transforms $\boldsymbol{H}^{-\infty}_A$ into $\boldsymbol{H}^{\infty}_A$, if and only if it belongs to $\boldsymbol{\Psi}^{A,-\infty}(\Xi)$.
\end{remark}

We use now Corollary \ref{105} and some ideas from the proof of Proposition 29.1.9 in \cite{Ho3} in order to study the fractional powers of an operator as in Corollary \ref{105}. We first remark from Corollary \ref{79} that (for the case $n\geq2$ and replacing if necessary $p$ by $-p$) $\mathfrak{Op}^A(p)$ is lower semibounded. Thus in this case (adding if necessary a sufficiently large constant) we may suppose that $p\geq1$ and $\mathfrak{Op}^A(p)\geq 1$. We can work with the usual Weyl quantization, because (having assumed that the magnetic field $B$ admits a vector potential with bounded derivatives of any strictly positive order) the two quantization are in a one-to-one correspondence that associates to elliptic magnetic operators, operators that are elliptic with respect to the metric $g^A$.

Given $p\in S^{m}_{\rho,\delta}(\Xi)$, resp. $p\in S^{A, m}_{\rho,\delta}$, we call a {\it principal symbol} of $\mathfrak{Op}^A(p)$, resp. $\mathfrak{Op}(p)$, any element $p_0\in S^{m}_{\rho,\delta}(\Xi)$, resp. $p_0\in S^{A,m}_{\rho,\delta}$, satisfying $p-p_0\in S^{m-(\rho-\delta)}_{\rho,\delta}(\Xi)$, resp. $p-p_0\in S^{A,m-(\rho-\delta)}_{\rho,\delta}$.

\begin{theorem}\label{107}
Let $m>0$, $p\in S^{A,m}_{\rho,\delta}$ a real elliptic symbol, such that $p\geq1$ and $P:=\mathfrak{Op}(p)\geq1$. Then for any $s\in\mathbb{R}$ we have 
$
P^s=\mathfrak{Op}(q_s)
$
for some $q_s\in S^{A,sm}_{\rho,\delta}$. Moreover $P^s$ admits $p^s$ as principal symbol.
\end{theorem}
\begin{proof}
For $s\in\mathbb{N}$ the Theorem results directly from the Weyl calculus. Corollary \ref{105} implies the case $s=-1$ and thus we conclude that the Theorem is true for any $s\in\mathbb{Z}$.

Taking into account the composition of symbols, we only have to prove the case $-1<s<0$. We consider the Cauchy formula for the function $f(z):=z^s$ on the domain
$
\left\{\,z\in\mathbb{C}\,\mid\,\Re(z)>0,\;0<\epsilon<|z|<R\,\right\}
$.
Letting $\epsilon\rightarrow0$ and $R\rightarrow\infty$ we get for any $u\in L^2(\mathcal{X})$
$$
P^su\;=\;-(2\pi i)^{-1}\int_{-i\infty}^{i\infty}dz\;z^s\left(P-z\right)^{-1}u.
$$
Then we notice that $p-z\,\in\,S^A_{\rho,\delta}(\mu^m+|z|)$ and $(p-z)^{-1}\,\in\,S^A_{\rho,\delta}((\mu^m+|z|)^{-1})$, uniformly for $z\in\,i\mathbb{R}$. Denoting by $\circ$ the usual Weyl composition of symbols, we get
$
(p-z)\circ(p-z)^{-1}=1-r_z
$
for some $r_z\,\in\,S^A_{\rho,\delta}(\mu^{-1})$. Let us denote by $q^{\circ j}$ the j-th power of $q$ with respect to the product $\circ$ and let $e_z\sim\sum_{j=0}^\infty r_z^{\circ j}$. Then $e_z\in S^A_{\rho,\delta}$, $e_z-1\in S^A_{\rho,\delta}(\mu^{-1})$ and
$
g_z:=(1-r_z)\circ e_z-1\in S^{A,-\infty}
$.
Moreover
$
h_z:=(p-z)^{-1}\circ e_z\in S^A_{\rho,\delta}((\mu^m+|z|)^{-1})
$ and
$
(p-z)\circ h_z=1+g_z
$.
We conclude that
$
(P-z)^{-1}=\mathfrak{Op}(h_z)-(P-z)^{-1}\mathfrak{Op}(g_z)=
$
$
\mathfrak{Op}\left[(p-z)^{-1}\right]+\mathfrak{Op}(\alpha_z)-(P-z)^{-1}\mathfrak{Op}(g_z),
$
where $\alpha_z:=(p-z)^{-1}\,\circ\,f_z\;\in\;S_{\rho,\delta}\left(\mu^{-1}(\mu^m+|z|)^{-1}\right)$. It follows that 
$$
P^s=\mathfrak{Op}(b_0)+\mathfrak{Op}(b_1)+R,
$$ 
where
$$
b_0:=-(2\pi i)^{-1}\int_{-i\infty}^{i\infty}dz\;z^s(p-z)^{-1}=p^s\;\in\;S_{\rho,\delta}^{A,sm},
$$
$$
b_1:=-(2\pi i)^{-1}\int_{-i\infty}^{i\infty}dz\;z^s\alpha_z,
$$
$$
R:=(2\pi i)^{-1}\int_{-i\infty}^{i\infty}dz\;z^s\left(P-z\right)^{-1}\mathfrak{Op}(g_z).
$$
Taking $z=it$ with $t\in\mathbb{R}$ and recalling that  $\alpha_z\in S^A_{\rho,\delta}(\mu^{-1}(\mu^m+|z|)^{-1})$ uniformly in $t\in\mathbb{R}$, we deduce that for any $\beta$, $\gamma$ in $\mathbb{N}^n$ the derivatives $\left(\partial_x^\beta\partial_\xi^\gamma b_1\right)(x,\xi)$ are bounded by
$$
C_{\beta\gamma}\,\int_{-\infty}^{\infty}dt\;|t|^s(\mu^m+|t|)^{-1}\,
\mu^{-\rho|\gamma|+\delta|\beta|-1}\;\leq\;C_{\beta\gamma}\,
\mu^{sm-1-\rho|\gamma|+\delta|\beta|},
$$
so that $b_1\,\in\,S_{\rho,\delta}^{A,sm-1}$. Recalling Theorem \ref{104} and Remark \ref{106}, we finish the proof by showing that the operator $R$ is ``smoothing''.

First we notice that for any $k\in\mathbb{Z}$, the norms $\|u\|_{\boldsymbol{H}_A^{km}}$ and $\|P^{k}u\|_{L^2}$ are equivalent on $\boldsymbol{H}_A^{km}$. Thus there exists a constant $C>0$ such that for any $t\in\mathbb{R}$, taking $z=it$ and $u\in\boldsymbol{H}_A^{km}$ we have
$$ 
\left\|(P-z)^{-1}u\right\|_{\boldsymbol{H}_A^{km}}
\leq\,C\,\left\|P^{k}(P-z)^{-1}
u\right\|_{L^2}\,\leq
$$
$$
\leq\,C\,(t^2+1)^{-1/2}\left\|P^{k}u\right\|_{L^2}\,\leq\,
C_1\,(t^2+1)^{-1/2}\|u\|_{\boldsymbol{H}_A^{km}}.
$$
Thus $\forall\,k\in\mathbb{Z}$ and $\forall\,t\in\mathbb{R}$, 
$
\left\|(P-z)^{-1}\right\|_{\mathcal{B}
(\boldsymbol{H}_A^{km})}\,\leq\,C_k\,(t^2+1)^{-1/2}
$.
As it is easy to see that the operator $\mathfrak{Op}(g_z)$ is ``smoothing'' uniformly in $z=it$ with $t\in\mathbb{R}$, we conclude that for any $k$ and $l$ in $\mathbb{Z}$ and any $u\in\mathcal{S}(\mathcal{X})$ we have
$$
\|Ru\|_{\boldsymbol{H}_A^{km}}\;\leq\;C(k,l)\,\|u\|_{\boldsymbol{H}^l_A}\,\int_{-\infty}^\infty dt\,|t|^s(1+|t|)^{-1}\;\leq\;C^\prime(k,l)\,\|u\|_{\boldsymbol{H}^l_A}.
$$
\end{proof}
\begin{remark}
Using Theorem \ref{107} and Proposition \ref{100} we see that the operators  $\mathfrak{Op}^A(<\xi>)$, $\mathfrak{Op}_A(<\xi>)=\mathfrak{Op}(\mu^A)$ and $\sqrt{(D-A)^2+1}$, are elliptic Weyl pseudodifferential operators of first order associated to the metric $g^A$ and having the same principal symbol $\mu^A$. Thus, all three define self-adjoint, lower semibounded operators in $L^2(\mathcal{X})$, having the same domain $\boldsymbol{H}^1_A(\mathcal{X})$ and differing only by bounded $L^2$ operators. Each one may be a candidate for a magnetic relativistic Schr\"{o}dinger Hamiltonian. Nevertheless, the last one cannot be obtained by a complete 'quantization' procedure applying to a larger class of classical observables, while the second one (although used in \cite{Ic1}, \cite{Ic2}, \cite{IT1}, \cite{IT2}, \cite{ITs1}, \cite{ITs2}, \cite{NU1}, \cite{NU2}, etc.) is not covariant for the gauge transformations. Thus, we consider that the only adequate one should be the first one.
\end{remark}

\section{The limiting absorption principle}

This section is devoted to the spectral analysis of operators of the form $\mathfrak{Op}(p)$, $\mathfrak{Op}_A(p)\equiv\mathfrak{Op}(p\circ\nu^A)$ and $\mathfrak{Op}^A(p)$, for an elliptic symbol $p\in S^m(\Xi)$, and a limiting absorption principle for this type of operators is obtained. 
The main tool we shall use is an abstract result belonging to the "conjugate operator method", (proved in \cite{ABG}). We shall also make use of some known properties of the Weyl calculus (\cite{Ho1}, \cite{Ho2}) and of the magnetic pseudodifferential calculus developed above.
The following hypothesis will be assumed all over this section:
\begin{hypothesis}\label{I}
There exists $\epsilon>0$ (that we can always suppose smaller then $n-1$), such that for any $\alpha\in\mathbb{N}^n$ there exists $C_\alpha>0$ for which $|(\partial^\alpha B_{jk})(x)|\leq C_\alpha<x>^{-1-\epsilon}$, for any $x\in\mathcal{X}$ and $j,k\in \{1,\ldots,n\}$.
\end{hypothesis}

Concerning the vector potential $A$ defining $B$, we shall suppose that it has been chosen to satisfy the conditions in Lemma \ref{80} below.

\begin{lemma}\label{80}
Suppose that Hypothesis \ref{I} is satisfied. Then:
\begin{itemize}
\item[a)]
there exists a vector potential $A$ such that $B=dA$ and for any multiindex $\alpha\in\mathbb{N}^n$, $|(\partial^\alpha A_{j})(x)|\leq C_\alpha<x>^{-\epsilon}$ for any $x\in\mathcal{X}$ and any $j$ in $\{1,\ldots,n\}$,
\item[b)]
if $|\alpha|\geq1$, then the above vector potential $A$ also satisfies
$
|(\partial^\alpha A_{j})(x)|\leq C_\alpha<x>^{-1-\epsilon}\ln(1+<x>)
$.
\end{itemize}
\end{lemma}
\begin{proof}
We choose the Coulomb gauge
$$
A_j(x):=\sum_{k=1}^n\int_\mathcal{X}dy\,(\partial_kE)(y)B_{kj}(x-y),
$$
where $E$ is the standard elementary solution of the Laplace operator on $\mathcal{X}$. 

a) It is evident by this definition that $A_j\in C^\infty(\mathcal{X})$ and $dA=B$. Moreover, for any $\alpha\in\mathbb{N}^n$,
$
\left|(\partial^\alpha A_j)(x)\right|\,\leq\,C\int_\mathcal{X}dy\,|y|^{1-n}(1+|x-y|)^{-1-\epsilon}
$.
For $|x|\leq1$ the above integral is bounded. For $|x|\geq1$ we make the change of variables $y=|x|\tilde{y}$ and get
$$
\left|(\partial^\alpha A_j)(x)\right|\leq\,C|x|^{-\epsilon}\int_\mathcal{X}
d\tilde{y}\,|\tilde{y}|^{1-n}\left|\frac{x}{|x|}-\tilde{y}
\right|^{-1-\epsilon},
$$
this last integral being uniformly bounded with respect to $|x|\geq1$.

b) Obviously it is enough to prove the estimation for the first-order derivatives, and to consider only the terms where these derivatives are acting on the elementary solution $E$. Using the results in \S 8 of Chapter II in \cite{Mi} we have
$$
\partial_lA_j(x)\,=\,\sum_{k=1}^n\int_\mathcal{X}dy\,(\partial_l\partial_kE)(x-y)\,B_{kj}(y)\;+\;n^{-1}B_{lj}(x),
$$
where the integral has to be interpreted as a {\it principal value}. We can write:
$$
B_{kj}(x)-B_{kj}(y)\,=\,\int_0^1dt\,<x-y\,,\,\left(\nabla B_{kj}\right)(x+t(y-x))>,
$$
so that
$$
|B_{kj}(x)-B_{kj}(y)|\leq\,C|x-y|<x>^{-1-\epsilon}<x-y>^{1+\epsilon}.
$$
For $|x-y|\leq1$ we obtain
$
|B_{kj}(x)-B_{kj}(y)|\,\leq\,C_1<x>^{-1-\epsilon}|x-y|
$.
We use then Theorem 1.7 from \S 7 of Chapter II in \cite{Mi} to estimate the singular integral and get
$$
\left|\int_\mathcal{X}dy\,(\partial_l\partial_kE)(x-y)\,B_{kj}(y) \right|\,
\leq\,C<x>^{-1-\epsilon}\ln(1+<x>).
$$
\end{proof}
Before formulating the main result of this section (Theorem \ref{89}) let us make some remarks. Any vector potential verifying the conditions in the above Lemma \ref{80} has the property $\partial^\alpha A_j\in L^\infty(\mathcal{X})$ for any $\alpha\in\mathbb{N}^n$ and any $1\leq j\leq n$. Thus we can apply the results of our previous Section \ref{S6}. Moreover it is easy to verify that in the present situation, all the magnetic Sobolev spaces $\boldsymbol{H}^s_A(\mathcal{X})$ (defined in Section \ref{S4}) coincide with the usual Sobolev spaces $H^s(\mathcal{X})\equiv \boldsymbol{H}^s_0(\mathcal{X})$.

If $p\in S^m(\Xi)$ is a real elliptic symbol and $m>0$, the operator $P:=\mathfrak{Op}(p)$ is self-adjoint in $L^2(\mathcal{X})$, having the domain $H^m(\mathcal{X})$. We shall denote its form domain by $\mathcal{G}:=\mathcal{D}(|P|^{1/2})=H^{m/2}(\mathcal{X})$. Let us still denote by $\mathcal{G}_{s,p}$ and $\mathcal{G}_{s,p}^*$ ($s\in\mathbb{R}$ and $1\leq p\leq\infty$) the spaces of the Besov scale associated to $\mathcal{G}$ and $\mathcal{G}^*\equiv H^{-m/2}(\mathcal{X})$ (see \cite{ABG}). Let us finally remark that for any $z\in\mathbb{C}_{\pm}$ we have $(P-z)^{-1}\in\mathcal{B}(\mathcal{G}^*;\mathcal{G})\subset\mathcal{B}
(\mathcal{G}_{1/2,1}^*;\mathcal{G}_{-1/2,\infty})$.

We shall denote by $g$ the metric $g_{_X}:=|dx|^2+<\xi>^{-2}|d\xi|^2$ and by $M_{m,\delta}$ (for $m$ and $\delta$ in $\mathbb{R}$) the weight function $M_{m,\delta}(X):=<x>^{-\delta}<\xi>^m$, for $X=(x,\xi)\in\Xi$.

\begin{theorem}\label{89}
Assume that the magnetic field $B$ satisfies Hypothesis \ref{I}. Let $p\in S^m(\Xi)$, with $m>0$, satisfying the conditions:
\begin{enumerate}
\item[i)] $p$ is real valued and elliptic;
\item[ii)] there exists $p_0\in S^m(\Xi)$ a real elliptic symbol depending only on the variable $\xi\in\mathcal{X}^*$, positive for $|\xi|$ large, and there exists $p_S\in S(M_{m,1+\epsilon},g)$ and $p_L\in S(M_{m-1,\epsilon},g)$ with $\partial_{x_j}p_L\in S(M_{m-1,1+\epsilon},g), 1\leq j\leq n$, such that $p=p_0+p_S+p_L$.
\end{enumerate}
Let $H$, $H_0$, respectively, the self-adjoint operators defined by $\mathfrak{Op}(p)$ and $\mathfrak{Op}(p_0)$ in $L^2(\mathcal X)$, both having domain $\boldsymbol{H}^{m}(\mathcal X)$. They have the following properties:
\begin{enumerate}
\item[a)] $\sigma_{\text{\sl ess}}(H)=\sigma_{\text{\sl ess}}(H_0)=\overline{p_0(\mathcal{X^*})}$.
\item[b)] The singular continuous spectrum of $H$ (if it exists) is contained in the set of critical values of $p_0$ defined as $\Lambda(p_0):=\{p_0(\xi)\;\mid\;p^\prime_0(\xi)=0\}$.
\item[c)] The eigenvalues of $H$ outside $\Lambda(p_0)$ have finite multiplicity and can accumulate only in $\Lambda(p_0)$.
\item[d)] {\sl (Limiting Absorption Principle)} The holomorphic function $\mathbb{C}_\pm\ni z\mapsto (H-z)^{-1}\in\mathcal{B}(\mathcal{G}^*_{1/2,1},\mathcal{G}_{-1/2,\infty})$ has a weak$^*$-continuous extension to\\ $\overline{\mathbb{C}_\pm}\setminus\left[\Lambda(p_0)\cup\sigma_{\text{\sl p}}(H)\right]$.
\item[e)] Properties a)-d) also hold if we replace $\mathfrak{Op}(p)$ with either $\mathfrak{Op}_A(p)$ or $\mathfrak{Op}^A(p)$.
\end{enumerate}
\end{theorem}

For the proof of this Theorem we shall need some auxiliary results.
\begin{lemma}\label{81}
Let Hypothesis \ref{I} be verified. We consider the symbol $p\in S^m(\Xi)$. 
\begin{enumerate}
\item[a)]
There exists $q\in S(M_{m-1,\epsilon},g)$ such that
$
\mathfrak{Op}^A(p)\,=\,\mathfrak{Op}(p)\,+\,\mathfrak{Op}(q).
$
\item[b)]
If, moreover, the symbol $p$ verifies for all $\alpha\in\mathbb{N}^n$  with $|\alpha|\geq1$, and for all $\beta\in\mathbb{N}^n$
\begin{equation}\label{81-b}
\left|(\partial^\alpha_x\partial^\beta_\xi p)(x,\xi)\right|\leq C_{\alpha,\beta}<x>^{-1}\ln(1+<x>)<\xi>^{m-|\beta|},
\end{equation}
then, for $\alpha$ and $\beta$ as above, we have
$\partial^\alpha_x\partial^\beta_\xi q\;\in\;S(M_\beta,g)$, with
$$
M_\beta(x,\xi):
=<x>^{-1-\epsilon}\ln(1+<x>)<\xi>^{m-1-|\beta|}.
$$
\end{enumerate}
\end{lemma}

\begin{proof}
a) For any $u\in\mathcal{S}(\mathcal{X})$ we have
$$
\mathfrak{Op}^A(p)u=
\int\!\!\!\int_{\mathfrak{X}\times\mathcal{X}^*}dy\dbar\xi\,
e^{i<x-y,\xi>}\,p\left(\frac{x+y}{2},\xi-\Gamma^A(x,y)\right)u(y).
$$
We can write
$
p\left(\frac{x+y}{2},\xi-\Gamma^A(x,y)\right)=p\left(\frac{x+y}{2},
\xi\right)+r(x,y,\xi)
$,
with
$
r(x,y,\xi):=-\left<\Gamma^A(x,y),\int_0^1d\tau\,(\partial_\xi p)\left(\frac{x+y}{2},\xi-\tau\Gamma^A(x,y)\right)\right>.
$
Thus $\mathfrak{Op}^A(p)=\mathfrak{Op}(p)+R$, where $R$ is defined by the integral kernel $K_r\in\mathcal{S}^*(\mathcal{X})$ given by the oscillatory integral
$
K_r(x,y)\,:=\,\int_{\mathcal{X}^*}\dbar\eta e^{i<x-y,\eta>}r(x,y,\eta).
$
We can write $R=\mathfrak{Op}(q)$, with
$$
q(x,\xi)=\int_\mathcal{X}dt\,e^{-i<t,\xi>}K_r\left(x+\frac{t}{2},x-
\frac{t}{2}\right)=
$$
$$
=\int\!\!\!\int_{\mathcal{X}\times\mathcal{X}^*}dt\,\dbar\eta
\,e^{i<t,\eta>}r\left(x+\frac{t}{2},x-\frac{t}{2},\eta+\xi\right).
$$
We set
$
\tilde{r}(x,z,\eta)\,:=\,r\left(x+\frac{z}{2},x-\frac{z}{2},
\eta\right)
$.
For any $\alpha\in\mathbb{N}^n$ we have $|(\partial^\alpha A)(x+sz)|\leq C_\alpha<x+sz>^{-\epsilon}$ and thus $|(\partial^\alpha A)(x+sz)|\leq C^\prime_\alpha<x>^{-\epsilon}<z>^{\epsilon}$ for any $x,z$ in $\mathcal{X}$ and any $s\in[-1/2,1/2]$. Also
$$
\left<\eta-\tau\int_{-1/2}^{1/2}dsA(x+sz)\right>^\lambda\leq C_{_\lambda}<\eta>^\lambda
$$
for any $(x,z)\in\mathcal{X}\times\mathcal{X}$, $\eta\in\mathcal{X}^*$, $\lambda\in\mathbb{R}$ and $\tau\in[0,1]$. From these inequalities and the formulae above, we obtain for any multiindices $\alpha,\beta,\gamma$ and for any $(x,z,\eta)\in\mathcal{X}\times\mathcal{X}
\times\mathcal{X}^*$
$$
\left|(\partial^\alpha_x\partial^\beta_z\partial^\gamma_\eta\tilde{r})(x,z,\eta)
\right|\,\leq\,C_{\alpha,\beta,\gamma}<x>^{-\epsilon}<z>^{\epsilon}
<\eta>^{m-1-|\gamma|}.
$$
For any two natural numbers $N_1$ and $N_2$ we get
$$
(\partial^\alpha_x\partial^\beta_\xi q)(x,\xi)\,=\,\int\!\!\!\int_{\mathcal X\times\mathcal{X}^*}dz\,\dbar\eta\,
e^{i<z,\eta>}<z>^{-2N_1}\times
$$
$$
\times(1-\Delta_\eta)^{N_1}<\eta>^{-2N_2}
(1-\Delta_z)^{N_2}(\partial^\alpha_x\partial^\beta_\eta\tilde{r})(x,z,\eta+\xi),
$$
so that we can deduce the estimation
$$
\left|(\partial^\alpha_x\partial^\beta_\xi q)(x,\xi) \right|\leq
$$
$$
\leq C_{\alpha,\beta}<x>^{-\epsilon}\int\!\!\!
\int_{\mathcal X\times\mathcal{X}^*}dz\,\dbar\eta\,<z>^{\epsilon-2N_1}<\eta>^{-2N_2}
<\eta+\xi>^{m-1-|\beta|}\,\leq
$$
$$
\leq
C^{\prime\prime}_{\alpha,\beta}<x>^{-\epsilon}
<\xi>^{m-1-|\beta|},
$$
by choosing $N_1$ and $N_2$ sufficiently large.

b) We can follow once again the proof of point (a), and remark that for $|\alpha|\geq1$ and $s\in[-1/2,1/2]$ we have
$$
\left|\partial^\alpha A(x+sz)\right|\leq
C^\prime_\alpha<x>^{-1-\epsilon}\ln(1+<x>)<z>^{1+\epsilon}\ln(1+<z>).
$$
The supplementary condition (\ref{81-b}) is needed when we estimate the derivatives $\partial^\alpha_x\tilde{r}$ (in the differentiation of $\partial_\eta p$ with respect to the first argument).
\end{proof}

\begin{lemma}\label{83}
Under Hypothesis \ref{I}, the operator $\mathfrak{Op}(q):\mathbf{H}^m(\mathcal{X})\rightarrow L^2(\mathcal{X})$ defined in Lemma \ref{81}(a) is compact.
\end{lemma}

\begin{proof}
For any $s\in\mathbb{R}$ the operator $<D>^s\equiv\mathfrak{Op}(p_s)\in\boldsymbol{\Psi}(p_s,g)$, and the operators $<D>^s$ and $<D>^{-s}$ are one the inverse of the other. If we denote $\tilde{p}_t(x,\xi):=<x>^t$, it follows that $\mathfrak{Op}(q)<D>^{-m}\in\boldsymbol{\Psi}(\tilde{p}_{-\epsilon}p_{-1},g)$, we also have
$\underset{|x|+|\xi|\rightarrow\infty}{\lim}\{<x>^{-\epsilon}<\xi>^{-1}\}=0$
and the Theorem 18.6.6 in \cite{Ho2} implies that $\mathfrak{Op}(q)<D>^{-m}$ is compact in $L^2(\mathcal{X})$. Thus $\mathfrak{Op}(q)=\left(\mathfrak{Op}(q)<D>^{-m}\right)<D>^m$ is compact as an operator from $\mathbf{H}^m(\mathcal{X})$ to $L^2(\mathcal{X})$.
\end{proof}

\begin{proposition}\label{84}
Under Hypothesis \ref{I}, if $p\in S^m(\Xi)$ is real and elliptic with $m>0$, it follows that $\mathfrak{Op}(q)$ (defined in Lemma \ref{81}(a)) is a relatively compact perturbation of $\mathfrak{Op}(p)$. In particular
$
\sigma_{\text{\sl ess}}[\mathfrak{Op}^A(p)]=\sigma_{\text{\sl ess}}[\mathfrak{Op}(p)].
$
Moreover, if $p(x,\xi)=p(\xi)$, then $\sigma_{\text{\sl ess}}[\mathfrak{Op}^A(p)]=\sigma_{\text{\sl ess}}[\mathfrak{Op}(p)]=\overline{p(\mathcal{X}^*)}$ and, if $\underset{|\xi|\rightarrow\infty}{\lim}p(\xi)=\infty$, then $\sigma_{\text{\sl ess}}[\mathfrak{Op}^A(p)]=\sigma_{\text{\sl ess}}[\mathfrak{Op}(p)]=[\gamma,\infty)$, with $\gamma:=\underset{\xi\in\mathcal{X}^*}{\inf}p(\xi)$.
\end{proposition}

\begin{proof}
Due to Theorem \ref{76} the operator $\mathfrak{Op}^A(p)$ is self-adjoint on the domain $\mathbf{H}^m(\mathcal{X})$, and the same is true for $\mathfrak{Op}(p)$. Lemma \ref{83} above implies that $\mathfrak{Op}(q)$ is a relatively compact perturbation of $\mathfrak{Op}(p)$. Then, if $p(x,\xi)=p(\xi)$ we see that $\mathfrak{Op}(p)$ is unitarily equivalent to the operator of multiplication with the function $p$ in $L^2(\mathcal{X}^*)$ and thus we have the second equality. For the last one just remark that $|p(\xi)\geq C|\xi|^m$ for $|\xi|\geq R$, so that $p$ will have constant sign for $|\xi|\geq R$.
\end{proof}

\begin{lemma}\label{87}
Let $m\in\mathbb{R}$ and $p\in S^m(\Xi)$ be given such that $p(x,\xi)=p(\xi)$; let also $\theta_r$ be a function in $C^\infty(\mathcal X)$ depending on the parameter $r\geq 1$ and such that for any $\alpha\in\mathbb{N}^n$ with $|\alpha|\geq 1$ satisfies
$
\left|(\partial_x^\alpha\theta_r)(x)\right|\leq\,C_\alpha\,r^{-1}$, $\forall x\in\mathbb{R}^n$, $\forall r\geq 1$.
Then $r(p\circ\theta_r-\theta_r\circ p)\in S^{m-1}(\mathbb{R}^n)$ uniformly in $r\geq 1$.
\end{lemma}

\begin{proof}
For any $u\in\mathcal{S}(\mathcal{X})$
$$
\left\{\left[\mathfrak{Op}(p),\mathfrak{Op}(\theta_r)\right]u\right\}(x)=
\int\!\!\!\int_{\mathcal{X}\times\mathcal{X}^*}dy\dbar\xi\; e^{i<x-y,\xi>}p(\xi)\left[\theta_r(y)-\theta_r(x)\right]u(y)\,=
$$
$$
=\,\sum_{j=1}^n\int\!\!\!\int_{\mathcal{X}\times\mathcal{X}^*}dy\dbar\xi\;
e^{i<x-y,\xi>}p_j(\xi)\lambda_j(x,y,r)u(y),
$$
where
$
p_j(\xi)=(D_jp)(\xi)\in S^{m-1}(\Xi)
$,
$
\lambda_j(x,y,r)=\int\limits_0^1ds\left(\partial_{x_j}\theta_r\right)(sx+(1-s)y)
$,
and for any $\alpha\in\mathbb{N}^n$, $\beta\in\mathbb{N}^n$, we have $r\partial_x^\alpha\partial_y^\beta\lambda_j\in L^\infty(\mathcal{X}\times\mathcal{X}\times[1,\infty))$.
Thus we have the formula $[\mathfrak{Op}(p),\mathfrak{Op}(\theta_r)]\,=\,r^{-1}\sum\limits_{j=1}^n\mathfrak{Op}(q_j)$ with 
$$
q_j(x,\xi,r):=\int\!\!\!\int_{\mathcal{X}\times\mathcal{X}^*} dy\dbar\eta\;e^{i<y,\eta>}p_j(\xi+\eta)\mu_j(x,y,r)
$$ 
and $\mu_j(x,y,r):=r\lambda_j(x+y/2,x-y/2,r)$. For $N_1$ and $N_2$ integers large enough 
$$
\left|\left(\partial_x^\alpha\partial_\xi^\beta q_j\right)(x,\xi,r)\right|\leq
$$
$$
\leq\,C_{\alpha\beta}<\xi>^{m-1-|\beta|}\int\!\!\!\int_{\mathcal X\times \mathcal{X}^*}dy\,
\dbar\eta\;<y>^{-2N_1}<\eta>^{-2N_2+|m-1-|\beta||}\leq
$$
$$
\leq C^\prime_{\alpha\beta}<\xi>^{m-1-|\beta|}.
$$
\end{proof}

\begin{lemma}\label{88}
Let $p\in S(M_{m,\delta},g)$ and $\delta>0$. Then there exists $q\in S^m(\Xi)$ such that $\mathfrak{Op}(p)=<Q>^{-\delta}\mathfrak{Op}(q)$.
\end{lemma}

\begin{proof}
The distribution kernel of the operator $<Q>^\delta\mathfrak{Op}(p)$ is given by
$
K(x,y):=\int_{\mathcal{X}^*}\dbar\eta\;e^{i<x-y,\eta>}<x>^\delta p\left(\frac{x+y}{2},\eta\right)
$.
We get $<Q>^\delta\mathfrak{Op}(p)=\mathfrak{Op}(q)$ with
$$q(x,\xi):
=\int\!\!\!\int_{\mathcal{X}\times\mathcal{X}^*}dy\dbar\eta\;e^{i<y,\eta-\xi>}
<x+y/2>^\delta p(x,\eta)=
$$
$$
=\int\!\!\!\int_{\mathcal X\times\mathcal{X}^*}dy\,\dbar\eta\;e^{i<y,\eta>}<x+y/2>^\delta p(x,\xi+\eta).
$$
For $N_1$ and $N_2$ large enough we get
$
\left|\left(\partial_x^\alpha\partial_\xi^\beta q\right)(x,\xi)\right|
\leq\,C_{\alpha\beta}<\xi>^{m-|\beta|}\int\!\!\!\int_{\mathcal X\times\mathcal{X}^*}dy\,
\dbar\eta\,<y>^{-2N_1+\delta}<\eta>^{-2N_2+|m-|\beta||}
\leq C^\prime_{\alpha\beta}<\xi>^{m-|\beta|}.
$
\end{proof}

\begin{proofTT} {\sl I. The case $H=\mathfrak{Op}(p)$}.
We are going to verify the hypothesis of Theorem 7.6.8. in \cite{ABG}, that directly implies the conclusion of our theorem; for the second equality in (a) we use Proposition \ref{84}.

\noindent
1. {\sl The difference $(H+i)^{-1}-(H_0+i)^{-1}$ is a compact operator in $L^2(\mathcal{X})$.}\\
We have 
$\mathfrak{Op}(p)=\mathfrak{Op}(p_0)+\mathfrak{Op}(p_S)+\mathfrak{Op}(p_L)$. Thus we can write
$$
(H+i)^{-1}-(H_0+i)^{-1}=-(H+i)^{-1}\{\mathfrak{Op}(p_S)+\mathfrak{Op}(p_L)\}
(H_0+i)^{-1}.
$$
From Lemma \ref{88} we see that $\mathfrak{Op}(p_S+p_L)=<Q>^{-\epsilon}\mathfrak{Op}(r)$ with $r\in S^m(\Xi)$. We remark that
$
(H_0+i)^{-1}$ is in $\mathcal{B}\left[L^2(\mathcal{X}),\boldsymbol{H}^m(\mathcal{X})
\right]$, $\mathfrak{Op}(r)$ is in $\mathcal{B}\left[\boldsymbol{H}^m(\mathcal{X}),L^2
(\mathcal{X})\right]$, $(H+i)^{-1}$ is in $\mathcal{B}\left[\boldsymbol{H}^{-m}(\mathcal{X}),L^2(\mathcal{X})\right]$
and $<Q>^{-\epsilon}$ is a compact operator in $\mathcal{B}(L^2(\mathcal{X}),\boldsymbol{H}^{-m}(\mathcal{X}))$, so that the difference of the resolvents is compact.

\noindent
2. {\sl For any $\rho\in C^\infty_0(\mathcal{X})$ with $\rho(x)=0$ for $x$ in a neighborhood of $0\in\mathcal{X}$, setting $\rho_r(x):=\rho(x/r)$, we have
$
\int_1^\infty dr\;\left\|\rho_r(Q)\mathfrak{Op}(p_S)\right\|_{_{\mathcal{B}(\mathcal{G},
\mathcal{G}^*)}}\;<\infty.
$}\\
By Lemma \ref{88} we have $\mathfrak{Op}(p_S)=<Q>^{-1-\epsilon}\mathfrak{Op}(q_0)$ with $q_0\in S^m(\Xi)$. Let us denote 
$
\widetilde{\rho_r}(x):=\rho\left(\frac{x}{r}\right)\left(\frac{r}{<x>}\right)^{1+\epsilon}.
$
Thus $\rho_r(Q)\mathfrak{Op}(p_S)=r^{-(1+\epsilon)}\widetilde{\rho_r}(Q)\mathfrak{Op}(q_0)$, and observing that
$$
\left\|\widetilde{\rho_r}(Q)\mathfrak{Op}(q_0)\right\|_{_{\mathcal{B}(\mathcal{G},
\mathcal{G}^*)}}\;=\;
\left\|<D>^{-m/2}\widetilde{\rho_r}(Q)\mathfrak{Op}(q_0)<D>^{-m/2}
\right\|_{_{\mathcal{B}(L^2(\mathcal{X}))}}
$$
it will be enough to prove that 
$$
\underset{r\geq 1}{\sup}\left\|<D>^{-m/2}\widetilde{\rho_r}(Q)\mathfrak{Op}(q_0)<D>^{-m/2}
\right\|_{_{\mathcal{B}(L^2(\mathcal{X}))}}\,<\,\infty.
$$
As the function $\widetilde{\rho_r}$ satisfies the hypothesis of Lemma \ref{87} we conclude that
$
r\left\{\mathfrak{p}_{-m/2}\circ\widetilde{\rho_r}-\widetilde{\rho_r}\circ
\mathfrak{p}_{-m/2}\right\}\;\in\;S^{-(1+m/2)}(\Xi),\;\text{uniformly in }
r\geq 1
$.
But 
$$
<D>^{-m/2}\widetilde{\rho_r}(Q)\mathfrak{Op}(q_0)<D>^{-m/2}=
$$
$$
=\widetilde{\rho_r}(Q)<D>^{-m/2}\mathfrak{Op}(q_0)<D>^{-m/2}+
$$
$$
+\left[<D>^{-m/2},\widetilde{\rho_r}(Q)\right]\mathfrak{Op}(q_0)<D>^{-m/2}
$$
and by the previous remark the second term above is a Weyl operator of order -1, uniformly for $r\geq 1$. The first term of the above sum is a Weyl operator of order 0, thus defining a bounded operator, uniformly in $r\geq 1$. 

\noindent
3. {\sl For any function $\theta\in C^\infty(\mathcal{X})$ with $\theta(x)=0$ on a neighborhood of $0\in\mathcal{X}$ and $\theta(x)=1$ in a neighbourhood of infinity, we have for any $j=1,\cdots,n$, 
$
\int_1^\infty \frac{dr}{r}\;\left\|\theta_r(Q)\left[Q_j,\mathfrak{Op}(p_L)\right]
\right\|_{_{\mathcal{B}(\mathcal{G},\mathcal{G}^*)}}\;<\infty$, for $\theta_r(x)=\theta(x/r)$}\\
In order to prove this estimation we notice that $\left[Q_j,\mathfrak{Op}(p_L)\right]=-\mathfrak{Op}(D_{\xi_j}p_L)$ and using Lemma \ref{88} we have $-\mathfrak{Op}(D_{\xi_j}p_L)=<Q>^{-\epsilon}\mathfrak{Op}(q_1)$ with $q_1\in S^{m-2}(\Xi)$. Thus it will be enough to prove that
$$
\underset{r\geq 1}{\sup}\,\left\|<D>^{-m/2}\varphi_r(Q)\mathfrak{Op}(q_1)<D>^{-m/2}\right\|_{_{\mathcal{B}(L^2(\mathcal{X}))}}\,<\,\infty,
$$
with $\varphi_r(x):=\theta(x/r)(r/<x>)^\epsilon$ and this follows by the same argument as in step (2) due to the fact that $\varphi_r$ also verifies the hypothesis of Lemma \ref{87}.

\noindent
4. {\sl For any test function $\theta\in C^\infty(\mathcal{X})$ with $\theta(x)=0$ on a neighborhood of $0\in\mathcal{X}$ and $\theta(x)=1$ in a neighbourhood of infinity, we have for any $j=1,\cdots,n$, 
$
\int_1^\infty \frac{dr}{r}\;\left\|\theta_r(Q)<Q>\left[D_j,\mathfrak{Op}(p_L)\right]
\right\|_{_{\mathcal{B}(\mathcal{G},\mathcal{G}^*)}}<\infty$, for $\theta_r(x)=\theta(x/r)$}\\
We start from the equality $\left[D_j,\mathfrak{Op}(p_L)\right]=\mathfrak{Op}(D_{x_j}p_L)$ and, using Lemma \ref{88}, we see that $\mathfrak{Op}(D_{x_j}p_L)=<Q>^{-(1+\epsilon)}\mathfrak{Op}(q_2)$ with $q_2\in S^{m-1}(\Xi)$. Thus everything goes on as before.

\noindent
{\sl II. The case $H=\mathfrak{Op}^A(p)$.} From Lemma \ref{81} we deduce the existence of a symbol $q\in S(M_{m-1,\epsilon},g)$ such that
for $0<\epsilon^\prime<\epsilon$ we have $\partial_{x_j}q\in S(M_{m-1,1+\epsilon^\prime},g)$, $1\leq j\leq n$ and $\mathfrak{Op}^A(p)=\mathfrak{Op}(p)+\mathfrak{Op}(q)=\mathfrak{Op}(p^\prime)$, where $p^\prime=p_0+p_S+p^\prime_L$ and $p^\prime_L=p_L+q\in S(M_{m-1,\epsilon},g)$, $\partial_{x_j}p^\prime\in S(M_{m-1,1+\epsilon^\prime},g)$, $1\leq j\leq n$. Since $p^\prime-p\in S(M_{m-1,\epsilon},g)$, we conclude that $p^\prime$ is elliptic. As $\mathfrak{Op}(p^\prime)=\mathfrak{Op}^A(p)$ is symmetric, we deduce that $p^\prime$ is real.

\noindent
{\sl III. The case $H=\mathfrak{Op}_A(p)$.} We start from the equality
$$
\left[\mathfrak{Op}(p\circ\nu^A)-\mathfrak{Op}(p)\right]u(x)=
\int\!\!\!\int_{\mathcal{X}\times\mathcal{X}^*}dy\dbar\xi\,e^{i<x-y,\xi>}r(x,y,\xi)u(y)
$$
for any $u\in\mathcal{S}(\mathcal{X})$, where $$r(x,y,\xi)=-\int_0^1d\tau\,<A((x+y)/2,(\partial_\xi p)((x+y)/2,\xi-\tau A(x+y)/2))>.$$ Repeating the proof of Lemma \ref{81} with $\Gamma^A(x,y)$ replaced by $A((x+y)/2)$, we conclude that there exists $q\in S(M_{m-1,\epsilon},g)$ such that 
$\partial_{x_j}q\in S(M_{m-1,\epsilon^\prime},g)$, $1\leq j\leq n$, $0<\epsilon^\prime<\epsilon$ and
$\mathfrak{Op}(p\circ\nu^A)=\mathfrak{Op}(p)+\mathfrak{Op}(q)$ and we are in the previous situation.
\end{proofTT}

For $t$ and $s$ in $\mathbb{R}$ let us denote by $\mathcal{H}^s_t$ the usual weighted Sobolev spaces, i.e. $\mathcal{H}^s_t=\{u\in\mathcal{S}^*(\mathcal{X})\mid\;<D>^s<Q>^tu\in L^2(\mathcal{X})\}$. We notice that $L^2_t\equiv\mathcal{H}^0_t$ and $\mathcal{G}=\mathcal{H}^{m/2}_0$. As shown in \cite{BGS}, for $\delta>0$ and $\gamma>\delta+1/2$, we have the continuous embedings $L^2_\gamma\subset\mathcal{H}^{-m/2}_{\delta+1/2}\subset\mathcal{G}^*_{1/2,1}$, the first being also compact. By duality we get the continuous embedings $\mathcal{G}_{-1/2,\infty}\subset\mathcal{H}^{m/2}_{-\delta-1/2}\subset L^2_{-\gamma}$, the last one being compact. One easily gets from the point (d) of our Theorem \ref{89} that the limiting absorption principle is valid in $\mathcal{B}(L^2_\gamma;L^2_{-\gamma})$ for the uniform topology, for any $\gamma>1/2$. 

\begin{remark}\label{7.8}
The limiting absorption principle is valid in $\mathcal{B}(\mathcal{H}^{-m/2}_\gamma;\mathcal{H}^{m/2}_{-\gamma})$ for the uniform topology, for any $\gamma>1/2$ (we may evidently suppose $\gamma\leq1$).
\end{remark}
To prove this fact we start with the following identity for $z\in\mathbb{C}_{\pm}$, consequence of the resolvent equation:
$$
(H-z)^{-1}=(H-i)^{-1}+(z-i)(H-i)^{-2}+(z-i)^2(H-i)^{-1}(H-z)^{-1}(H-i)^{-1}.
$$
As $(H-i)^{-1}\in\mathcal{B}(\mathcal{H}^{-m/2}_0;\mathcal{H}^{m/2}_0)$, the desired result will follow if we prove that $(H-i)^{-1}\in\mathcal{B}(\mathcal{H}^{-m/2}_\gamma;\mathcal{H}^{m/2}_\gamma)$ for any $\gamma\in[-1,1]$. In order to verify this relation, one may proceed as in the proof of our Lemma \ref{87}, and show that for any function $\varphi\in C^\infty(\mathcal{X})$ with $\partial^\alpha\varphi\in L^\infty(\mathcal{X})$ for $|\alpha|\geq1$, the commutator $[\mathfrak{Op}(p),\varphi(Q)]$ is a Weyl operator with symbol of class $S^{m-1}(\Xi)$. It follows that for any $u\in\mathcal{S}(\mathcal{X})$
$$
\varphi(Q)(H-i)^{-1}u=(H-i)^{-1}(\varphi u)+Tu,
$$
with $T\in\mathcal{B}(\mathcal{H}^{-m/2}_0;\mathcal{H}^{m/2}_0)$. This equality may then be extended to those elements of $\mathcal{H}^{-m/2}_0$ which verify $\varphi u\in\mathcal{H}^{-m/2}_0$. Choosing $\varphi(x):=<x>^\gamma$ with $0\leq\gamma\leq1$, we deduce that $(H-i)^{-1}\in\mathcal{B}(\mathcal{H}^{-m/2}_\gamma;\mathcal{H}^{m/2}_\gamma)$ for any $\gamma\in[0,1]$. By duality we obtain the same statement for $\gamma\in[-1,0]$.

\noindent
{\it Exemple 1 The magnetic relativistic Schr\"{o}dinger Hamiltonian $\mathfrak{Op}^A(<\xi>)$.} We consider the situation of Theorem \ref{89} with $p(\xi)=p_0(\xi)=<\xi>$, $p_S=p_L=0$. In this case we have $\overline{p_0(\mathcal{X}^*)}=[1,\infty)$ and $\Lambda(p_0)=\{1\}$. Thus this operator has no singular continuous spectrum.

\noindent
{\it Exemple 2 The operator $\mathfrak{Op}_A(p)\equiv\mathfrak{Op}(\mu^A)$.} Recall that $\mu^A(x,\xi)=<\xi-A(x)>$. We use again Theorem \ref{89}, with 
$p(\xi)=p_0(\xi)=<\xi>$. T. Umeda (\cite{Um}) has applied the Enss method to this operator obtaining properties (a), (b) and (c) from our Theorem \ref{89}, but not a limiting absorption principle. Besides, the hypothesis in \cite{Um} are less general then ours: he imposes restrictions on the vector potential $A$ of the form $|\partial^\alpha A_j(x)|\leq C_\alpha<x>^{-1-\epsilon}$ for any $\alpha\in\mathbb{N}^n$ (for $1\leq j\leq n$ and $x\in\mathcal{X}$) and $\epsilon>0$. We are making hypothesis only on the magnetic field $B$ and the only properties of $A$ that we use in the proof of Theorem \ref{89} are those deduced in Lemma \ref{80} from our Hypothesis \ref{I}.

\begin{remark}
An interesting result has been obtained by T. Ichinose and H. Tamura \cite{IT2}, showing that we have $\mathfrak{Op}(\mu^A)\geq 1$; some straightforward modifications of the proof in \cite{IT2} shows that $\mathfrak{Op}^A(<\xi>)\geq 1$. Thus under our hypothesis one has $\sigma(\mathfrak{Op}_A(<\xi>))=\sigma(\mathfrak{Op}^A(<\xi>))=[1,\infty)$.
\end{remark}

\noindent
{\it Exemple 3.} Our arguments may also be applied to the Schr\"{o}dinger operator $H=(D-A)^2$, taking $p(\xi)=p_0(\xi)=|\xi|^2$, $p_S=p_L=0$. In this case Theorem \ref{89} does not bring anything new (the situation may be understood from the one without magnetic field), this type of results being known for much more general (singular) magnetic filds, but also of the "short-range" type (see \cite{BMP}). This situation is a consquence of the fact that there exist magnetic fields which verify our Hypothesis \ref{I} with $\epsilon\leq0$ and such that $(D-A)^2$ has dense pure spectrum in an interval of $\mathbb{R}$ (see \cite{CFKS}), and thus Theorem \ref{89} clearly may not be applied.

{\bf Acknowledgements:} We thank the University of Geneva, where a large part
of  this work was done and especially 
Werner Amrein for his kind hospitality. We also acknowledge partial support
from the CERES Programme (contract no. 4-187/2004) and the CNCSIS
grant no. 13A (contract no. 27694/2005).

\end{document}